\documentclass[11pt]{article}
\makeatletter
\def\LaTeX{\leavevmode L\raise.42ex
\hbox{\kern-.3em\size{\sf@size}{0pt}\selectfont A}\kern-.15em\TeX}
\makeatother

\newcommand{\BibTeX}{{\rm B\kern-.05em{\sci\kern-.025emb}\kern-.08em\TeX}}

\usepackage{amssymb,amsfonts,amsmath,latexsym,graphicx,color}

\makeatletter

\def\theequation{\thesection.\@arabic\c@equation}
\makeatother

\newcommand{\lbl}[1]{\label{#1}}

\def\theequation{\arabic{section}.\arabic{equation}}

\newtheorem{theorem}{Theorem}[section]
\newtheorem{proposition}{Proposition}[section]
\newtheorem{lemma}{Lemma}[section]
\newtheorem{corollary}{Corollary}[section]
\newtheorem{remark}{Remark}[section]

\newcommand{\R}{{\mathbb R}}
\newcommand{\N}{{\mathbb N}}
\newcommand{\RR}{{\mathcal R}}
\newcommand{\equ}[1]{(\ref{#1})}
\newcommand{\ve}{\varepsilon}

\newcommand{\bb}{{\bf b}}
\newcommand{\LL}{{\mathcal L}}
\newcommand{\WW}{{\tt W}}
\newcommand{\XX}{{\mathcal X}}
\newcommand{\MM}{{\mathcal M}}

\newcommand{\BB}{{\mathcal B}}
\newcommand{\qed}{\framebox{\vspace{0.15truein}}}
\newcommand{\pd}[2]{\frac{\partial #1}{\partial #2}}
\newcommand{\pdd}[2]{\frac{\partial^2{#1}}{\partial{#2}^2}}

\newcommand{\la}{{\lambda}}
\newcommand{\ee}{{\bf e}}
\newcommand{\ddelta}{{\mathrm \delta}}

\newcommand{\ttt}{\tilde }
\newcommand{\LLL}{{\tt L}}
\newcommand{\KKK}{{\tt K}}
\newcommand{\TT}{{\cal T}}

\newcommand{\tX}{ {\tt X}  }
\newcommand{\tx}{ {\tt x}  }
\newcommand{\tz}{ {\tt z}  }
\newcommand{\tZ}{ {\tt Z}  }

\newcommand{\A}{\alpha }
\newcommand{\ff}{{\bf f}}

\renewcommand{\SS}{{\mathcal S}}
\renewcommand{\AA}{{\mathcal A}}
\newcommand{\proof}[1]{\noindent{\bf Proof{#1}.}}

\definecolor{green2}{rgb}{0,0.6,0}

\newcommand\runninghead[2]{\pagestyle{myheadings}\markboth{{\footnotesize\it{\quad #1}\hfill}}{{\footnotesize\it{#2\hfill\quad}}}}\headsep=40pt
\runninghead{M. del Pino, M. Kowalczyk, F.Pacard, J.
Wei\qquad}{\qquad \qquad Multiple-end solutions of autonomous elliptic problems}

\begin{document}
\title{The Toda system and multiple-end  solutions of autonomous planar elliptic problems}

\author{Manuel del Pino
\footnote{{Manuel del Pino - Departamento de
Ingenier\'{\i}a  Matem\'atica and CMM, Universidad de Chile,
Casilla 170 Correo 3, Santiago,
Chile.} Email: {delpino@dim.uchile.cl}}
\medskip
Micha{\l} Kowalczyk
\footnote{{Micha{\l} Kowalczyk - Departamento de
Ingenier\'{\i}a  Matem\'atica and CMM, Universidad de Chile,
Casilla 170 Correo 3, Santiago, Chile.}
Email {kowalczy@dim.uchile.cl}}
\medskip
Frank Pacard
\footnote{{Franck Pacard -
Universit\'e Paris 12 and Institut Universitaire de France}
Email {pacard@univ-paris12.fr}}
\medskip
Juncheng Wei
\footnote{{Juncheng  Wei - Department of Mathematics, Chinese University of Hong Kong, Shatin, Hong Kong
} Email{wei@math.cuhk.edu.hk}}}


\date{} \maketitle

\abstract{We construct a new class of positive solutions for the classical elliptic problem $$ \Delta u  -u + u^p   = 0, p\geq 2, \quad\mbox{in}\ \R^2. $$
We establish a deep  relation between them and the following Toda
system
\[c^2 f_j'' =   e^{f_{j-1}- f_{j}} - e^{f_{j}-
f_{j+1}}\quad\hbox{in }\R,\quad  j=1,\ldots, k.\]
We show that  these solutions have the approximate form  $ u(x,z)
\sim \sum_{j=1}^k w(x- f_{j}(z))$ where $w$ is the unique even,
positive, asymptotically vanishing solution of $w'' -w + w^p =0$ in $\R$.
Functions
 $f_j(z)$, representing the multiple ends of $u(x,z)$,  solve the  aforementioned  Toda  system, they are  even, asymptotically linear, with
$$f_0 \equiv -\infty < f_1\ll\cdots\ll f_k < f_{k+1} \equiv +\infty.$$

The solutions of the elliptic problem we construct  have their counterpart in the theory of constant mean curvature
surfaces. An analogy can also be made between their construction  and the  gluing of constant scalar curvature Fowler
singular metrics in the sphere.}

\bigskip
\bigskip

\setcounter{section}{0}



\tableofcontents

\bigskip
\bigskip

\setcounter{equation}{0}
\section{Introduction and statement of main results}
 This paper deals with the classical semilinear elliptic
problem
\begin{equation}
\Delta u  -u + u^p   = 0,\quad u>0, \quad\hbox{in } \R^N
\label{sch}\end{equation} where $p>1$. Equation \equ{sch}  arises
for instance as the standing-wave problem for the standard
nonlinear Schr\"odinger equation
$$
i\psi_t = \Delta_y \psi  + |\psi|^{p-1}\psi,
$$
typically $p=3$, corresponding to that of solutions of the form
$\psi(y,t) = u(y)e^{-it}$. It also arises in  nonlinear models in
Turing's theory biological theory of pattern formation \cite{turing}  such as the
Gray-Scott or Gierer-Meinhardt systems, \cite{grayscott,gm}.
The solutions of \equ{sch} which decay to zero at infinity are well understood. Problem \equ{sch} has a  radially symmetric
solution $w_N(y)$ which approaches $0$ at infinity provided that
$$
1< p < \begin{cases}  \frac{N+2}{N-2} & \hbox{ if }\  N\ge 3,\medskip\\   +\infty \ \  & \hbox{ if }\ N=1,2,\\
\end{cases}
$$
see
 \cite{strauss,berestyckilions}. This solution
is unique \cite{kwong}, and actually any
positive solution to \equ{sch} which vanishes at infinity must be
radially symmetric around some point\cite{GNN}.

\medskip
Problem  \equ{sch} and its variations  have been broadly treated in
the PDE literature in the last two decades. These variations are mostly of one of the two types:
(\ref{sch}) is changed to  a non-autonomous problem with a potential depending on the space variable;  or (\ref{sch}) is considered in a bounded domain under suitable boundary conditions. Typically, in both versions
 a small parameter is introduced  rendering (\ref{sch})
 a singular  perturbation problem. We refer the reader to the works
\cite{abc,bw,cotirabinowitz,df1,df2,dkw1,fw,gw,kw,linniwei,m1,m2,ni,nt1,nt2} and
references therein. Many constructions in the literature refer to
``multi-bump solutions'', built by perturbation of a sum of copies
of the basic radial bump $w_N$ suitably scaled, with centers
adjusted in equilibrium under appropriate constraints on the
potential or the geometry of the underlying domain.

\medskip  Much less is known about solutions to this
equation in entire space which do not vanish at infinity (while
they are all known to be bounded, see \cite{quittner}). For example, the solution $w_N$ of \equ{sch} in $\R^N$
induces a solution in $\R^{N+1}$ which only depends on $N$ variables. This solution vanishes
asymptotically  in all but one variable. For simplicity, we
restrict ourselves to the case $N=2$, and consider positive
solutions $u(x,z)$ to problem \equ{sch} which are even in $z$ and
vanish as $|x|\to +\infty$, namely
\begin{equation}
u(x,z) = u(x,-z) \quad\hbox{for all } (x,z) \in \R^2,
\label{even}\end{equation} and
\begin{equation}
\lim_{|x|\to +\infty } u(x,z) = 0 \quad\hbox{for all } z \in \R.
\label{vanish}\end{equation}

A canonical example is thus built from the one-dimensional  bump $w_1$, which we denote in the sequel just by $w$, namely the unique solution of
the ODE
\begin{align}
w'' - w + w^p =0, &\quad w>0,\quad\hbox{ in } \R,
\label{ode11}\\
w'(0)=0,&\quad w(x) \to 0 \quad\hbox{ as } |x| \to +\infty ,
\label{ode12}\end{align}
 corresponding in phase plane to a homoclinic orbit for the equilibrium 0.
Using this  function  we  can define  a family of solutions $u$ of
equation \equ{sch}  with the properties \equ{even}-\equ{vanish} setting $u(x,z) := w(x-a)$, $a\in \R$. By analogy with the
above terminology, we may call  these solutions "single
bump-lines". A natural question is whether a solution that
satisfies \equ{even}-\equ{vanish} and which is in addition even in $x$
must equal $w(x)$.  The solution $w$  of \equ{sch} was found to be
isolated by Busca and Felmer in
 \cite{buscafelmer} in a uniform topology which avoids oscillations at infinity. On the other hand, a second class of solutions which are even both in $z$ and $x$  was discovered by Dancer in
\cite{dancer} via local bifurcation arguments. They constitute a
one-parameter family of solutions which are periodic in the $z$
variable and originate from   $w(x)$. Let us briefly review their
construction: we consider problem \equ{sch} with $T$-periodic
conditions in $z$,
\begin{equation}
u(x,z+T) = u(x,z) \quad\hbox{for all } (x,z) \in \R^2,
\label{per}\end{equation}
and regard $T>0$ as a bifurcation parameter.
The linearized operator around the single bump line
is
$$
L(\phi) = \phi_{zz} + \phi_{xx} +(pw^{p-1} -1)\phi .
$$
It is well known that the eigenvalue problem
\begin{equation}\phi_{xx} +(pw^{p-1} -1)\phi = \la \phi, \label{eigen}\end{equation}
has a unique positive eigenvalue  $\la_1$ with $Z(x)$ a positive eigenfunction.
We observe that the operator $L$ has a bounded element of its kernel given by
$$
Z(x) \cos ( \sqrt{\la_1} z) ,
$$
which turns out to be the only one which is  even, both in $x$ and
$z$ variable, and  in addition $T= \frac{2\pi}{\sqrt{\la_1}}$-periodic
in $z$. Crandall-Rabinowitz bifurcation theorem can then be
adapted to yield existence of a continuum of solutions bifurcating  at this
value of $T$, periodic in $z$ with period $T_\delta = \frac{2\pi}{\sqrt{\la_1}} + O(\delta)$. They are uniformly close to $w(x)$ and  their  asymptotic formula is:
$$
w_\delta (x,z) = w(x) + \delta Z(x) \cos ( \sqrt{\la_1} z) +
O(\delta^2) e^{-|x|}.
$$
We refer to the functions $w_\delta$ in what follows as {\em Dancer solutions}.

\bigskip
The purpose of this paper is to construct a new type of solutions of (\ref{sch}) in $\R^2$ that have  multiple ends in the form of  {multiple bump-lines}, and  satisfy in addition  properties \equ{even}-\equ{vanish}. To explain let us consider the function
 $$
w_*(x) = \sum_{j=0}^k w(x-a_j).
$$
While for numbers $a_1 \ll a_2 \ll \cdots \ll a_k$, $w_*$ appears to be a very good approximation to a solution of the ODE \equ{ode11}, the only solutions
of \equ{ode11} which go to zero at infinity must be single translates of $w(x)$.
 Our main
result in this paper is that, in the two dimensions,  there exist  solutions $u(x,z)$ with a profile of this type for
each fixed $z$. Of course in this case the numbers $a_j$
must be replaced by non-constant functions   $f_j(z)$.  In
addition, it will be necessary to use as basic cells not just the
standard bump-lines, but rather the {\em wiggling bump lines} found by
Dancer.

 Thus, what we actually look for is a solution $u(x,z)$
which is close, up to lower order terms, to a multi bump-line of
the form
\begin{equation}
w_*(x,z) = \sum_{j=0}^k w_{\delta_j} (x-f_j(z),z ),
\label{form}\end{equation} for suitable small numbers $\delta_j$
and  even functions
$$f_1(z) \ll f_2(z) \ll \cdots \ll f_k(z), $$
which have uniformly small derivatives.  The functions $f_j$ cannot be arbitrary and they turn out to
satisfy (asymptotically) a second order
system of differential equations, the {\em Toda system}, given by
\begin{equation}
c_p^2 f_j'' =   e^{f_{j-1}- f_{j}} - e^{f_{j}-
f_{j+1}}\quad\hbox{in }\R,\quad  j=1,\ldots, k,
\label{toda}\end{equation} with the conventions $f_0 = -\infty$,
$f_{k+1} =+\infty$, where $c_p$ is an explicit positive constant
that will be specified later.  In agreement with the symmetry
requirement \equ{even} we consider even solutions of system
\equ{toda}, namely
\begin{equation} f_j(z) = f_j(-z)\quad \hbox{for all } z\in \R.
\label{even1}\end{equation} We observe  that for an even
solution $\ff=(f_1,\dots,f_k)$ of this system, function  $\ff_\alpha$ defined by
\begin{equation}
{\bf f}_\alpha=(f_{\alpha1},\dots,f_{\alpha k}), \quad f_{\alpha j} (z) \, := \, f_{1j}(\alpha z )\, + \, (j- \frac {k+1}
2)\log\frac 1\alpha \label {falpha}\end{equation} is also an even
solution of the system.
Our main result in this paper asserts that given such an $\ff$
whose values at $z=0$ are ordered, for all sufficiently small $\alpha$ there exists a multi
bump-line solution of
approximate form \equ{form} with  $\ff$ replaced with  $ \ff_\alpha$
 and suitably chosen
small numbers $\delta_j$ dependent on $\alpha$. Thus, we fix
numbers
\begin{equation}
 a_1 <a_2 < \cdots < a_k,\qquad \sum_{i=1}^k a_i=0,
 \label{aj}\end{equation}
 (as we will see shortly the latter condition can be assumed without loss of   generality)
and consider the unique solution $\ff$ of system \equ{toda} for
which
\begin{equation}
f_{j}(0) = a_j,\quad f_{j}'(0) = 0,\quad j=1,\ldots, k ,
 \label{f1}\end{equation}
and their associated scalings $\ff_{\A}$. As an explicit example,
we immediately check that for $k=2$ we have
\begin{align}
\begin{aligned}
f_{1}(z) &=- \frac{1}{2}\Big(\log(2\lambda^{-2}c_p) - \log \frac {1}{2 \cosh^2(\lambda z/2)}\Big) ,\\
 f_{2}(z)& =
\frac{1}{2}\Big( \log(2\lambda^{-2}c_p)- \log \frac {1}{2 \cosh^2(\lambda z/2)}\Big), \label{k2}\end{aligned}\end{align}
where
$$
\lambda=\sqrt{\frac{2c_p}{e^{a_2-a_1}}}.
$$
As we
will see later, in the case of a general $k$ the functions $f_{j}$
are asymptotically linear: the limits $\nu_j= f'_{1j}(+\infty),$
exist and
\begin{equation} \nu_1<\nu_2<\cdots < \nu_k ,\quad  \sum_{j=0}^k \nu_j = 0. \label{nu}\end{equation}
For instance, for $k=2$ we have  $\nu_2 = -\nu_1 =  \frac {1} {2}\sqrt{\frac{2c_p}{e^{a_2-a_1}}}$ according to formula \equ{k2}. Besides, we
have (globally) for $\A$ small
$$ f_{\A 1}(z) \ll f_{\A 2} (z)\ll \cdots \ll  f_{\A k}(z) ,\qquad
f_{\A j}' (+\infty )= \nu_j\A,$$ and
$$
f_{\A j} (z) \, = \, \nu_j\A\, |z| + b_j \, + \, (j- \frac {k+1}
2)\log\frac 1\alpha  + O( e^{-\vartheta\A |z|} ),\quad\hbox{as } |z| \to
+\infty,
$$
for certain scalars $b_j$ and $\vartheta>0$.  These are standard  facts about the Toda system
that  can be found for instance  in \cite{konst}.
Thus, each of the multiple ends of $u(x,z)$  is a  bump-line
that is  nearly straight but bent, with an angle slightly distinct than the angles of the other ends.
The Toda system is a classical  model
 describing scattering of $k$  particles distributed on a straight  line, which  interact
only with their closest neighbors with a forces given by a potential depending on the
exponentials of their mutual distances. Here the $z$
variable is interpreted as time. In this context, $\ff_\alpha$
corresponds to a setting in which the particles starting  from the rest,
scatter at slightly different, nearly
constant small velocities whose average is zero. The latter fact
follows the identity $\sum_{j=1}^k f_{\A j}''(z) = 0$, which also
implies conservation of the center of mass used in (\ref{aj}).

\medskip

\medskip
Our main result is:



\begin{theorem}\label{teo1}
Assume that $N=2$ and $p \geq 2$. Given $k\ge 2$ and numbers $a_j$ as
in $\equ{aj}$,  for any sufficiently small number
$\alpha>0$, there exists a solution $u_\alpha $ of equation
$\equ{sch}$ which satisfies conditions $\equ{even}$-$\equ{vanish}$, and
that has the form
$$
u_\A(x,z) = \sum_{j=1}^k w_{\delta_j} (x- f_{\alpha j} (z),z )\,
(1 + o(1)).
$$
Here $\ff_\alpha$ is the scaling $\equ{falpha}$ of $\ff$, the
unique solution of $\equ{toda}$ satisfying $\equ{f1}$, and
$o(1)\to 0$, $\delta_j \to 0$ as $\alpha\to 0$.

\end{theorem}

\begin{remark}{\em
By no means we have intended to state Theorem \ref{teo1} in its most general form. For instance, the same result holds for
the more general problem
\begin{equation}
\Delta u  + g(u)  = 0,\quad u>0, \quad\hbox{in } \R^2  .
\label{sch1}\end{equation}
where $g$ is of class $C^2$, $g(0)=0$, $g'(0)<0$, $\int_0^c g(s)ds = 0$ for some $c>0$. In this case a homoclinic solution
$w$ analogous to that of \equ{ode11}-\equ{ode12} exists.
On the other hand, for problem \equ{sch} in $\R^N$, $N\ge 3$, we may take as the basis of the construction
the unique radially symmetric decaying solution $w_{N-1}(x)$ of \equ{sch} in $\R^{N-1}$,  provided that $p< \frac{N+1}{N-3}$, for which Dancer solutions $w_{N-1,\delta}(x)$ are equally available. In such a case we look for solutions close to a function of the form
$$
\sum_{j=1}^k w_{N-1,\delta_j}(x_1,x_2,\ldots,x_{N-2}, x_{N-1} -f_{\A j}(z) ,z),
$$
which are radial in the first $N-2$ variables. In both extensions
the necessary changes in the constructions are  straightforward,
so that for simplicity we only consider here the
case $N=2$ for \equ{sch}.
 }
\end{remark}

\begin{remark}{\em
The structure of equation \equ{sch} makes it natural to look for solutions
$u$ which are even both in $x$ and $z$ variables. It turns out that
if the numbers $a_j$ in \equ{aj} satisfy the symmetry requirement
$$ a_j = - a_{k+1 -j}, \quad j=1,\ldots, k, $$ then the solutions $\ff_\A$ satisfy similarly
 \begin{equation}
f_{\A j} = -f_{\A k+1 -j}, \quad \hbox{for all $j=1,\ldots, k$}. \label{simf}\end{equation}
Indeed, in such a case the construction of the solutions $u_\A$ in Theorem \ref{teo1}
yields that they are even in both variables:
\begin{equation}
u_\A (-x,z) = u_\A(x,z) = u_\A (x,-z), \quad\hbox{for all } (x,z) \in \R^2,
\label{even3}\end{equation}
 In particular for  $k$ even this solution satisfies
\begin{equation}
\lim_{x\to \pm\infty }  u_\alpha (x,z ) = 0, \quad\hbox{for all } \,
z\in \R,\qquad
\lim_{z\to \pm\infty }  u_\alpha (x,z ) = 0, \quad\hbox{for all } \,
x\in \R.
\label{shno}\end{equation} By the well known result by
Gidas, Ni and Nirenberg \cite{GNN} a positive solution of equation
\equ{sch} that satisfies \equ{even3} and the limit conditions \equ{shno} {\em uniformly} must be radially symmetric around
the origin. Theorem \ref{teo1} shows that uniformity cannot be relaxed in this
classical result.

}
\end{remark}

\bigskip
One of the striking features of the existence  result in  Theorem~\ref{teo1}, which is a purely PDE result,   is that its counterparts can be found in geometric framework. Indeed,  there are many examples  where correspondence between  solutions of (\ref{sch}) and those of some geometric problem can be drawn. To illustrate this, we will concentrate on  what is perhaps the most spectacular  one: the analogy between the theory of complete constant mean curvature surfaces in Euclidean $3$-space and the study of entire solutions of (\ref{sch}). For simplicity we will  restrict ourselves to  constant mean curvature surfaces in $\mathbb R^3$ which have embedded coplanar ends. In the following we will draw parallels
between these geometric objects and families of solutions of (\ref{sch}).

\medskip

Embedded constant mean curvature surfaces of revolution were found by Delaunay in the mid 19th century \cite{Del}.  They constitute a smooth one-parameter family of singly periodic surfaces $D_\tau$, for $\tau \in (0,1]$, which interpolate between the cylinder $D_1 = S^1(1) \times \mathbb R$ and the singular surface $ D_0 : = \lim_{\tau\rightarrow 0} D_\tau$, which is the union of an infinitely many spheres of radius $1/2$ centered at each of the points $(0,0,n)$ as $n \in \mathbb Z$. The Delaunay surface $D_\tau$ can be parametrized by
\[
X_\tau (x, z) = (\varphi (z)  \, \cos x , \varphi (z)  \, \sin x  , \psi (z) ) \in D_\tau \subset \mathbb R^{3} ,
\]
for $ (x, z) \in {\mathbb R} \times \mathbb R / 2 \pi \mathbb Z$. Here the function $\varphi$ is  smooth solution of
\[
(\varphi')^2 + \left( \frac{\varphi^2 + \tau}{2} \right)^2 =  \varphi^2,
\]
and the function $\psi$ is defined by
\[
\psi' = \frac{\varphi^2 + \tau}{2} .
\]

As already mentioned, when $\tau =1$, the Delaunay surface is nothing but a right circular cylinder $D_1 = S^1 (1) \times \mathbb R$, with the unit circle as the cross section. This cylinder is clearly invariant under the continuous group of vertical translations, in the same way that the single
bump-line solution of (\ref{sch})  is  invariant under a one parameter group of translations. It is then natural to agree on the correspondence between
\[
\fbox{
$\begin{array}{c}
\mbox{The cylinder} \\
D_1 = S^1\times {\mathbb R}
\end{array}
$} \longleftrightarrow
\fbox{
$\begin{array}{c}
\mbox{The single bump-line}\\
(x,z) \longmapsto w (x)
\end{array}$}.
\]
Let us denote by $w_{2}$ the unique radially symmetric, decaying solution of (\ref{sch}). Inspection of the  other end of the Delaunay family, namely when the parameter $\tau$ tends to $0$, suggests the correspondence between
\[
\fbox{
$\begin{array}{c}
\mbox{The sphere} \\
S^1(1/2)
\end{array}
$} \longleftrightarrow
\fbox{
$\begin{array}{c}
\mbox{The radially symmetric solution} \\
(x,z) \longmapsto w_2 (\sqrt{x^2+ z^2})
\end{array}$}.
\]
To justify this correspondence, let us observe that on the one hand, as the parameter $\tau$ tends to $0$,  the surfaces $D_\tau$ resemble a sequence of spheres of radius $1/2$ arranged along the $x_3$-axis which  are connected together by small catenoidal necks. On the other hand an analogous
 solution of (\ref{sch})
   can be built as follows. Let $S_R = \R \times (0,R)$ and consider a least energy (mountain pass) solution
 in $H^1(S_R)$ for the  the energy
 $$ \frac 12 \int_{S_R} |\nabla u|^2 + \frac 12 \int_{S_R} u^2 -  \frac 1{p+1} \int_{S_R} u^{p+1} , $$
 for large $R>0$, which we may assume to be even in $x$ and with maximum located at the origin. For $R$ very large, this solution, which satisfies
 zero Neumann boundary conditions, resembles half of the unique radial, decaying  solution $w_{2}$ of \equ{sch}. Extension by successive even reflections in $z$ variable yields a solution to \equ{sch} which resembles a periodic array of radially symmetric solutions of (\ref{sch}), with a very large period, along the $z$-axis. While this is not known, these solutions may be understood as a limit of the branch solutions constructed by Dancer.

\medskip
More generally, there is a natural correspondence between
\[
\fbox{
$\begin{array}{c}
\mbox{Delaunay surfaces} \\
D_\tau
\end{array}
$} \longleftrightarrow
\fbox{
$\begin{array}{c}
\mbox{Dancer solutions}\\
(x,z) \longmapsto w_\delta (x,z)
\end{array}$}.
\]
To give further credit to this correspondence, let us recall that the Jacobi operator about the cylinder $D_1$ corresponds to the linearized mean curvature operator when nearby surfaces are considered as normal graphs over $D_1$. In the above parameterization, the Jacobi operator reads $J_1 = \frac{1}{\varphi^2} \, \left(  \partial_x^2 + \partial_z^2 + 1 \right)$.  In this geometric context, it plays the role of the linear operator defined in section 2 which is the linearization of (\ref{sch}) about the single bump-line solution $w$. Hence we have the correspondence
\[
\fbox{
$\begin{array}{c}
\mbox{The Jacobi operator} \\
J_1 =\frac{1}{\varphi^2} \, \left(  \partial_x^2 + \partial_z^2 + 1 \right)
\end{array}
$} \longleftrightarrow
\fbox{
$\begin{array}{c}
\mbox{The linearized operator}\\
L = \partial_x^2 + \partial_z^2 -1 + p \, w^{p-1}
\end{array}$}.
\]
In our construction, the polynomially bounded kernel of the linearized operator $L$ plays a crucial role. Similarly, the polynomially bounded kernel of the Jacobi operator $J_1$ has some geometric interpretation. Let us recall that we only consider surfaces whose ends are coplanar, the Jacobi fields  associated to the action of rigid motions are then given by
\[
(x,z) \longmapsto \cos x \qquad \mbox{and} \qquad (x,z) \longmapsto z \, \cos x,
\]
which correspond respectively to the action of translation and the action of the rotation of the axis of the Delaunay surface $D_1$.  Clearly, these Jacobi fields are the counterpart of  the elements of the kernel of $L$ which are given by
\[
(x,z) \longmapsto \partial_x w (x) \qquad \mbox{and} \qquad (x,z) \longmapsto z \, \partial_x w (x),
\]
since the latter are also generated using the invariance of the problem with respect to the same kind of rigid motions.

\medskip
Two additional Jacobi fields  associated to $J_1$ are given by
\[
(x,z) \longmapsto  \cos z \qquad \mbox{and} \qquad (x,z) \longmapsto \sin z,
\]
which are associated to the existence of the family $D_\tau$ as $\tau$ is close to $1$, as can be easily seen using a bifurcation analysis, in a similar  way that the functions
\[
(x,z) \longmapsto  Z(x) \, \cos (\sqrt{\lambda_1} z) \qquad \mbox{and} \qquad (x,z) \longmapsto  Z(x) \, \sin (\sqrt{\lambda_1} z),
\]
are associated to the existence of Dancer solutions when the parameter $\delta$ is close to $0$. These two bifurcation results have their origin in the fact that we have the correspondence between
\[
\fbox{
$\begin{array}{c}
\mbox{The ground state $1$ of } \\
\partial_x^2  + 1
\end{array}
$} \longleftrightarrow
\fbox{
$\begin{array}{c}
\mbox{The first eigenfunction  $Z(x)$ of }\\
\partial_x^2 -1 + p \, w^{p-1}
\end{array}$},
\]
both of them associated to negative eigenvalues.  The fact that the least eigenvalue of these operators is negative is precisely the reason why a bifurcation analysis can be performed and gives rise to the existence of Delaunay surfaces close to $D_1$ or Dancer's solutions close to the bump-line $w$.

\medskip

With these analogies in mind, we can now {\em translate} our main result into the constant mean curvature surface framework. The result of Theorem~\ref{teo1} corresponds to the connected sum of finitely many copies of the cylinder $S^1(1) \times \mathbb R$ which have a common plane of symmetry. The connected sum construction is performed by inserting small catenoidal necks between two consecutive cylinders and this can be done in such a way that the ends of the resulting surface are coplanar.  Such a result, in the context of constant mean curvature surfaces, follows at once from \cite{Maz-Pac-Pol}. It is observed that, once the connected sum is performed the ends of the cylinder have to be slightly bent and moreover,  the ends cannot be kept asymptotic to the ends of right cylinders but have to be asymptotic to Delaunay ends with parameters close to $1$, in agreement with the result of Theorem~\ref{teo1}.

\medskip

However there is a major difference.  The Toda system which governs the level sets has found no analogy in the constant mean curvature surfaces. This is mainly due to the strong interactions in the elliptic equations.

\medskip

Another (older) construction of complete noncompact constant mean curvature surfaces was performed by N. Kapouleas  \cite{Kap} (see also \cite{Maz-Pac-2}) starting with finitely many halves of Delaunay surfaces with parameter $\tau$ close to $0$ which are connected to a central sphere.  The corresponding solutions of (\ref{sch}) have recently been constructed by A. Malchiodi in \cite{malchiodi}. It should be clear that many more examples of solutions of (\ref{sch}) can be found using the above correspondence and we shall return to this in a forthcoming paper.

\medskip

It is well known that the story of complete constant mean curvature surfaces in $\mathbb R^3$  parallels that of  complete locally conformally flat metrics with constant, positive scalar curvature. Therefore, it is not surprising that there should be a correspondence between these objects in conformal geometry and solutions of (\ref{sch}).  For example, Delaunay surfaces and Dancer solutions should now be replaced by Fowler solutions which correspond to constant scalar curvature metrics on the cylinder $\mathbb R \times S^{n-1}$ which are conformal to the product metric $dz^2 + g_{S^{n-1}}$, when $n \geq 3$. These are given by
\[
v^{\frac{4}{n-2}} ( dz^2 + g_{S^{n-1}}),
\]
where $z \longmapsto v(z)$ is a smooth positive solution of
\[
(v')^2 -  v^2 + \tfrac{n-2}{n} \,  v^{\frac{2n}{n-2}} = - \tfrac{2}{n}  \, \tau^2.
\]
When $\tau =1$ and $v \equiv 1$ the solution is a straight
cylinder while as $\tau$ tends to $0$ the metrics converge on
compacts to the round metric on the unit sphere.  The connected
sum construction for such Fowler type metrics was performed by R.
Mazzeo, D. Pollack and K. Uhlenberk \cite{Maz-Pol-Uhl} (where it
is called the dipole construction).  N. Kapouleas' construction
mentioned above was initially performed by R. Schoen \cite{Sch}
(see also R. Mazzeo and F. Pacard \cite{Maz-Pac-2}).

 \setcounter{equation}{0}

    \section{Linear theory}
    In this section we will consider the basic linearized operator. The developments presented here are crucial for our paper later on.
    By $w$ we will denote the homoclinic solution to $u''-u+u^p=0$ such that $w'(0)=0$. Let
    $$
    L_0(\phi)=\phi_{xx}+f'(w)\phi, \quad f'(w)=pw^{p-1}-1.
    $$
    We recall some well known facts about $L_0$. First notice that $L_0(w_x)=0$ i.e. has  one dimensional kernel. Second we observe that
    $$
    \lambda_1=\frac{1}{4}(p-1)(p+3), \qquad Z=\frac{w^{(p+1)/2}}{\sqrt{\int_{\R} w^{p+1}}},
    $$
    correspond, respectively,  to principal eigenvalue and eigenfunction. Except for $\lambda_1$ and $\lambda_2=0$
    the rest of the spectrum of $L_0$ is negative. As a consequence of these facts we observe that problem
 \begin{equation}
    L_0(\phi)-\xi^2\phi=h,
  \label{eqxi}
  \end{equation}
    is uniquely solvable whenever $\xi\neq \pm\sqrt{\lambda_1}, 0$ for $h\in L^2(\R)$.  Actually, rather standard argument,
    using comparison principle  and the fact that  $L_0$ is of the form
    $$
    L_0(\phi)=\phi_{xx}-\phi+q(x)\phi,\quad |q(x)|\leq Ce^{-c|x|},
    $$
    can be used to show that whenever $h$ is for instance a compactly supported function then the solution of \equ{eqxi}
    is an exponentially decaying function.

    Let us  consider operator
    $$
    L\phi=L_0(\phi)+\phi_{zz},
    $$
    defined in the whole plane $(x,z)\in \R^2$.
    Equation $L(\phi)=0$,  has $3$ obvious bounded solutions
    $$
    w_x(x),\quad  Z(x)\cos{\sqrt{\lambda_1}z}, \quad Z(x)\sin(\sqrt{\lambda_1}z).
    $$
Our first result shows that converse is also true.
 \begin{lemma}\label{lemma 1}
    Let $\phi$ be a bounded solution of the problem
    \begin{equation}
    L(\phi) = 0 \quad\hbox{in } \R^2.
    \label{3}\end{equation}
    Then $\phi(x,z)$ is a linear combination of the functions
    $w_x(x)$,
    $Z(x)\cos (\sqrt{\la_1} z )$,
    and
    $Z(x)\sin (\sqrt{\la_1} z )$.
 \end{lemma}

    \noindent{\bf Proof.}
    Let assume that $\phi$ is a bounded function that
    satisfies
    \begin{equation}
    \phi_{zz} + \phi_{xx} + (pw^{p-1} -1 ) \phi = 0 .
    \label{4}\end{equation}
    Let us consider the Fourier transform of $\phi(x,z)$ in the $z$ variable,
    $\hat \phi (x,\xi)$ which is by definition the distribution defined as
    $$
    \langle \hat \phi (x,\cdot ) , \mu \rangle_\R = \langle\phi (x,\cdot ) , \hat
    \mu\rangle_\R\, =\, \int_\R
     \phi (x,\xi )  \hat \mu(\xi) d\xi,
    $$
    where $\mu(\xi)$ is any smooth rapidly decreasing function.
    Let us consider a smooth
    rapidly decreasing function of the two variables $\psi(x,\xi)$.
    Then from equation \equ{4} we find
    $$
    \int_\R
    \langle\hat \phi (x,\cdot ) ,
    \psi_{xx} -\xi^2 \psi + (pw^{p-1} -1 ) \psi \rangle_\R dx = 0.
    $$
    Let  $\varphi(x)$ and $\mu(\xi)$ be  smooth
    and compactly supported functions  such that
    $$
    \{\sqrt{\la_1},
    -\sqrt{\la_1}, 0 \} \cap Supp\, (\mu)= \emptyset .
    $$
    Then we can solve the equation
    $$
    \psi_{xx} -\xi^2 \psi + (pw^{p-1} -1 ) \psi = \mu(\xi)\varphi(x), \quad
    x\in \R,
    $$
    uniquely as a smooth, rapidly decreasing function $\psi(x,\xi)$ such that
    $\psi(x,\xi) = 0 $ whenever
    $\xi \not\in Supp\, (\mu)$. We conclude then that
    $$
    \int_\R
    \langle\hat \phi (x,\cdot ) ,    \mu \rangle_\R \varphi (x) \, dx = 0,
    $$
    so that for all $x\in \R$,
    $< \hat \phi (x,\cdot ) ,    \mu >= 0 $, whenever
    $
    \{\sqrt{\la_1},
    -\sqrt{\la_1}, 0 \} \cap Supp\, (\mu)= \emptyset ,
    $
    in other words
    $$
    Supp\, (\hat \phi (x, \cdot))\subset
    \{\sqrt{\la_1},
    -\sqrt{\la_1}, 0 \}.
    $$
By distribution theory we find then that
$\hat \phi (x, \cdot)$ is a linear combination, with coefficients
depending on $x$, of derivatives up to a finite order of Dirac masses
supported in
$\{\sqrt{\la_1},
-\sqrt{\la_1}, 0 \}.
$
Taking inverse  Fourier transform, we get that
$$
\phi(x,z) =
p_0(z,x) +
p_1(z,x) \cos(\sqrt{\la_1}z) +
p_2(z,x) \sin(\sqrt{\la_1}z),
$$
where $p_j$ are polynomials in $z$ with coefficients depending on $x$. Since $\phi$ is bounded these polynomials are of zero  order,
 i.e. $p_j(z,x)\equiv p_j(x)$, and the bounded functions $p_j$ must satisfy
the equations
$$
L_0(p_0) =0 ,\quad
L_0(p_1) -\la_1 p_1  =0 ,\quad
L_0(p_2) -\la_1 p_2  =0 ,
$$
from where the desired result thus follows.
\qed

  Let   $B(\phi)$ be an operator of the form
    $$
    B(\phi) = b_{1}\partial_{xx}\phi +
    b_{2}\partial_{xz}\phi +
    b_{3}\partial_{x}\phi +
    b_{4}\partial_{z}\phi +
    b_{5}\phi ,
    $$
    where the coefficients $b_i$ are small functions. In the sequel we  will denote $\bb =(b_1,\ldots
,b_5)$ and assume that
    \begin{equation}
    \|\bb\|\equiv \sum_{j=1}^5 \|b_j\|_\infty +\|\nabla b_1\|_\infty +\|\nabla
    b_2\|_\infty < \delta_0, \label{small}\end{equation} where the small
    number $\delta_0$ will be subsequently fixed.
The linear theory used in this paper is based on  a priori estimates for the solutions of the following  problems
\begin{equation}
B(\phi)+ L(\phi) = h ,\quad\hbox{in } \R^2.
\label{11}\end{equation} The results of Lemma \ref{lemma 1} imply
that such estimates without imposing extra conditions on $\phi$
may not exist. The form of the bounded solutions of $L(\phi)=0$
suggests  the following orthogonality conditions:
\begin{equation}
\int_\R \phi(x,z) w_x(x)d\mu(x) = 0 = \int_\R \phi(x,z)
Z(x)d\mu(x), \quad\hbox{for all }z\in  \R,
\label{22}\end{equation} where $d\mu(x)$ is a fixed measure in
$\R$ absolutely continuous with respect to the Lebesque measure.
In the sequel we will in particular  consider
$d\mu(x)=\rho(x)\,dx$ where $\rho$ is a compactly supported
cut-off function, however our next result  applies for a general
$d\mu(x)$ as well. With these restrictions imposed we have the
following result concerning a priori estimates for this problem.
\begin{lemma}\label{lemma2}
There exist constants $\delta_0$ and $C$ such that if the bound
$\equ{small}$ holds and $h\in L^\infty (\R^2)$, then any bounded
solution $\phi$ of problem $\equ{11}$-$\equ{22}$ satisfies
$$
\|\phi\|_\infty \le C
\|h\|_\infty.
$$
\end{lemma}
\noindent{\bf Proof.} We will argue by contradiction.  Assuming
the opposite means that  there are sequences $b_j^n$, $\phi_n$,
$h_n$ such that
$$
\sum_{j=1}^5 \|b_j^n\|_\infty + \|\nabla b_1^n\|_\infty  + \|\nabla
b_2^n\|_\infty \to 0,
$$
$$
\|\phi_n\|_\infty = 1 ,\quad
\|h_n\|_\infty
\to 0,
$$
and
\begin{equation}
B_n(\phi_n)+ L(\phi_n) = h_n, \quad\hbox{in } \R^2 ,
\label{111}\end{equation}
\begin{equation}
\int_\R \phi_n(x,z) w_x(x)d\mu(x) = 0= \int_\R \phi_n(x,z)
Z(x)d\mu(x), \quad\hbox{for all }z\in  \R .
\label{221}\end{equation} Here
$$
B_n(\phi) = b_{1}^n\partial_{xx}\phi +
b_{2}^n\partial_{xz}\phi +
b_{3}^n\partial_{x}\phi +
b_{4}^n\partial_{z}\phi +
b_{5}^n\phi .
$$
Let us assume that
$(x_n,z_n) \in \R^2$ is such that
$$ |\phi_n (x_n,z_n) |\to 1.
$$
We claim that the sequence $x_n$ is bounded. Indeed, if not, using the fact that $L\phi =\Delta \phi-\phi +O(e^{-c|x|})\phi$ and employing elliptic estimates we find that the sequence of functions
$$
\tilde
\phi_n (x,z) =
\phi_n (x_n+x,z_n+z),
$$
converges, up to a subsequence, locally uniformly to a solution $\tilde \phi$
of the equation
$$
\Delta \tilde \phi
-\tilde \phi = 0, \quad\hbox{in } \R^2 ,
$$
whose absolute value attains its maximum  at $(0,0)$, This implies $\tilde\phi \equiv 0$, so that $x_n$ is indeed bounded.
Let now
$$
\tilde
\phi_n (x,z) =
\phi_n (x,z_n+z).
$$
Then $\tilde \phi_n$ converges uniformly over compacts to
a bounded, nontrivial solution $\tilde \phi$ of
$$
L(\tilde\phi)= 0
\quad\hbox{in } \R^2 ,
$$
$$
\int_\R \tilde \phi(x,z) w_x(x)d\mu(x) = 0 =\int_\R
\tilde\phi(x,z) Z(x)d\mu(x), \quad\hbox{for all }z\in  \R .
$$
Lemma \ref{lemma 1} then implies $\tilde \phi \equiv 0$, a
contradiction and the proof is concluded. \qed

\bigskip
Using Lemma \ref{lemma2} we can also find a priori estimates with norms
involving exponential weights. Let us consider
the norm
$$
\|\phi\|_{\sigma, a}  \equiv \| e^{\,{\sigma |x|}+a|z|} \phi
\|_\infty .
$$
where  numbers $\sigma, a\geq 0$ are fixed and will be
subsequently adjusted. In the case $a=0$  we have the following a
priori estimates.

\begin{corollary}\label{corollary 1}
There are numbers $C$ and $\delta_0$ as in Lemma $\ref{lemma2}$
for which, if $\|h\|_{\sigma, 0}  <+\infty$, $\sigma\in [0,1)$,
then a bounded solution $\phi$ of \equ{11}--\equ{22} satisfies
\begin{equation}
\|\phi\|_{\sigma, 0}  +
\|\nabla \phi\|_{\sigma, 0}
\le C \|h\|_{\sigma, 0}.
\label{sigma0 est}
\end{equation}
\end{corollary}
\noindent{\bf Proof.} Again we concentrate on estimates for the problem  \equ{11}--\equ{22}.
We already know that
$$
\|\phi\|_\infty \le C \|h\|_{\sigma, 0}.
$$
We set $\tilde\phi=\phi\|h\|_{\sigma, 0}^{-1}$. Then we have
$$
(L+B)(\tilde\phi)=\tilde h, \quad \mbox{where}\ \|\tilde
h\|_{\sigma, 0}\leq 1,
$$
and also $\|\tilde\phi\|_{\infty}\leq C$. Let us fix a number
$R_0>0$ such that for $x>R_0$ we have
$$
pw^{p-1} (x) < \frac{1-\sigma^2}{2},
$$
which is always possible since $w(x)= O(e^{-c|x|})$. For an arbitrary number $\rho>0$ let us set
$$
\bar \phi (x,z ) = \rho[ \cosh (z/2) + e^{\sigma x} ] + M
e^{-\sigma x},
$$
where $M$ is  to be chosen. Then we find that, reducing $\delta_0$
if necessary,
$$
(L+B)(\bar \phi ) \le - \frac {M (1-\sigma^2)}{4}e^{\,-\sigma x}, \quad \hbox{for } x > R_0 .
$$
Thus
$$
(L+B)(\bar \phi ) \le \tilde h,  \quad \hbox{for } x > R_0,
$$
if
$$
\frac {M(1-\sigma^2)}{4} \ge  \|\tilde h\|_{\sigma, 0} =1.
$$
If we also also assume
$$
M e^{\,-\sigma {R_0}} \ge  \|\tilde \phi\|_\infty ,
$$
we conclude from maximum principle that $\tilde \phi \le \bar
\phi$. Letting $\rho \to 0$ we then get by fixing $M$,
$$
\tilde \phi  \le M  e^{\,-\sigma x },   \quad \hbox{for } x
>0,
$$
hence
$$
\phi\leq M\|h\|_{\sigma,0}e^{\,-\sigma x },   \quad \hbox{for }
x>0.
$$
In a similar way we obtain the  lower bound
$$
\phi\geq -M\|h\|_{\sigma,0}e^{\,-\sigma x },   \quad \hbox{for }
x>0.
$$
Finally,  the same argument  for $x<0$ yields
$$
\|\phi\|_{\sigma, 0}   \le C\|h\|_{\sigma, 0}  ,
$$
while from local elliptic estimates we find
$$
\|\nabla \phi\|_{\sigma, 0}   \le C\|h\|_{\sigma, 0}  ,
$$
and the proof is concluded. \qed

\bigskip


When $a>0$ in the definition of the norm $\|\cdot\|_{\sigma, a}$
then we have the following a priori estimates.

\begin{corollary}\label{corollary 2}
There are numbers $C$,  $\delta_0$ as in Lemma $\ref{lemma2}$, and
$\alpha_0>0$ for which, if $\|h\|_{\sigma, a} <+\infty$,  $\sigma
\in (0,1)$, $a\in [0,a_0)$, then a bounded solution $\phi$ to
problem $\equ{11}$-$\equ{22}$ satisfies
$$
\|\phi\|_{\sigma, a} + \|\nabla \phi\|_{\sigma, a} \le C_\sigma
\|h\|_{\sigma, a}  .
$$
\end{corollary}
\noindent{\bf Proof.} We already know that
$$
\|\phi\|_{\sigma, 0}  + \|\nabla \phi\|_{\sigma, 0}  \le
C\|h\|_{\sigma, a}.
$$
Then we may write
$$
\psi(z) = \int_\R \phi^2(x,z)\, dx,
$$
and differentiate twice weakly to get
$$
\psi''(z) = 2\int_\R \phi_z^2(x,z)\, dx + 2\int_\R \phi_{zz}
\phi(x,z)\, dx . $$ We have
\begin{align}
\int_\R \phi_{zz} \phi\, dx  =
 \int_\R \phi_x^2 \, dx
+ \int_\R (1-pw^{p-1}) \phi^2\, dx  - \int_\R B(\phi)\phi +
\int_\R h\phi .
\end{align}
Integrating by parts once in $x$ we find
\begin{align}
\begin{aligned}
\Big| \int_\R  B(\phi)\phi
\Big|&=\Big|\int_\R[-(b_1\phi)_x\phi_x-(b_2\phi)_x\phi_z +
b_3\phi_x\phi+b_4\phi_z\phi+b_5\phi^2]\Big|\\
&\le C\delta_0 \int_\R (\phi_z^2 + \phi_x^2 + \phi^2 ) \, dx .
\label{b est}
\end{aligned}
\end{align}
Because of the orthogonality conditions (\ref{22}) we also have
that for a certain $\gamma >0$,
 $$
\int_\R \phi_x^2 \, dx + \int_\R (1-pw^{p-1}) \phi^2\, dx  \ge
\gamma \int_\R ( \phi_x^2 + \phi^2 ) \, dx .
$$
Hence, reducing $\delta_0$ if necessary, we find that for a
certain  constant $C>0$
$$
\psi''(z) \ge \frac \gamma 4 \psi(z)  - C \int_\R h^2(x,z)\,dx ,
$$
so that
$$
-\psi''(z) + \frac \gamma 4 \psi(z) \le\frac{C}{\sigma} e^{-2a|z|}
\|h\|_{\sigma,a}^2 .
$$
Since we also know that $\psi $ is bounded by:
$$
|\psi(z)|\leq \frac{C}{\sigma} \|h\|^2_{\sigma,0}, $$ we can use a
barrier of the form $\psi^+(z)=M\|h\|_{\sigma, a}^2 e^{-2a z}+\rho
e^{2a z}$, with $M$ sufficiently large and $\rho >0$ arbitrary, to
get the bound $0\leq \psi\leq \psi^+ $ for $z\geq 0$ and any $a <
\frac{\sqrt{\gamma}}{4}\equiv a_0$. A similar argument can be used
for $z<0$.  Letting $\rho \to 0$ we get then
$$
\int_\R \phi^2(x,z)\, dx\le  C_\sigma e^{-2a|z|} \|h\|_{\sigma,
a}^2, \quad a < a_0.
$$
Elliptic estimates yield that for $R_0$ fixed and large
$$
|\phi(x,z)| \le C_\sigma e^{-a|z|} \|h\|_{\sigma,a} \quad\hbox{for
} |x| < R_0 .
$$
The corresponding estimate in the complementary region can be
found by barriers. For instance in the quadrant $\{x>R_0,z>0\}$ we
may consider a barrier of the form
$$
\bar\phi(x,z) = M \|h\|_{\sigma,a}e^{-( \sigma x  + a z )} +
                \rho  e^{\frac x2 + \frac z2 },
$$
with $\rho>0$ arbitrarily small. Fixing $M$ depending  on
 $R_0$ we find the desired estimate for
$\|\phi\|_{\sigma,a}$ letting $\rho\to 0$. Arguing similarly in
the remaining quadrants is similar. The corresponding bound for
$\|\nabla \phi\|_{\sigma,a}$ is then deduced from local elliptic
estimates. This concludes the proof. \qed



\bigskip

Notice that for a  general right hand side $h$ equation of the
form $L(\phi)+B(\phi)=h$ with the  orthogonality conditions
imposed as above  does not have a solution. On the other hand the
problem
  \begin{equation}
    L(\phi)+B(\phi) = h + c(z) w_x + d(z) Z, \quad\hbox{in } \R^2,
    \label{1}\end{equation}
    under orthogonality conditions
    \begin{equation}
    \int_\R \phi(x,z) w_x(x)d\mu(x) = 0= \int_\R \phi(x,z) Z(x)d\mu(x),
    \quad\hbox{for all }z\in  \R.  \label{2}\end{equation}
has a solution in the sense that for given $h$ one can find
$(\phi, c, d)$ satisfying (\ref{1})--(\ref{2}).


\begin{corollary}\label{est full problem}  There exist $C>0$, $\eta_0>0$,  $\delta_0>0$ as in Lemma $\ref{lemma2}$, and
$\alpha_0>0$ for which, if $\|h\|_{\sigma, \A} <+\infty$,  $\sigma
\in (0,1)$, $a\in [0,a_0)$, and $d\mu(x)=\rho(x)\,dx$ is such that
\begin{equation}
\int_\R
e^{\,-\sigma|x|}[(|w_x|+|Z|)|\rho_{xx}|+2(|w_{xx}|+|Z_{x}|)|\rho_{x}|]\,dx<\eta_0,
\label{meas1}
\end{equation}
  then a bounded solution $\phi$ to problem $\equ{1}$-$\equ{2}$
satisfies
\begin{equation}
\|\phi\|_{\sigma,a} + \|\nabla \phi\|_{\sigma,a} \le C
\|h\|_{\sigma,a}  . \label{esti} \end{equation} Moreover we have
\begin{align}
|c(z)| +     |d(z)| \le C \|h\|_{\sigma,a} e^{\,-a|z|}.
\label{cdesti}
\end{align}

\end{corollary}
\noindent{\bf Proof.} To find a priori estimate (\ref{esti}) we
have to find bounds for the coefficients $c(z)$ and $d(z)$.
Testing  equation (\ref{1}) against $w_x$ and integrating with
respect to $d\mu(x)$ we get
\begin{align*}
\int_\R \phi_{zz} w_x\,d\mu(x)   + \int_\R L_0(\phi) w_x\,d\mu(x)
+ \int_\R B(\phi) w_x\,d\mu(x) &   = \int_\R h w_x\,d\mu(x) \\
&\qquad   +
c(z) \int_\R w_x^2\,d\mu(x).  \end{align*} Let us assume that
$\|h\|_{\sigma, 0} < +\infty$ and that $\phi$ is a bounded
solution. Integrating by parts and using  $L_0(w_x)=0$ and the
orthogonality condition (\ref{2}) we get
\begin{equation}
 c(z) \int_\R
w_x^2\rho\,dx=  \int_\R B(\phi) w_x\rho\,dx +\int_\R
\phi(2w_{xx}\rho_x+w_x\rho_{xx})\,dx- \int_\R h w_x\,dx.
\label{full 1}
\end{equation}
To estimate term $\int_\R B(\phi)w_x\rho,dx$ we  use an argument
similar to that of (\ref{b est}), and to estimate $\int_\R
\phi(w_x\rho)_{xx}\,dx$ we use (\ref{meas1})  to get
$$
|c(z)| \le C\int_\R |hw_x|  + C(\delta_0+\eta_0) (
\|\nabla\phi\|_{\sigma, 0}  + \|\phi\|_{\sigma, 0} ) .
$$
From the a priori estimates (Corollary \ref{corollary 1}) applied
to (\ref{1}) we know that
$$
\|\nabla\phi\|_{\sigma, 0}  + \|\phi\|_{\sigma, 0}  \le C(
\|h\|_{\sigma, 0}  + \|c\|_\infty  + \|d\|_\infty) ,
$$
since
$$\|w_x e^{\,\sigma|x|}\|_\infty, \|Z e^{\,
\sigma|x|}\|_\infty <\infty, $$ for $\sigma\in [0,1)$. Thus,
reducing $\delta_0, \eta_0$  if necessary, we find
$$
\|c\|_\infty  \le C\big(\|h\|_{\sigma, 0}  + (\delta_0+\eta_0)
\|d\|_\infty \big) .
$$

Testing the equation against $Z$ and using exactly the same
argument we find
$$
\|d\|_\infty  \le C(\|h\|_{\sigma, 0}  + (\delta_0+\eta_0)
\|c\|_\infty ) .
$$
Hence
$$
\|d Z \|_{\sigma, 0}   + \|cw_x\|_{\sigma, 0}   \le
C\|h\|_{\sigma, 0},
$$
and the estimate
$$
\|\nabla\phi\|_{\sigma, 0}  + \|\phi\|_{\sigma, 0}  \le C
\|h\|_{\sigma, 0}.
$$
follows.

\bigskip
Finally, if we additionally have $ \|h\|_{\sigma, a} < +\infty$,
we obtain that
$$
\int_\R|hw_x| \le C \|h\|_{\sigma,a} e^{\,-a|z|}.
$$
The same procedure above and the a priori estimates found in
Corollary \ref{corollary 2} then yield
$$
|c(z)| +     |d(z)| \le C \|h\|_{\sigma,a} e^{\,-a|z|},
$$
from  where the relation (\ref{esti}) immediately follows. \qed

Concerning the existence of bounded solutions of
(\ref{1})--(\ref{2}) we have:
\begin{proposition}\label{proposition1} There exists numbers $C>0$,
$\delta_0>0$, $\eta_0>0$ such that whenever  bounds $\equ{small}$,
$(\ref{meas1})$  hold, then given $h$ with $\|h\|_{\sigma, a} <
+\infty$, $\sigma\in (0,1)$, $a\in [0, a_0)$, there exists a
unique bounded solution $\phi=T(h)$ to problem $\equ{1}$-$\equ{2}$
which defines a bounded linear operator of $h$ in the sense that
$$
\|\nabla\phi\|_{\sigma, a} + \|\phi\|_{\sigma, a} \le C
\|h\|_{\sigma, a}.
$$
\end{proposition}
\noindent{\bf Proof.} We will first consider solvability of the
following problem
\begin{equation}
(L+B)(\phi)= h, \quad\hbox{in } \R^2, \label{cor 31}\end{equation}
in the  space $V$,  where $\psi\in V$ if  $\|\psi\|_{\sigma,
0}<\infty$, $\sigma\in (0,1]$ and
\begin{equation}\int_\R \psi(x,z) w_x(x)\rho(x)\,dx = 0= \int_\R \psi(x,z)
Z(x)\rho(x)\,dx, \quad\hbox{for all }z\in  \R, \label{cor 32}
\end{equation}
where the density $\rho(x)$ satisfies the hypothesis of Corollary
\ref{est full problem}. We claim that given $h\in V$ there exists
a unique solution $\phi$ of (\ref{cor 31}) in $V$. We will argue
by approximations. Let us replace $h$ by the function
$h(x,z)\chi_{(-R, R)} (z)$ extended $2R$-periodically to the whole
plane. With this right hand side we can give to the problem
(\ref{cor 31}) a weak formulation in the subspace of $H^1_R\subset
H^1(\R^2)$ of functions that are $2R$-periodic in $z$. To be more
precise let
$$
[\psi, \eta]=\int_{-\infty}^\infty\int_{-R}^R
\nabla\psi\cdot\nabla\eta\,dzdx +\int_{-\infty}^\infty\int_{-R}^R
\psi\eta\,dzdx.
$$
By $W$ we will denote the subspace of functions in $H^1_R$ that
satisfy (\ref{cor 32}).  Then (\ref{cor 31}) can be written in the
form
\begin{equation}
-[(A+K)(\phi),\psi]=\int_{-\infty}^\infty\int_{-R}^R h\psi, \quad
\psi\in W, \label{cor 33}
\end{equation}
where $A:W\to W$ is defined by
$$
[A(\phi),\psi]=\int_{-\infty}^\infty\int_{-R}^R
(\nabla\phi\nabla\psi+b_1\phi_x\psi_x+b_2\phi_x\psi_z)\,dzdx+\int_{-\infty}^\infty\int_{-R}^R(1-b_5)\phi\psi\,dzdx,
\psi\in W,
$$
and $K: W\to W$ is a linear operator defined by
$$
[K(\phi), \psi]=\int_{-\infty}^\infty\int_{-R}^R
[(b_{1x}+b_{2z}-b_3)\phi_x-b_{4}\phi_z]\psi\,dzdx
-p\int_{-\infty}^\infty\int_{-R}^R \phi\psi w^{p-1}\,dzdx, \
\psi\in W.
$$
Using (\ref{small}), (\ref{meas1}),  and the fact that
$\|w^{(p-1)}\|_{(p-1),0}<\infty$ one can show that the operator
$A$ is invertible and  the operator $K$ is compact.

From  Fredholm alternative and Lemma \ref{lemma 1} which in
addition  can be extended periodically to a unique solution
$\phi\in V$, of (\ref{cor 31}) with $h$ replaced by
$h(x,z)\chi_{(-R, R)} (z)$. Letting $R\to +\infty$ and using the
uniform a priori estimates valid for the approximations completes
the proof of the claim.

The existence of a solution to (\ref{1})--(\ref{2}) as well as the
rest of the Proposition  follows  from this claim. Indeed, given
$h$ such that $\|h\|_{\sigma, 0}<\infty$ by $\Pi_V(h)$ we will
denote the orthogonal projection of $h$ onto $V$ (in the sense of
$L^2(\rho\,dx)$ as indicated by (\ref{cor 32})). Using the claim
we can solve then the following problem
$$
(L+B)(\phi)=\Pi_V(h).
$$
Now we only need to chose functions $c(z)$, $d(z)$ such that
$$
(I-\Pi_V)[(L+B)(\phi)]=(I-\Pi_V)(h)+c(z)w_x+d(z) Z.
$$
This ends the proof.

\qed



\bigskip
\begin{remark}\label{remark0}{\em
Incidentally,  Lemma \ref{lemma 1}  helps us to sketch  a proof of
existence of Dancer solutions (essentially that in \cite{dancer})
as follows. Consider the bifurcation problem
$$
\la \Delta v -   v +   v^p  = 0,  \quad \hbox{in } \R^2,
$$
with $\la>0$, which we write, for $w_\la  (x) = w(\sqrt{\la}x)$,
as
$$
\la \Delta \phi + (pw_\la^{p-1} -1)\phi  + N_\la (\phi) = 0,
$$
with $N(\phi) = (w_\la +\phi)^p - w_\la ^p - pw_\la ^{p-1}\phi $.
We consider the space $\XX$ of all functions $\phi$ with
$\|\phi\|_{\sigma, 0}  <+\infty $ which are $T =
2\pi/\sqrt{\la_1}$-periodic and even in the $z$-variable. The
operator $\la \Delta - 1$ has a bounded inverse in $\XX$ and thus
the problem gets rewritten as
$$ \phi + ( \la \Delta -1)^{-1} (pw_\la^{p-1}\phi  + N_\la (\phi) )=0. $$
The derivative of this operator in $\phi$ at $\la=1$ and $\phi =0$
is just $$I+ (\Delta -1)^{-1} (pw^{p-1}),$$ which has the
form: $I+K$, where  $K$ is a  compact operator.
Crandall-Rabinowitz theorem thus yields that bifurcation takes
place at $\la =1$ since Lemma \ref{lemma 1}  implies that this
linearization  has a simple eigenvalue at $\la =1$ with
eigenfunction $Z(x)\cos (\sqrt{\la_1}z)$, and a branch of
nontrivial solutions in $\XX$ constituted by a smooth curve
$\delta \mapsto ( \la (\delta), \phi(\delta))$ with $\phi(0) =0$
and $\la(0) =1$, and $\partial_\delta \phi(0) = Z(x)\cos
(\sqrt{\la_1}z)$. Scaling out $\la(\delta)$ we obtain then
solutions $w_\delta$ to the problem for $\la =1$ with period
$2\pi/\sqrt{\la_\delta\la }$ and even in the $z$-variable. With
the aid of barriers, we find that these solutions have an
expansion of the form
$$
w_\delta (x,z) = w(x) + \delta Z(x)\cos (\sqrt{\la_1}z) +
O(\delta^2) e^{-|x|},
$$
for all small $\delta$. By the standard theory we have that
derivatives of this solutions satisfy the corresponding relations
\begin{align}
\begin{aligned}
w_{\delta, x } (x,z)&= w_x  (x) +\delta Z_x (x
) \cos (\sqrt{\lambda_1} z) + O(\delta^2 e^{\,-|x|}),\\
w_{\delta, x x } (x,z)&= w_{x x } (x ) +\delta Z_{x x } (x ) \cos
(\sqrt{\lambda_1} z) + O(\delta^2 e^{\,-|x|}),
\\
w_{\delta,z} (x,z)&=  -\sqrt{\lambda_1} \delta Z (x )
\sin (\sqrt{\lambda_1}z ) + O(\delta^2 e^{\,-|x|}),\\
w_{\delta, zz} (x,z)&= -\lambda_1 \delta Z (x ) \cos
(\sqrt{\lambda_1} z) + O(\delta^2 e^{\,-|x |}),
\\
w_{\delta, xz } (x,z)&= -\sqrt{\lambda_1} Z_x  (x ) \sin
(\sqrt{\lambda_1}z)+O(\delta^2 e^{\,-|x|}).
\end{aligned}
\label{wdelta}
\end{align}

}
\end{remark}
\setcounter{equation}{0}
\section{Linear theory for  multiple bump-lines }

In this section we will choose a collection of  approximations
$\WW$ which will constitute the basis of our construction of a
solution to Problem \equ{sch}. Let $\alpha$ be a small positive
number and $k$ be a fixed positive integer. Let us consider
functions $f_j\in C^2(\R)$, $j=1,\dots, k$, ordered in the sense
that
\begin{equation}
f_1(z)<f_2(z)<\cdots < f_k(z), \quad\forall \, z\in \R.
\label{co1}\end{equation} These functions represent approximate
locations of the bump lines. Asymptotically these lines are
straight lines with slopes proportional to $\alpha$ and whose
mutual distances are at least $O(d_*)$, where
$$
d_*\equiv \log\frac1{\alpha}.
$$
Thus, in addition to (\ref{co1}), we assume  that for certain
large, fixed constant $M>0$ these functions satisfy
\begin{align}
f_{j+1} (z) -f_{j}(z) \ge 2d_* - M,  \quad j=1,\ldots , k-1,
\label{co3}\end{align} Here and in  what follows we  use the
following  weighted norm for those functions that depend on $z$
only
$$
\|g\|_{\theta_0\A} : = \|e^{\,\theta_0\alpha|z|} g\|_\infty,
$$
where $g:\R\to \R^m$ and $\theta_0>0$ is a fixed number to be
determined later.
Notice that the  factor $\theta_0$ in front of the exponent
$\alpha|z|$ in the definition of the norm will be taken small but  fixed independent
on $\alpha$. For convenience we will denote
$$
\ff = (f_1,f_2,\ldots, f_k).$$  We assume that
\begin{align}
\|\ff''\|_{\theta_0\alpha} \le M\A^2. \label{co3a}
\end{align}
We will further  suppose that the locations of the bump lines are
symmetric with respect to $x$ axis, which in terms of $\ff$ means,
\begin{equation}
\ff(z) = \ff(-z), \quad \forall \, z\in \R.
\label{co2}\end{equation}
 Let us
observe
 that under these assumptions the numbers
$$ \beta_j := f_j'(+\infty ),
$$
are well defined and \begin{equation}|\beta_j| \le M\alpha
.\label{co3aa}\end{equation} In fact, from
(\ref{co3a})--(\ref{co2}) we see that there exist numbers $B_j$
such that we have
$$
f_j(z) = \beta_j|z|   +  B_j + O(e^{-\theta_0\alpha|z|}), \quad
\hbox{ as }|z|\to +\infty.
$$
From (\ref{co3})--(\ref{co2}) it follows also that:
\begin{equation} \beta_{j+1}-\beta_{j}>0, \quad B_{j+1}-B_{j}> 2d_*- M, \quad j=1,\dots,k,
\label{co4}\end{equation} and
\begin{equation}
\||f_j'|-\beta_j\|_{\theta_0\alpha}\leq C\alpha, \qquad j=1,
\dots, k. \label{co3aaa}
\end{equation}
 In addition  we will assume that with  for some $\theta_1$, $\theta_1>\theta_0>0$ with $\theta_0$ defined  as
above  we have
\begin{equation}
\beta_{j+1}-\beta_{j}>\theta_1\alpha. \label{cobeta}
\end{equation}
With functions $f_j$ we will associate a change of
variables defined as follows: for any even function $f(z)$ such that
$$
\|f''\|_{\theta_0\alpha}\leq M\alpha^2, \quad \beta=f'(\infty),
$$
we set
\begin{equation}
 \tX (x, z)= \frac{ x-f(z)}{ \sqrt{1+ \big(\beta \eta(\alpha|z|)\big)^2}},\quad  \tZ(x, z)=
 |z| \sqrt{1+ \big(\beta \eta(\alpha|z|)\big)^2} +\frac{\beta \eta (x-f(z))}{ \sqrt{1+
 \big(\beta
 \eta(\alpha|z|)\big)^2}},
 \label{defxz}
\end{equation}
where $\eta$ is  a smooth cut-off function such that  $\eta(t) =
0$ if $t<T_1$,  $\eta(t)=1$, $t>T_2$, where $T_1, T_2$ will be
chosen later independent on $\alpha$. With $f=f_j$ and
$\beta=\beta_j$, satisfying (\ref{co1})--(\ref{cobeta}),  we will
denote $\tX=\tX_j$, $\tZ=\tZ_j$.

 Let us consider  Dancer
solutions $w_{\delta}(x,z)$, even in $z$, which can be expanded at
main order as
$$w_\delta (x,z) = w(x) + \delta Z(x)\cos ({\sqrt{\la_1}} z) + O(\delta^2) e^{-|x|},\quad |\delta|<\delta_0 .$$
Given $\beta\in \R$ and $\delta$, $|\delta|<\delta_0$ we define
the function
\begin{equation}
w_{\delta,\beta} (x,z) =w_\delta(\tX, \tZ), \quad\mbox{whith} \
(\tX, \tZ) \ \mbox{defined in (\ref{defxz})} .
\label{wbetadelta}\end{equation} With this definition
$w_{\delta,\beta}$ is an exact solution of the problem in the
region where $\eta(\A|z|)= 0$.  If we set $f(z) = \beta |z| +B $
then the function $w_{\delta,\beta} (x,z)$ is also an exact
solution of the problem
 where $\eta(\A|z|)= 1$. Indeed, restricted to the half-plane $\R^2_+=\{z>0\}$, the latter situation corresponds simply to rotating
and translating the axis of the solution $w_\delta $ from the $z$
axis to the line $x = \beta z + B$.

We will  also consider functions $e_j\in C^2(\R)$, $j=1, \dots,
k$, such that
\begin{equation}
\|\ee\|_{\theta_0\A} +\|\ee'\|_{\theta_0\A} +
\|\ee''\|_{\theta_0\A} \le M\A^2, \label{co2a}
\end{equation}
where
$$
\ee = (e_1,e_2,\ldots, e_k).$$ As in (\ref{co2}) we shall assume
that  functions $e_j$ are even,
\begin{equation}
\ee(z) = \ee(-z), \quad \forall \, z\in \R.
\label{co2b}\end{equation}

\bigskip
We are now ready to set up our first approximation.  Given
functions $f_j$, $e_j$ satisfying (\ref{co1})--(\ref{co2}),
(\ref{co2a})--(\ref{co2b}),  and small numbers $\delta_j$,
$j=1,\dots, k$ we let
\begin{equation}
\WW(x,z) \, =\, \sum_{j=1}^k w_{\delta_j,\beta_j} (x,z)
+\sum_{j=1}^k e_j(z)\,Z_j(x,z) \, , \label{WW}\end{equation} where
$\beta_j = f_j'(+\infty)$, $j=1,\dots, k$, and
$$
Z_j(x,z)=Z(\tX_j).
$$
In the sequel we assume
\begin{equation}
\|\vec{\delta} \| \le M\A , \quad \vec{\delta} = (\delta_1,\ldots,
\delta_k). \label{co5}\end{equation} We want to develop a theory
similar to that in the previous section now for the operator
$$
\LL(\phi) = \Delta \phi +  (p\WW^{p-1} -1 ) \phi .
$$
More precisely we define the weighted norm
$$
\|\phi\|_{\sigma, \theta_0\alpha, *} := \Big\| \, \Big(
\sum_{j=1}^k e^{\,-\sigma|x-f_j(z)| - \theta_0\alpha|z|}\Big)^{-1}
\, \phi \, \Big\|_\infty,
$$
and search for a bounded left inverse for a projected problem for
the operator $\LL$ in the space of functions whose
$\|\cdot\|_{\sigma, \theta_0\alpha, *}$ norm is finite. Thus we
consider the problem

\begin{equation}
 \LL(\phi)  = h +\sum_{j=1}^k  c_j(z) \eta_jw_{j,x} + d_j(z)
\eta_j Z_j, \quad\hbox{in } \R^2, \label{01}\end{equation} where
now we {\em do not assume } necessarily orthogonality conditions
on $\phi$. Here and in what follows we denote
\begin{align}
w_{j}(x,z) &= w(x-f_j),
\\ w_{j,x} (x,z) &= w'(x-f_j),
\\
Z_j(x,z)&= Z(x-f_j).
\end{align}
We will now define several cut-off functions that will be
important in the sequel. To this end let $\eta_a^b(s)$ be a smooth
function with $\eta_a^b (s) = 1$ for $|s|<a$ and $=0$ for $|s|>b$,
where $0<a<b<1$. Then, with $ d_*= \log \frac 1\A$, we set
\begin{align}
\begin{aligned}
\rho_j(x,z)&=\eta_a^b\Big( \frac{|\tX_j|}{d_*}\Big), \quad
a=\frac{2^5-1}{2^5}, b = \frac{2^6-1}{2^6},\\
\eta^-_j(x,z)&=\eta_a^b\Big( \frac{|\tX_j|}{d_*}\Big), \quad
a=\frac{2^6-1}{2^6}, b = \frac{2^7-1}{2^7},\\
\eta_j(x,z)&=\eta_a^b\Big( \frac{|\tX_j|}{d_*}\Big), \quad
a=\frac{2^7-1}{2^7}, b = \frac{2^8-1}{2^8},\\
\eta_j^+(x,z)&=\eta_a^b\Big( \frac{|\tX_j|}{d_*}\Big), \quad
a=\frac{2^8-1}{2^8}, b = \frac{2^9-1}{2^9}.\end{aligned}
\label{etaj}\end{align} We will prove the following result:
\begin{proposition}\label{proposition2}
There exist positive constants $\A_0$, $\delta_0>0$ such that if
$\equ{small}$ is satisfied, and  if $\A\in[0,\A_0)$, $\sigma\in
(0,1)$, and constraints $\equ{co1}$-$\equ{co4}$, $\equ{co5}$ hold,
then problem $\equ{01}$ has a solution $\phi= \TT(h)$ which
defines a linear operator of $h$ with $\|h\|_{\sigma,
\theta_0\alpha, *} <+\infty$ and satisfies the estimate
$$
\|\phi\|_{\sigma, \theta_0\alpha,*} + \|\nabla \phi\|_{\sigma,
\theta_0\alpha, *} \le C_{\sigma}\|h\|_{\sigma,\theta_0 \alpha,
*}.
$$
In addition function $\phi$ satisfies the following orthogonality
conditions
\begin{equation}
\int_\R \phi(x,z) w_{j,x}(x,z)\rho_j(x,z)\,dx = 0= \int_\R
\phi(x,z) Z_j(x,z)\rho_j(x,z)\,dx, \quad\hbox{for all }z\in \R .
\label{phiort}
\end{equation}

Besides, the coefficients $c_j(z)$ and $d_j(z)$ in $\equ{01}$ can
be estimated as
\begin{align}
\sum_{j=1}^k(|c_j(z)|+|d_j(z)|)&\leq C\|h\|_{\sigma,
\theta_0\alpha, *} e^{\,-\theta_0\alpha|z|}.
\label{c11}\end{align}
\end{proposition}
\proof{} The main idea in the proof of this proposition is to
decompose problem (\ref{01}) into {\it interior} problems that can
be handled with the help of the theory developed in the previous
section and an {\it exterior} problem and then {\it glue} the
solutions of the subproblems.

Form the definition of the functions $\eta_j$, $\eta^\pm_j$ we
have
\begin{equation}\eta_j\eta_j^- = \eta_j^-, \quad \eta_j^+\eta_j =
\eta_j\label{etas}.
\end{equation}
We search for a solution of \equ{01} of the form
$$
\phi = \sum_{j=1}^k \eta_j \phi_j    + \psi       .
$$
Substituting this expression into equation \equ{01} and arranging terms
we find
\begin{align*}
&  \sum_{j=1}^k \eta_j [\Delta \phi_j + (1-p\WW^{p-1})\phi_j  -
c_j(z) w_{j,x} - d_j(z) Z_j - h ] \\
&\quad + \Big[\Delta \psi    + (1-p\,\WW^{p-1})\psi- \Big(1-
\sum_{j=1}^k \eta_j\Big)h  - \sum_{j=1}^k (2\nabla \eta_j\nabla
\phi_j + \Delta \eta_j \phi_j)\Big]=0 .\end{align*} We will denote
$$
h_j=\eta^+_j h, \quad r_j=p\eta^+_j(w_{j}^{p-1}-\WW^{p-1}).$$ Let
us observe that in the support of $\eta_j$ we have, using
(\ref{etas}),
$$\eta_j h_j=\eta_j h,\quad \eta_j r_j= p\eta_j(w_{j}^{p-1}-\WW^{p-1}) ,$$
Then we find a solution to problem \equ{01} if we solve the
following linear system of equations
\begin{align}
 \Delta \phi_j    -  (1-pw_j^{p-1})\phi_j +r_j\phi_j=
 h_j -p\WW^{p-1}\eta_j^-\psi+
 c_j(z) w_{j,x} + d_j(z) Z_j,
\label{s1} \end{align} in  $\R^2$, for $j=1,\ldots, k$, and
\begin{align}
\Delta \psi    - \Big[1-p\Big(1- \sum_{j=1}^k \eta_j^-\Big)
\WW^{p-1}\Big]\psi = \Big(1- \sum_{j=1}^k \eta_j\Big) h -
 \sum_{j=1}^k (2\nabla \eta_j\nabla  \phi_j + \Delta \eta_j \phi_j
 ),
\label{s2} \end{align} in  $\R^2$. To solve equations \equ{s1} we
denote $\tilde\phi_j=\phi_j+\eta_i^-\psi$ and use
(\ref{s1})--(\ref{s2}) to write the equation for $\tilde\phi_j$
\begin{align}
\begin{aligned}
\Delta\tilde \phi_j    -  (1-pw_j^{p-1})\tilde\phi_j
+r_j\tilde\phi_j&=h_j+\eta_j^-\psi\Big[p(w^{p-1}-\WW^{p-1})+r_j+p\Big(1-\sum\eta_m^-\Big)\WW^{p-1}\Big]\\
&\quad +c_jw_{j,x}+d_jZ_j.\end{aligned}
\label{s3}
\end{align}
We observe that  equation (\ref{s2}) written in terms of
$\tilde\phi_j$ has form
\begin{align}
\Delta \psi    - \Big[1-p\Big(1- \sum_{j=1}^k \eta_j^-\Big)
\WW^{p-1}\Big]\psi = \Big(1- \sum_{j=1}^k \eta_j\Big) h -
 \sum_{j=1}^k (2\nabla \eta_j\nabla  \tilde\phi_j + \Delta \eta_j \tilde\phi_j
 ),
\label{s4} \end{align} since for instance $\nabla\eta_j\nabla
(\eta_j^-\psi) \equiv 0$.

Let us denote
$$
L_j(\phi)=\Delta \phi   -  (1-pw_j^{p-1})\phi,
$$
and  consider first the auxiliary problem
\begin{equation}
 L_j(\phi) = \tilde h + c(z) w_{j,x} + d(z) Z_j, \quad\hbox{in } \R^2,
\label{1j}\end{equation} under orthogonality conditions
\begin{equation}
\int_\R \phi(x,z) w_{j,x}(x,z)\rho_j(x,z)\,dx = 0= \int_\R
\phi(x,z) Z_j(x,z)\rho_j(x,z)\,dx, \quad\hbox{for all }z\in  \R .
\label{2j}\end{equation}
For future references observe that if $\beta_j$ is sufficiently
small then
\begin{align} \rho_j\eta^-_j=\rho_j, \quad \rho_j(1-\eta_j^-)\equiv
0. \label{rhoj}
\end{align}

 We want to solve
(\ref{1j})--(\ref{2j}) using Proposition \ref{proposition1}. To
this end we consider the natural change of coordinates
$$ x\mapsto\tx\equiv  x-f_j  ,\quad z\mapsto  \tz. $$
and set
$$
\phi (x,z) = \ttt \phi(\tx,\tz).
$$
Direct computation then shows that problem \equ{1j}-\equ{2j} is
equivalent to
\begin{equation}
L(\ttt \phi) +B_j(\ttt \phi)= \ttt h + c(\tz) w_\tx(\tx) + d(\tz)
Z(\tx), \quad\hbox{in } \R^2. \label{1jj}\end{equation} under
orthogonality conditions
\begin{equation}
\int_\R \ttt \phi(\tx ,\tz) w_\tx(\tx)\rho(\tx)\,d\tx = 0= \int_\R
\phi(\tx ,\tz) Z(\tx)\rho(\tx)\,d\tx, \quad\hbox{for all }\tz\in
\R,\ \label{2jj}\end{equation} where
$$
\rho(\tx)=\eta_a^b\Big(\frac{\tx}{d^*}\Big),\quad
a=\frac{2^5-1}{2^5}, b=\frac{2^6-1}{2^6} . $$ Here
$$L(\ttt\phi) = \ttt\phi_{\tz\tz} + \ttt\phi_{\tx\tx} + (pw^{p-1} -1 )\ttt\phi,
$$
and
\begin{equation}
B_j(\ttt \phi) =  \Big(\pd{\tx}{z}\Big)^2\ttt \phi_{\tx\tx}
+2\Big(\pd{\tx}{z}\Big)\ttt \phi_{\tx\tz}
+\Big(\pdd{\tx}{z}\Big)\ttt \phi_\tx,\quad
\pd{\tx}{z}=-f_j'(\tz), \quad \pdd{\tx}{z}=-f''_j(\tz).
\label{Bj}\end{equation} The operator $B_j$ satisfies the
assumptions of Proposition \ref{proposition1}, since from
(\ref{co3a})   we have using the notation of (\ref{small}) and
denoting the vector of the coefficients of $B_j$ by ${\mathbf
b}_j$:
\begin{equation}
\|{\mathbf b}_j\|\leq C(\|f_j'\|_\infty+\|f''_j\|_\infty)\leq
C\alpha. \label{Bj1}
\end{equation}
Problem (\ref{1jj})--(\ref{2jj})  has a unique bounded solution
$\tilde \phi = \ttt T_j(\ttt \phi ) $ where $\ttt T_j$ is the
linear operator $T$ predicted by the proposition for $B=B_j$.
Besides, if $\|\ttt h\|_{\sigma, \A} < +\infty$ then we have the
estimate
$$ \|\nabla_{\tx ,\tz } \ttt\phi \|_{\sigma,\theta_0 \A}+
\|\ttt\phi \|_{\sigma, \theta_0\A} \le C \|\ttt h \|_{\sigma,
\theta_0\A} .
$$
going  back  to the original variables we see then that there is a
unique solution to \equ{1j}-\equ{2j},  $\phi = T_j(h)$ where $T_j$
is a linear operator. In addition we have the estimate
$$ \|\nabla \phi \|_{\sigma,\theta_0 \alpha,j}+
\|\phi \|_{\sigma, \theta_0\alpha, j} \le C \| h \|_{\sigma,
\theta_0\alpha, j},
$$
where  $$ \|\phi \|_{\sigma, \theta_0\alpha,j}\, =\,  \| \,
e^{\sigma |x-f_j(z)| + \theta_0\alpha |z| }\, \phi\, \|_\infty .
$$ Moreover, the coefficients $c$ and $d$ are estimated using
(\ref{cdesti}) by
\begin{align}
|c(z)|+|d(z)|\leq  C\|h\|_{\sigma, \theta_0\alpha, j}
e^{\,-\theta_0\alpha|z|}. \label{c1}\end{align}

\medskip
Let us observe that  given $\psi$ we can recast the equations
\equ{s1} for
 $\tilde\phi_j$  as a system of the form
\begin{align}
\begin{aligned}
\tilde\phi_j +T_j(r_j\tilde \phi_j)&=  T_j\Big( h_j
-p\WW^{p-1}\eta_j^-\psi\Big[p(w^{p-1}-\WW^{p-1})+r_j+p\Big(1-\sum\eta_m^-\Big)\WW^{p-1}\Big]\Big),\\
&
\qquad\qquad\qquad\qquad j=1,\ldots , k,\end{aligned} \label{s11}\end{align}

We will solve next equation \equ{s4} for $\psi$ as a linear
operator
$$ \psi = \Psi(\Phi, h),$$ where $\Phi$ denotes the
 $k$-tuple $\Phi = (\tilde\phi_1,\ldots ,\tilde\phi_k)$. To  this end  let us
consider first the problem
\begin{equation}
  \Delta \psi    - (1- \theta )\psi =
 g  \quad\hbox{in } \R^2.
\label{eqpsi}\end{equation}
where
$$
\theta = p\Big(1- \sum_{j=1}^k \eta_j^-\Big) \WW^{p-1}.
$$
Observe that if the number $d_*$ is large enough then  $\theta$ is
uniformly small, indeed $\theta=o(1)$ as $\alpha\to 0$. Let us
assume that $g$ satisfies
$$
|g(x,z)| \le  A\sum_{j=1}^k e^{-\mu|x-f_j(z)|-\theta_0\alpha|z| },
$$
for some $0\le \mu <1$. Then, given that  if $d_*$ is sufficiently
large then number $\theta$ is small, and also that (\ref{co3a})
holds, the use of barriers and elliptic estimates proves that this
problem has a unique bounded solution with
$$
|\nabla \psi(x,z)|+  |\psi(x,z)| \le  C\,  \sum_{j=1}^k
e^{-\mu|x-f_j(z)| -\theta_0\alpha|z|} .
$$
Thus if we take
$$
g =
(1- \sum_{j=1}^k \eta_j) h -
 \sum_{j=1}^k (2\nabla \eta_j\nabla  \tilde\phi_j + \Delta \eta_j \tilde\phi_j ),
$$
we clearly have that
\begin{align*}
|g(x,z)| \le & C\Big[\|h \|_{\sigma, \theta_0\alpha, *} +
o(1)\sum_{j=1}^k \big(\|\tilde\phi_j\|_{\sigma, \theta_0\alpha, j}
+ \|\nabla \tilde \phi_j\|_{\sigma, \theta_0\alpha,
j}\big)\Big]\\
&\qquad \times\sum_{j=1}^k e^{-\mu|x-f_j(z)| -\theta_0\alpha|z|},
\end{align*}
and hence equation \equ{s4} has a unique bounded solution
$$
\psi = \Psi (\Phi, h),
$$
which defines a linear operator in its argument and satisfies the
estimate
\begin{align}
\begin{aligned}
|\Psi(\Phi,h)| \le &C\Big[\|h \|_{\sigma, \theta_0\alpha, *} +
o(1)\sum_{j=1}^k \big(\|\tilde\phi_j\|_{\sigma, \theta_0\alpha, j}
+ \|\nabla\tilde \phi_j\|_{\sigma, \theta_0\alpha,
j}\big)\Big]\\
&\qquad\times\sum_{j=1}^k e^{-\mu|x-f_j(z)| -\theta_0\alpha|z|}.\end{aligned}
 \label{estpsi1}\end{align} In addition, we find that
\begin{equation}
\|\Psi ( \Phi ,h ) \|_{\sigma, \theta_0\alpha, *} \le  C\Big[
\|h\|_{\sigma, \theta_0\alpha, *}+o(1)\sum_{j=1}^k
\big(\|\tilde\phi_j\|_{\sigma,\theta_0 \alpha, j} + \|\nabla
\tilde\phi_j\|_{\sigma,\theta_0 \alpha, j}\big)\Big].
\label{estpsi2}\end{equation}

\medskip
Now we have the ingredients to solve the full system
\equ{s3}-\equ{s4}. Accordingly to \equ{s11} we obtain a solution
if we solve the system in $\Phi$

\begin{equation}
\tilde\phi_j + T_j\big(\eta_j^-\Psi(\Phi, 0
)\chi_j+r_j\tilde\phi_j \big) =
  T_j\big( h_j -\eta_j^-\Psi(0, h )\chi_j \big), \quad
j=1,\ldots , k, \label{ss}\end{equation} where
$$
\chi_j=\Big[p(w^{p-1}-\WW^{p-1})+r_j+p\Big(1-\sum\eta_m^-\Big)\WW^{p-1}\Big],\quad
j=1, \dots, k.
$$
 We consider this system
defined in the space $X$ of all $C^1$ functions $\Phi$ such that
the norm
$$
\|\Phi\|_X := \sum_{j=1}^k \|\nabla \tilde\phi_j\|_{\sigma,
\theta_0\alpha, j} +\| \tilde\phi_j\|_{\sigma, \theta_0\alpha, j},
$$
is finite. System \equ{ss} can be written as
$$
\Phi + \AA(\Phi) = \BB (h ),
$$
where $\AA$ and $\BB$ are linear operators. Thanks to the
estimates for the operators $T_j$ and the bound \equ{estpsi1} we
see that
$$
\| \BB (h)\|_X \le C \|h\|_{\sigma, \theta_0\alpha, *}.
$$
On the other hand we have that
$$
\| \AA (\Phi)\|_X \le C\Big[\sum_{j=1}^k \|\eta_j^-\Psi(\Phi, 0
)\chi_j\|_{\sigma, \theta_0\alpha,
j}+\sum_{j=1}\|r_j\tilde\phi_j\|_{\sigma,\theta_0 \alpha, j}\Big].
$$
Using estimates \equ{estpsi1} and \equ{estpsi2}  we  find
\begin{equation}\sum_{j=1}^k \|\eta_j^-\Psi(\Phi, 0 )\chi_j\|_{\sigma,
\theta_0\alpha, j}\leq o(1)\|\Phi\|_{X}. \label{est chij}
\end{equation} From the definition of $r_j$ and (\ref{co5}) we
get $\|r_j\|_\infty=o(1)$ which implies
$$
\sum_{j=1}^k\|r_j\tilde\phi_j\|_{\sigma,\theta_0 \alpha, j}\leq
o(1)\|\Phi\|_{X}.
$$
Summarizing the last  estimates we obtain
$$
\|\AA(\Phi)\|_X\leq o(1)\|\Phi\|_X,
$$
hence the operator $\AA$ is a  uniformly small operator in the
norm $\|\cdot\|_X$ provided that $\A$ is sufficiently small. We
conclude that system \equ{ss} has a unique solution $\Phi=
\Phi(h)$, which  in addition is a linear operator of $h$ such that
$$
\|\Phi(h) \|_X \le C\| h\|_{\sigma, \theta_0\alpha, *}.
$$
Thus we get a solution to problem \equ{01} by setting
\begin{equation}
\phi = \sum_{j=1}^k \eta_j\tilde \phi_j(h)
+\Big(1-\sum\eta_j^-\Big) \Psi(\Phi(h),h).
\label{phi}\end{equation} Using (\ref{rhoj})  we get
$\rho_j\phi=\rho_j\tilde\phi_j$ hence from and (\ref{2j}) we
obtain (\ref{phiort}). Estimate (\ref{c11}) follows directly from
(\ref{c1}). The proof of the proposition is complete. \qed

\setcounter{equation}{0}
\section{The nonlinear projected problem}\label{sec nonlinear}
Let us recall that our goal is to find  a solution of the problem
$$
S[u]:= \Delta u + u^p - u = 0, \quad p\geq 2, \quad\hbox{in }\R^2,
$$
which is close to the function $\WW$   defined in \equ{WW}.  We
will denote by $E$   the error of approximation at $\WW$:
$$
E:= S[\WW] = \Delta \WW +  \WW^p - \WW .
$$
A main observation we make is that, under the assumptions
\equ{co1}-\equ{co4}, \equ{co5}, $\WW$ is such that
\begin{equation}
E_* : = \|S[\WW]\|_{\sigma, \theta_0\alpha, *} \le
C\A^{2-2\sigma}, \label{e*}\end{equation}
where $\theta_0$ is defined  in (\ref{cobeta}). We will postpone
the proof of (\ref{e*}) for now and for the purpose of the present
section we will simply accept it as a fact.
\medskip
We look for a solution to our problem in the form
$$
u = \WW + \phi,
$$
where $\phi$ is a small perturbation of $\WW$. Thus the equation
for $u$ is equivalent to
\begin{equation}
\LL(\phi)  = E - N(\phi), \quad\hbox{in }\R^2,
\label{nl}\end{equation} where
\begin{equation}
\LL(\phi) := \Delta \phi + (p\, \WW^{p-1} - 1)\phi,
\label{Lb}\end{equation} and
\begin{equation}   N(\phi)
= ( \WW +\phi )^p -  \WW^{p}- p\,\WW^{p-1}\phi.
\label{EN}\end{equation} Rather than solving problem \equ{nl}
directly we consider an intermediate projected version of it,
\begin{equation}
\LL(\phi) = E - N(\phi) +\sum_{j=1}^k  c_j(z) \eta_jw_{j,x} +
d_j(z) \eta_j Z_j,  \quad\hbox{in }\R^2, \label{nl1}\end{equation}
where $\eta_j$ are the cut-off functions defined  by \equ{etaj}
for problem \equ{01}.

\medskip
We will establish next that this nonlinear problem is solvable
with similar estimates to those obtained in Proposition in
\ref{proposition2} for problem \equ{01} with $h=E$, namely  we
show  the following result.

\begin{proposition}\label{proposition 3}
There exist positive numbers $\A_0, \delta_0$, $\sigma_0$, such that for any
number $\sigma\in (0,\sigma_0)$,
$\alpha\in [0,\A_0)$
 and any $\ff$, $\ee$ and $\ddelta$
satisfying constraints $\equ{co1}$-$\equ{co4}$, $\equ{co5}$,
Problem $\equ{nl1}$ has a solution $\phi$ with $\|\phi\|_{\sigma,
\theta_0\alpha, *} \le C_{\sigma'}\A^{2-2\sigma'}$, where
$\sigma'\in(\sigma,3\sigma/2)$ such that
\begin{equation}
\int_\R \phi(x,z) w_{j,x}(x,z)\rho_j(x,z)\,dx = 0= \int_\R
\phi(x,z) Z_j(x,z)\rho_j(x,z)\,dx, \quad\hbox{for all }z\in \R .
\label{phiort2}
\end{equation}

The coefficients $c_{j}$ and $d_j$ can be estimated as follows:
\begin{equation}
\sum_{j=1}^k (|c_j(z)|+|d_j(z)|)\leq C_{\sigma'}\A^{2-2\sigma'}
e^{\,-\theta_0\alpha|z|}. \label{cdj}\end{equation}
\end{proposition}
\proof{} Let us observe that we obtain a solution of the  problem
\equ{nl1} if we solve the fixed point problem
\begin{equation}
\phi    = \TT \big( E - N(\phi)\big):= \mathcal{M}(\phi),
\label{nl2}\end{equation} where $\TT$ is the operator found in
Proposition \ref{proposition2}. Assume that $\|\phi_j\|_{\sigma,
\theta_0\alpha,  *} < 1$, $j=1,2$. We have that
 $$ |\, N(\phi_2)- N(\phi_1)\, | \le  C(|\phi_1| +
|\phi_2|)\,|\phi_2 - \phi_1|, $$ and hence
\begin{equation}
\|\, N(\phi_2)- N(\phi_1)\, \|_{\sigma, \theta_0\alpha,*} \,\le \,
C(\, \|\phi_1\|_{\sigma, \theta_0\alpha, *} + \|\phi_2\|_{\sigma,
\theta_0\alpha, *})\, \|\phi_2 - \phi_1\|_{\sigma,
\theta_0\alpha,*} , \label{ll1}\end{equation} in particular
$$
\|\, N(\phi)\, \|_{\sigma, \theta_0\alpha, *} \,\le \,
\|\phi\|_{\sigma, \theta_0\alpha,  *}^{2} .
$$
Then, from \equ{e*}, the following holds: for each fixed  $\sigma'\in (\sigma,3\sigma/2)$ there exists  a number
 $\nu>0$ such that for all small $\A$ the operator $\mathcal M $ is
a contraction mapping in a region of the form
$$
\BB = \{ \phi \ / \ \|\phi\|_{\sigma, \theta_0\alpha, *} \le  \nu
\A^{2-2\sigma'} \},
$$
and hence a solution of the fixed point problem \equ{nl2} in $\BB$
exists. Furthermore  $\phi$ solves \equ{nl1}, and by
(\ref{phiort}) we find that $\phi$ satisfies (\ref{phiort2}).
The proof of the proposition is
complete.  \qed

\medskip

\setcounter{equation}{0}
\section{Estimates of the error of the initial
approximation}\label{sec error} In this section we will show
estimate (\ref{e*}) announced above for the error
$$
E = S[\WW] =\Delta \WW  - \WW+ \WW^p.
$$

We will denote
\begin{equation}
\label{e1e2}
E_{1}=\Delta\WW, \quad E_2=-\WW+\WW^p,
\end{equation}
 and let for $j=1, \dots, k$,
$$
 \tX_j (x, z)= \frac{ x-f_j(z)}{ \sqrt{1+ \big(\beta_j \eta(\alpha|z|)\big)^2}},\quad  \tZ_j(x, z)=
 |z| \sqrt{1+ \big(\beta_j \eta(\alpha|z|)\big)^2} +\frac{\beta_j \eta (x-f_j(z))}{ \sqrt{1+
 \big(\beta_j
 \eta(\alpha|z|)\big)^2}}.
$$
With a $j$ fixed  we will set for brevity of the notation $f_j=f$,
$\delta_j=\delta$, $\beta_j=\beta$, $\tX_j=\tX$, $\tZ_j=\tZ$ and
let
\begin{equation}
\label{wdef} w_{\delta, \beta} (x-f(z), z):= w_\delta (\tX  (x,
z), \tZ(x, z)).
\end{equation}
%
%

We observe that $\tX(x,z)=\tX(x,-z)$ and $\tZ(x,z)=\tZ(x,-z)$ and
therefore we only need to consider the error of the approximation
assuming $z\geq 0$. For brevity we will denote
\begin{align*}
a_0&=\sqrt{1+(\beta \eta)^2},&
a_1&=\frac{1}{a_0}=\frac{1}{\sqrt{1+(\beta
\eta)^2}}, \\
a_2&=a_1'=\frac{-\alpha\beta^2\eta\eta'}{a_0^3},&
a_3&=a_2'=\frac{-\alpha^2\beta^2[(\eta')^2+\eta\eta'']}{a_0^3}+\frac{3\alpha^2\beta^4\eta\eta'}{a_0^5}.
\end{align*}
Using this notation we can write, whenever $z\geq 0$
$$
\tX=(x-f)a_1, \quad \tZ= za_0+\beta\eta \tX.
$$
We have that
\begin{align*}
\tX _x &=a_1,&
\tX_z &= -{f' }a_1 +\tX a_2, \\
\tZ_x&= {\beta \eta}a_1,& \tZ_z&=a_0+\alpha\beta^2 z\eta'\eta a_1
+\alpha\beta\eta' \tX+\beta\eta\tX_z\\
&\   &\    & =a_0+\alpha\beta^2 z\eta'\eta
a_1-f'\beta\eta a_1+\alpha\beta\eta' \tX+\beta\eta a_2\tX^2,
\end{align*}
and therefore
\begin{align*} A_{11}:= \tX _x^2 +\tX _z^2&=
a_1^2\big(1+(f')^2\big) -2f'a_1a_2\tX+a_2^2\tX^2\\
&=1+a_1^2\big( (f')^2-\beta^2\eta^2\big)+O(\alpha^3)(1+\tX^2)e^{\,-\theta_0\alpha|z|},\\
A_{12}=A_{21}:=\tX_x\tZ_x+\tX_z\tZ_z&= \beta\eta
a_1^2-f'a_1(a_0+\alpha\beta\eta'\tX+\beta\eta\tX_z)+a_2\tX\tZ_z\\
&=\beta\eta-f'+a_1^2\big((f')^2-\beta^2\eta^2\big)+O(\alpha^3)(1+\tX^4)e^{\,-\theta_0\alpha|z|},\\
A_{22}:=\tX^2_z+\tZ_z^2&=
1+2\beta\eta(\beta\eta-f')+\big((f')^2-\beta^2\eta^2)+2\alpha\beta^2
z \eta'\eta+2a_0\alpha\beta\eta' \tX\\
&\qquad+
O(\alpha^3)(1+\tX^4)e^{\,-\theta_0\alpha|z|}.
\end{align*}
 Similarly, we get
\begin{align*}
l_1&:=\Delta\tX= -f'' a_1+a_2(-f'+\tX_z)+\tX a_3,\\
l_2&:= \Delta \tZ=2\beta^2\alpha\eta\eta'a_1+\beta^2\alpha
z(\eta'\eta a_1)' +\beta\eta\Delta
\tX+2\alpha\beta\tX_z\eta'+\alpha^2\beta\eta''\tX .\end{align*}


Denoting ${\mathbf A}=(A_{ij}-\delta_i^j)$, ${\mathbf
l}=(l_1,l_2)$ we have for any function $u=u(\tX(x,z),\tZ(x,y))$,
$$
\Delta_{x,z} u=\Delta_{\tX,\tZ}u+Tr\big({\mathbf A}\cdot
Hess_{\tX,\tZ}(u)\big)+{\mathbf l}\cdot\nabla_{\tX,\tZ} u,
$$
From (\ref{co3a})--(\ref{co3aaa}) we also have the relations
\begin{align*}
{\mathbf
A}=\left[\begin{array}{cc}O(\alpha^2)+O(\alpha^3)(1+\tX^2)
&  O(\alpha)+O(\alpha^3)(1+\tX^4)\\
O(\alpha)+O(\alpha^3)(1+\tX^4) &
O(\alpha^2)(1+\tX^4)\end{array}\right]e^{\,-\theta_0\alpha|z|}.
\end{align*}
and
$$
{\mathbf l}= \big(O(\alpha^2)+O(\alpha^4)\tX,
O(\alpha^3)+O(\alpha^3)\tX\big)e^{\,-\theta_0\alpha|z|}.
$$
Using now the fact that $$ \WW(x,z)=\sum_{j=1}^k
w_{\delta_j}(\tX_j,\tZ_j)+\sum_{j=1}^k e_j(z)Z(\tX_j), $$ for each
$j=1, \dots,k$ and (\ref{wdelta}) and (\ref{co2a}) we  obtain
\begin{equation}
\Delta\WW=\sum_{j=1}^k\Delta_{\tX_j,\tZ_j}w_{\delta_j}(\tX_j,
\tZ_j)+E_{12}, \label{e11}
\end{equation}
where
$$
\|E_{12}\|_{\sigma, \alpha,*}\leq C\alpha^2.
$$
%
We now turn to computing $E_{2}$. For brevity we will set:
$f(u)=-u+u^p$, $p\geq 2$.  We fix a $j$, and if $1< j< k$ we
define sets
$$
A_j=\left\{\, (x,z)\,\left| \, \frac{f_{j-1}(\ve z)+f_j(\ve
z)}{2}\leq x<\frac{f_{j+1}(\ve z)+f_{j}(\ve z)}{2}\,\right.
\right\}\, ,
$$
while when $j=1$ or $j=k$ we set
\begin{align*}A_1&=\left\{\, (x,z)\,\left| \, -\infty< x<\frac{f_{2}(\ve z)+f_{1}(\ve z)}{2}\,
\right.\right\}, \\ A_k&=\left\{\, (x,z)\,\left| \,
\frac{f_{k-1}(\ve z)+f_k(\ve z)}{2}\leq x<\infty\, \right.\right\}
.
\end{align*}
For $x\in A_j$, with $1\leq j\leq k$, fixed we write
\begin{align*}
E_2&=f(\WW)=f(w_{\delta_j,\beta_j})+[f(\WW)-f(w_{\delta_j, \beta_j})]\\
&=f(w_{\delta_j,\beta_j})+f'(w_{\delta_j,\beta_j})(\WW-w_{\delta_j,
\beta_j} )\\
&  \qquad+\frac{1}{2}f''(w_{\delta_j,\beta_j})(\WW-w_{\delta_j,
\beta_j} )^2+E_{2j}
\\
&=\sum_{i=1}^{k}
f(w_{\delta_i,\beta_i})+\left[f'(w_{\delta_j,\beta_j})(\WW-w_{\delta_j,\beta_j})-\sum_{i\neq
j} f(w_{\delta_i, \beta_i})\right]
\\
&  \qquad +\frac{1}{2}f''(w_{\delta_j,\beta_j})(\WW-w_{\delta_j,
\beta_j} )^2+ E_{2j}.
\end{align*}
We first assume that $p >2$.
Since $p>2$ we have
\begin{align*}
|E_{2j}|&\leq C\Big(\max_{i\neq
j}e^{\,-\frac{1}{2}(p-3)|f_i-f_j|}\Big)\Big[\sum_{i\neq
j}|w_{\delta_i,\beta_i}|^3+\sum_{i=1}^k |e_i|^3|Z(\tX_i)|^3\Big]\\
&\leq
C\alpha^{(p-3)}e^{\,-\frac{1}{2}(p-3)\theta_0\alpha|z|}\Big[\sum_{i\neq
j}e^{\,-\sigma|x-f_i|}
e^{\,-\frac{1}{2}(3-\sigma)|f_j-f_i|}+\alpha^3\sum_{i=1}^k
e^{\,-3\theta_0\alpha|z|}e^{\,-3|x-f_i|}\Big]\\
&\leq C\Big(\alpha^{p-\sigma}
e^{\,-\frac{1}{2}(p-\sigma)\theta_0\alpha|z|}\sum_{i=1}^k
e^{-\sigma|x-f_i|}\Big).
\end{align*}
Thus when $\sigma < p-2$ then
$$
\|E_{2j}\|_{\sigma,\theta_0\alpha,*}\leq C\alpha^{p-\sigma}.
$$

Going back to estimation of $E_2$, we have
\begin{align*}
&f'(w_{\delta_j,\beta_j})(\WW-w_{\delta_j,\beta_j})-\sum_{i\neq j}
f(w_{\delta_i,\beta_i})\\&\qquad =\sum_{i\neq
j}[pw_{\delta_j,\beta_j}^{p-1}w_{\delta_i,\beta_i}-w_{\delta_i,\beta_i}^p]+f'(w_{\delta_j,\beta_j})\sum_{i=1}^k
e_iZ(\tX_i)\\ &\qquad=\max_{i\neq
j}O(e^{\,-(p-1)|x-f_j|-|x-f_i|})+\max_{i\neq
j}O(e^{\,-p|x-f_i|})\\
&\qquad\qquad +O(\alpha^2) e^{\,-\theta_0\alpha|z|}e^{\,-|x-f_j|}\sum_{i=1}^k
e^{\,-(1+\mu)|x-f_i|}.
\end{align*}
When $ 0<\sigma<p- 2$ we have for $i\neq j$
\begin{align*}
(p-1)|x-f_j|+|x-f_i|&= \sigma|x-f_j|+(p-1-{\sigma})|x-f_j|+|x-f_i|\\
&\geq \sigma|x-f_j|+\min\{1, (p-1-\sigma)\}(|x-f_j|+|x-f_i|)
\\
&\geq \sigma|x-f_j|+|f_j-f_i|\\
&\geq \sigma|x-f_j|+(|\beta_i-\beta_j||z|+2d_*-M),
\end{align*}
and
$$
p|x-f_i|\geq
\sigma|x-f_i|+\frac{p-\sigma}{2}(|\beta_j-\beta_i||z|+2d_*-M).
$$
Using $|\beta_i-\beta_j|\geq\theta_0\alpha$ we get for $x\in A_j$
and $p> 2$
\begin{equation}
\Delta\WW+f(\WW)=O(\alpha^{2-2\sigma})
e^{\,-\sigma|x-f_i(z)|-\theta_0\alpha|z|}+E_{2j}, \label{sw111}
\end{equation}
hence in $A_j$ we have
\begin{equation}
\label{e-est}
\|E\|_{\sigma,\theta_0\alpha, *}\leq C\alpha^{2-2\sigma}.
\end{equation}
When $p=2$,  we trivially have
$$
0=\|E_{2j}\|_{\sigma,\theta_0\alpha,*}\leq C\alpha^{p-\sigma},
$$
and
\begin{align*}
|x-f_j|+|x-f_i|&= \sigma|x-f_j|+(1-{\sigma})|x-f_j|+|x-f_i|\\
&\geq \sigma |x-f_j|+ (1-\sigma) |f_j-f_i|\\
&\geq \sigma |x-f_j|+ (1-\sigma) (\theta_1 \alpha |z|+2d_*-M).
\end{align*}
Since $\theta_1>\theta_0$, we may choose $\sigma$ small such that
$(1-\sigma) \theta_1 >\theta_0$  and hence (\ref{e-est}) also holds.
The ends the proof of the estimate. \qed

\setcounter{equation}{0}
\section{Dependence of the solution on the parameters}

Now we will study the
dependence of the function $\phi$ found in the proposition above
on the parameters ${\bf f}$, ${\bf e}$ and $\delta$. More
specifically we are interested in establishing the Lipschitz
character of $\phi$ as a function of variables $({\bf f}'', {\bf
e}'')$. We will begin with the observation that the error term
$S[\WW]$ can be written as follows:
$$
S[\WW]=E_1({\bf f}'', {\bf e}''; {\bf f};\vec{\delta}) +E_2({\bf
f}', {\bf e}', {\bf f}, {\bf e};\vec{\delta}),
$$
where for functions ${\mathbf g}$, ${\mathbf h}$ which are even
and a vector parameter $\vec{\delta}=(\delta_1,\dots,\delta_k)$
such that
\begin{equation}
\|{\mathbf g}''\|_\alpha, \|{\mathbf h}''\|_\alpha\leq
C\alpha^2,\qquad \|\vec{\delta}\|\leq C\alpha, \label{lips00}
\end{equation}
we have defined:
\begin{align}
\begin{aligned}
E_1({\bf g}'', {\bf h}''; {\bf f};\vec{\delta})&=-\sum_{i=1}^k
g''_i a_{i1}(1,\beta_i)\cdot\nabla_{\tX_i,\tZ_i}
w_{\delta_i,\beta_i}+\sum_{i=1}^k h'' Z(\tX_i) , \\
 E_2({\bf f}', {\bf e}'; {\bf f}, {\bf e}; \vec{\delta})&=
S[\WW]+\sum_{i=1}^k f''_i
a_{i1}(1,\beta_i)\cdot\nabla_{\tX_i,\tZ_i}
w_{\delta_i,\beta_i}-\sum_{i=1}^k e''
Z(\tX_i).\end{aligned}\label{lips1}
\end{align}
Notice that $E_2$ depends only on ${\mathbf f}, {\mathbf f}',
{\mathbf e}, {\mathbf e}'$, $\vec{\delta}$ and it has exactly the
same expression as in the previous section.
Given ${\mathbf g}$, ${\mathbf h}$, ${\bf{\tilde f}}$,
${\bf{\tilde e}}$,  ${\mathbf f}$, ${\mathbf e}$, $\vec{\delta}$
let us now consider the following problem
\begin{equation}
\LL(\phi) = E_1({\bf g}'', {\bf h}''; {\bf
f};\vec{\delta})+E_2({\bf\tilde f}', {\bf\tilde e}'; {\bf f}, {\bf
e};\vec{\delta}) - N(\phi) +\sum_{j=1}^k c_j(z) \eta_jw_{j,x} +
d_j(z) \eta_j Z_j, \quad\hbox{in }\R^2. \label{lips4}\end{equation}
The solution of this problem  can be obtained by the argument of
Section \ref{sec nonlinear}. In fact the analog of Proposition
\ref{proposition 3} yields $\|\phi\|_{\alpha, \sigma, *}\leq
C\alpha^{2-\sigma'}$, where $\sigma'\in (\sigma,3\sigma/2)$.

We will consider functions ${\mathbf g}_i$, ${\mathbf h}_i$,
${\bf{\tilde f}}_i$, ${\bf{\tilde e}}_i$,  ${\mathbf f}_i$,
${\mathbf e}_i$, and vectors $\vec{\delta}_i$, $i=1,2$ such that
\begin{align}
\begin{aligned}
\|{\mathbf g}''_1-{\mathbf
g}''_2\|_{\theta_0\alpha}+\alpha\|{\bf{\tilde f}}'_1-{\bf{\tilde
f}}'_2\|_{\theta_0\alpha}+\alpha^2\|{\mathbf f}_1-{\mathbf
f}_2\|_{\infty} &\leq C\alpha^2,\\
\|{\mathbf h}''_1-{\mathbf h}''_2\|_{\theta_0\alpha}+\|{\bf{\tilde
e}}'_1-{\bf{\tilde e}}'_2\|_{\theta_0\alpha}+\|{\mathbf
e}_1-{\mathbf
e}_2\|_{\theta_0\alpha}&\leq C\alpha^2,\\
\|\vec{\delta}_1-\vec{\delta}_2\|\leq C\alpha.
\end{aligned}
\label{lips0}
\end{align}
First we want to show that functions $E_1$ and $E_2$ are Lipschitz
functions of its variables. To make it more precise we will
distinguish the norms taken with respect to different functions
${\bf f}$. Thus we will denote in this section:
$$
\|\phi\|_{\sigma, \theta_0\alpha, {\bf f}} := \Big\| \, \Big(
\sum_{j=1}^k e^{\,-\sigma|x-f_j(z)| - \theta_0\alpha|z|}\Big)^{-1}
\, \phi \, \Big\|_\infty, \quad {\bf f}=(f_1,\dots, f_k).
$$
\begin{lemma}\label{lem lipsch}
Under the assumptions (\ref{lips0}) we have the following
estimates:
\begin{align}
\begin{aligned}
&\|E_1({\bf g}_1'', {\bf h}_1;{\bf f}_1;\vec{\delta}_1)-E_1({\bf
g}_2'', {\bf h}_2;{\bf f}_2;\vec{\delta}_2)\|_{\sigma,\theta_0
\alpha, {\bf f}_1}\\
&\quad\leq C(\|{\bf
g}_1''-{\bf g}''_2\|_\alpha+\|{\bf h}_1''-{\bf h}_2''\|_{\theta_0\alpha})
+C(\max_{i=1,2}\|{\bf
g}''_i\|_{\theta_0\alpha}+\max_{i=1,2}\|{\bf
h}''_i\|_{\theta_0\alpha})\|{\bf f}_1-{\bf f}_2\|_{\infty}\\
&\qquad +C\max_{i=1,2}\|{\bf
g}''_i\|_{\theta_0\alpha}\|\vec{\delta}_1-\vec{\delta}_2\|,
\end{aligned}\label{lips0a} \\
\begin{aligned}
&\|E_2({\bf\tilde f}'_1, {\bf\tilde e}'_1; {\bf f}_1, {\bf
e}_1;\vec{\delta}_1)-E_2({\bf\tilde f}'_2, {\bf\tilde e}'_2; {\bf
f}_2, {\bf e}_2;\vec{\delta}_2)\|_{\sigma, \theta_0\alpha, {\bf
f}_1}\\
&\quad\leq C(\alpha^{1-2\sigma}\|{\bf \tilde f }_1'-{\bf \tilde
f}'_2\|_{\theta_0\alpha}+\|{\bf \tilde e}_1'-{\bf \tilde e
}_2'\|_{\theta_0\alpha}) + C(\alpha^{2-2\sigma}\|{\bf f}_1-{\bf
f}_2\|_{\infty}+\|{\bf e}_1-{\bf e}_2\|_{\theta_0\alpha})\\
&\qquad + C\alpha^{2-2\sigma}\|\vec{\delta}_1-\vec{\delta}_2\|.
\label{lips0b}
\end{aligned}
\end{align}
\end{lemma}
\proof{} First notice that to prove Lipschitz estimates
(\ref{lips0a})--(\ref{lips0b}) it suffices to show corresponding
inequalities varying one component  at a time (for instance
considering ${\bf g}_i$, $i=1,2$  and fixing the rest of the
parameters). When ${\bf\tilde f}$ and ${\bf f}$ are fixed then the
estimates are rather standard. We will therefore concentrate on
the case when all parameters except ${\bf\tilde f}$ and ${\bf f}$
are fixed.

We will indicate the dependence of the change of variables $\tX_j,
\tZ_j$ on ${\bf f}$ by writing $\tX_j({\bf f}), \tZ_j({\bf f})$.
Similarly, we will write ${\bf A}_{j}({\bf f}, {\bf\tilde f}')$
for the components of the matrix ${\bf A}_j$. Finally we will set
$$
{\bf\tilde l}_j({\bf f}, {\bf\tilde f}')={\bf l}_j({\bf f},
{\bf\tilde f}', {\bf g}'')+g''_i a_{i1}(1,\beta_i).
$$

We will gather first some estimates that are important in the
proof. We have
\begin{align}
\begin{aligned}
|w(\tX_j({\bf f}_1)-w'(\tX_j({\bf f}_2)|+|w(\tX_j({\bf
f}_1)-w'(\tX_j({\bf f}_2)|&\leq C |w(\tX_j({\bf
f}_1)||{\bf f}_1-{\bf f}_2|,\\
|{\bf\tilde l}_j({\bf f}_1, {\bf\tilde
f}')-{\bf\tilde l}_j({\bf f}_2,
{\bf\tilde f}')|&\leq C\alpha^2|{\bf f}_1-{\bf f}_2|,\\
|{\bf\tilde l}_j({\bf f}, {\bf\tilde f}'_1)-{\bf \tilde l}_j({\bf
f}, {\bf\tilde f}'_2)|&\leq C\alpha|{\bf\tilde f}'_1-{\bf\tilde
f}'_2|. \label{lips0c}
\end{aligned}
\end{align}

From the definition of $E_1$ and (\ref{lips0c}) one gets
(\ref{lips0a}). Estimate (\ref{lips0b}) is somewhat tedious but
rather standard. For instance let us consider a typical term
coming from $\Delta \WW$:
\begin{align*}
&|Tr\big({\mathbf A}_\cdot
Hess_{\tX_j,\tZ_j}(w_{\delta_j,\beta_j})\big)({\bf f}_1,
{\bf\tilde f}'_1)-Tr\big({\mathbf A}_j\cdot
Hess_{\tX_j,\tZ_j}(w_{\delta_j,
\beta_j})\big)({\bf f}_2, {\bf\tilde f}'_2)|\\
&\leq|Tr\big({\mathbf A}_j\cdot
Hess_{\tX_j,\tZ_j}(w_{\delta_j,\beta_j})\big)({\bf f}_1,
{\bf\tilde f}'_1)-Tr\big({\mathbf A}_j\cdot
Hess_{\tX_j,\tZ_j}(w_{\delta_j, \beta_j})\big)({\bf f}_2, {\bf\tilde f}'_1)|\\
&\qquad+|Tr\big({\mathbf A}_j\cdot Hess_{\tX_j,\tZ_j}(w_{\delta_j,
\beta_j} )\big)({\bf f}_2, {\bf\tilde f}'_2)-Tr\big({\mathbf
A}_j\cdot Hess_{\tX_j,\tZ_j}(w_{\delta_j, \beta_j})\big)({\bf
f}_2, {\bf\tilde f}'_1)|\\
&\qquad=B_1+B_2.
\end{align*}
Using the (\ref{lips0c}) we can estimate the first term above by
\begin{align*}
B_1&\leq C|{\mathbf A}_{j,11}({\bf f}_1, {\bf\tilde
f}'_1)-{\mathbf A}_{j,11}({\bf f}_2, {\bf\tilde
f}'_1)||w(\tX_j)({\bf f}_1)|
\\
&\qquad +C|{\mathbf A}_{j,11}({\bf f}_2, {\bf\tilde
f}'_1)||w''(\tX_j)({\bf f}_1)-w''(\tX_j)({\bf f}_2)|
\\
 &\qquad +C\delta_j|{\mathbf A}_j({\bf f}_1, {\bf\tilde
f}'_1)-{\mathbf A}_j({\bf f}_2, {\bf\tilde f}'_1)||Z(\tX_j)({\bf
f}_1)|\\
&\qquad + C\delta_j|{\mathbf A}_j({\bf f}_2, {\bf\tilde
f}'_1)||Z(\tX_j)({\bf f}_1)-Z(\tX_j)({\bf
f}_2)|\\
&\leq C\alpha^2|{\bf f_1}-{\bf
f}_2|e^{\,-\theta_0\alpha|z|}\big|w\big(\tX_j({\bf
f}_1)\big)\big|^{\sigma}.
\end{align*}
In an analogous way we can estimate $B_2$. Another typical term we
have to deal with is (considered in the set $A_j$ defined in the
previous section) the following
\begin{align*}
&|f'(w_{\delta_j,\beta_j})(\WW-w_{\delta_j, \beta_j} )({\bf
f}_1)-f'(w_{\delta_j,\beta_j})(\WW-w_{\delta_j, \beta_j} )({\bf
f}_2)|\\
&\leq |f'(w_{\delta_j,\beta_j})({\bf
f}_1)-f'(w_{\delta_j,\beta_j})({\bf f}_2)||(\WW-w_{\delta_j,
\beta_j} )({\bf f}_1)|\\
&\qquad +|f'(w_{\delta_j,\beta_j})({\bf f}_2)||(\WW-w_{\delta_j,
\beta_j} )({\bf f}_1)-(\WW-w_{\delta_j, \beta_j} )({\bf f}_2)|\\
&\leq C\alpha^{2-\sigma}e^{\,-\theta_0\alpha|z|}|{\bf f}_1-{\bf
f}_2|\big|w\big(\tX_j({\bf f}_1)\big)\big|^\sigma+ C\max_{i\neq
j}e^{\,-|f_{1j}-f_{1i}|}|{\bf
f}_1-{\bf f}_2|\big|w\big(\tX_j({\bf f}_1)\big)\big|^\sigma\\
&\leq C\alpha^{2-\sigma}e^{\,-\theta_0\alpha|z|}|{\bf f}_1-{\bf
f}_2|\big|w\big(\tX_j({\bf f}_1)\big)\big|^\sigma.
\end{align*}

Other terms in the definition of
$E_2$ are handled similarly and from this we get (\ref{lips0b}).
We leave the details  to the reader. \qed

\bigskip

In what follows we will emphasize the dependence of $\phi$ and
$c_j, d_j$ on parameters by writing
\begin{align}
\begin{aligned}
\phi^{(i)}&=\phi({\bf g}''_i, {\bf h}''_i;{\bf \tilde f}'_i,
{\bf\tilde e}'_i; {\bf f}_i, {\bf e}_i; \vec{\delta}_i), \\
c^{(i)}_j&=c_j^{(i)}({\bf g}''_i, {\bf h}''_i;{\bf \tilde f}'_i,
{\bf\tilde e}'_i; {\bf f}_i, {\bf e}_i;\vec{\delta}_i), \\
d_j^{(i)}&=d_j^{(i)}({\bf g}''_i, {\bf h}''_i; {\bf \tilde f}'_i,
{\bf\tilde e}'_i;{\bf f}_i, {\bf e}_i; \vec{\delta}_i).
\end{aligned}\label{lips01}
\end{align}

\begin{proposition}\label{proposition 3a}
Let $\phi$,  be the solution of (\ref{lips4}). Then for
$j=1,\dots, k$ functions $\phi$, $c_j$, $d_j$ are continuous with
respect to the parameters ${\bf g}''$, ${\bf h}''$,${\bf{\tilde
f}}'$, ${\bf{\tilde e}}'$,  ${\mathbf f}$, ${\mathbf e}$. Moreover
assuming (\ref{lips0}), in the notation of (\ref{lips01})  we have
the following estimates
\begin{align}
\begin{aligned}
&\|\phi^{(1)}-\phi^{(2)}\|_{ \sigma, \theta_0\alpha,{\bf f}_1}+
\|\nabla\phi^{(1)}-\nabla\phi\|_{ \sigma,\theta_0\alpha {\bf f}_1}\\
&\qquad\leq C \alpha^{-2\sigma'}(\|{\bf g}_1''-{\bf
g}_2''\|_{\theta_0\alpha}+\|{\bf h}_1''-{\bf
h}_2''\|_{\theta_0\alpha})\\
&\qquad\qquad +C\alpha^{-2\sigma'}(\alpha\|{\bf \tilde f}_1'-{\bf \tilde
f}_2'\|_{\theta_0\alpha}+\|{\bf \tilde e}_1'-{\bf \tilde
e}_2'\|_{\theta_0\alpha})
\\
&\qquad\qquad +C\alpha^{-2\sigma'}(\alpha^2\|{\bf f}_1-{\bf
f}_2\|_{\infty}+\|{\bf e}_1-{\bf e}_2\|_{\theta_0\alpha})
\\
&\qquad\qquad +C\alpha^{2-2\sigma'}\|\vec{\delta}_1-\vec{\delta}_2\|,
\label{lips5}
\end{aligned}
\end{align}
where $\sigma'\in (\sigma, 3\sigma/2)$.
\end{proposition}
\proof{} First observe that it is sufficient to prove Lipschitz
dependence considering dependence of $\phi$ on different
components taken separately. For instance let us fix ${\bf \tilde
f}$, ${\bf \tilde e}$, ${\bf f}$ and ${\bf e}, \vec{\delta}$ and
let vary ${\bf g}_i$, ${\bf h}_i$. Then the
 proof of the  result follows directly from Proposition
(\ref{proposition2}) and Proposition (\ref{proposition 3}) applied
to the equation for the difference $\phi({\bf g}_1'', {\bf
h}_1'';{\bf\tilde f}', {\bf \tilde e}';  {\bf f}, {\bf
e};\vec{\delta})-\phi({\bf g}_2'', {\bf h}_2'';{\bf\tilde f}',
{\bf \tilde e}'; {\bf f}, {\bf e};\vec{\delta})$. We leave the
details to the reader.

We will now consider Lipschitz dependence of $\phi$ on ${\bf\tilde
f}', {\bf\tilde e}'$, ${\bf e}$. Let us fix ${\mathbf g}$ and
${\mathbf h}$ as above as well ${\mathbf f}$ and $\vec{\delta}$.
We will denote by $\phi^{(i)}$ the solution of the following
problem in $\R^2$:
\begin{equation}
\LL(\phi^{(i)}) = E_1({\bf g}'', {\bf h}''; {\bf
f};\vec{\delta})+E_2({\bf{\tilde  f}}'_i, {\bf{\tilde e}}'_i; {\bf
f}, {\bf{e}}_i;\vec{\delta}) - N(\phi^{(i)}) +\sum_{j=1}^k
c^{(i)}_j(z) \eta_jw_{j,x} + d^{(i)}_j(z) \eta_j Z_j . \label{lips6}\end{equation}
It  follows from Lemma \ref{lem lipsch} that
\begin{align}
\begin{aligned}
&\|E_2({\bf{\tilde  f}}'_1, {\bf{\tilde e}}'_1; {\bf f},
{\bf{e}}_1;\vec{\delta})-E_2({\bf{\tilde  f}}'_2, {\bf{\tilde
e}}'_2; {\bf f}, {\bf{e}}_2;\vec{\delta})\|_{\sigma,
\theta_0\alpha, {\bf f}_1}\\
\quad&\leq
C\alpha^{-2\sigma'}(\alpha\|{\bf{\tilde f}}'_1-{\bf{\tilde
f}}'_2\|_{\theta_0\alpha}+\|{\bf{e}}_1-{\bf{e}}_2\|_{\theta_0\alpha})+\|{\bf{\tilde
e}}'_1-{\bf{\tilde e}}'_2\|_{\theta_0\alpha}). \end{aligned}
\label{lips8}
\end{align}
Applying  now the theory developed for the linear problem for the
difference $\phi^{(1)}-\phi^{(2)}$ we get
\begin{align}
\begin{aligned}
&\|\phi^{(1)}-\phi^{(2)}\|_{\sigma,\theta_0\alpha, {\bf f}_1
}+\|\nabla\phi^{(1)}-\nabla\phi^{(2)}\|_{\sigma,
\theta_0\alpha,{\bf f}_1}\\&\quad \leq
C\alpha^{-2\sigma'}(\alpha\|{\bf{\tilde f}}'_1-{\bf{\tilde
f}}'_2\|_{\theta_0\alpha}+\|{\bf{e}}_1-{\bf{e}}_2\|_{\theta_0\alpha})+\|{\bf{\tilde
e}}'_1-{\bf{\tilde e}}'_2\|_{\theta_0\alpha}).\end{aligned} \label{lips9}
\end{align}

Finally, it remains to consider the Lipschitz dependence of
solutions on ${\bf e}$ and ${\bf f}$. This case is somewhat more
complicated since the norms involved depend on the variables
$\tX_j$, which in turn depend on the functions $f_j$.

Now, by $\phi^{(i)}$ we will denote solutions of
\begin{align}
\begin{aligned}
\LL_{{\bf f}_i}(\phi^{(i)}) &= E_1({\bf g}'', {\bf h}''; {\bf
f}_i;\vec{\delta})+E_2({\bf{\tilde  f}}', {\bf{\tilde e}}'; {\bf
f}_i, {\bf{e}};\vec{\delta}) - N(\phi^{(i)})\\
&\qquad +\sum_{j=1}^k
c^{(i)}_j(z) \eta_j^{(i)}w^{(i)}_{j,x} + d^{(i)}_j(z) \eta^{(i)}_j
Z^{(i)}_j, \quad\hbox{in }\R^2,\end{aligned} \label{lips10}\end{align} with
fixed functions ${\bf{g}}$, $\bf{h}$, ${\bf{\tilde f}}$,
${\bf{\tilde e}}$, and function ${\bf f}_i$ ,   $i=1,2$ such that
\begin{equation} \|{\bf{f}}_1-{\bf{f}}_2\|_{\infty}\leq C.
\label{lips11}\end{equation} The equation for the difference
$\phi^{(1)}-\phi^{(2)}\equiv \tilde \phi$ can be written in the
form:
\begin{align}
\begin{aligned}
{\mathcal L}_{{\bf f}_1}\tilde \phi&=E_1({\bf g}'', {\bf h}'';
{\bf f}_1;\vec{\delta})-E_1({\bf g}'', {\bf h}''; {\bf
f}_2;\vec{\delta})\\
&\qquad +E_2({\bf{\tilde f}}', {\bf{\tilde e}}'; {\bf
f}_1, {\bf{\tilde e}};\vec{\delta}) -E_2({\bf{\tilde  f}}',
{\bf{\tilde e}}'; {\bf f}_2, {\bf{\tilde e}};\vec{\delta})
\\
&\qquad+ p[(\WW^{(1)})^{p-1}-(\WW^{(2)})^{p-1}]\phi^{(2)} +[N(\phi^{(1)})-N(\phi^{(2)})]\\
&\qquad+\sum_{j=1}^k
c_j^{(2)}(\eta_j^{(1)}w_{j,x}^{(1)}-\eta_j^{(2)}w_{j,x}^{(2)})+\sum_{j=1}^k
d^{(2)}_j(\eta_j^{(1)}Z_{j}^{(1)}-\eta_j^{(2)}Z_{j}^{(2)})\\
&\qquad +\sum_{j=1}^k
(c_j^{(1)}-c_j^{(2)})\eta_j^{(1)}w_{j,x}^{(1)}+\sum_{j=1}^k
(d_j^{(1)}-d^{(2)}_j)\eta_j^{(1)}Z_{j}^{(1)}.
\end{aligned}
\label{lips12}
\end{align}
Observe that if we denote ${\bf \tilde r}={\bf f}_1-{\bf f}_2$
then we have $\|{\bf \tilde r}\|_{\theta_0\alpha}\leq C$ and
$$
|\WW^{(1)}-\WW^{(2)}|\leq C|{\bf \tilde r}|\sum_{j=1}^k
e^{\,-|x-f_{1j}|}.
$$
It follows
\begin{align}
\|\WW^{(1)}-\WW^{(2)}\|_{\sigma, \infty, {\bf f}_{1}}\leq C\|{\bf
\tilde r}\|_{\infty}. \label{lips13}
\end{align}
Similar estimates hold for other terms involved in the right hand
side of (\ref{lips12}).
Now, we will further decompose
$$
\tilde\phi=\tilde\phi_0+\sum_{j=1}^k\tilde\phi_j\rho_j^{(1)} w_{j,x}^{(1)},
$$
where
$$
\tilde\phi_j=\frac{\int_{\R}\phi^{(2)}(\rho_j^{(1)}w_{j,x}^{(1)}-\rho_j^{(2)}w_{j,x}^{(2)})\,dx}{\int_{\R}
(\rho_j^{(1)}w_{j,x}^{(1)})^2\,dx}.
$$
The Lipschitz character of functions $\tilde\phi_j$, $j=1,\dots,k$ can be established directlly from the definition. On the other hand function $\tilde\phi_0$ satisfies an equation similar to (\ref{lips12}) together with the orthogonality conditions.
This, estimates in Proposition
\ref{proposition 3} and Lemma \ref{lem lipsch} yield finally
\begin{align}
\|\phi^{(1)}-\phi^{(2)}\|_{\sigma, \theta_0\alpha, {\bf f}_{1}}+
\|\nabla\phi^{(1)}-\nabla\phi^{(2)}\|_{\sigma, \theta_0\alpha,
{\bf f}_{1}}\leq C\alpha^{2-2\sigma'}\|{\bf\tilde
r}\|_{\infty}e^{\,-\theta_0\alpha|z|}. \label{lips14}
\end{align}

\setcounter{equation}{0}
\section{Projection of the error on the elements of the kernel}
Let us now compute the projections of the error on the elements of
the approximate kernel. Thus for each $j=1, \dots, k$ we will
compute \begin{align*}
\Pi_{f_j}&:=\int_\R S[\WW] w_x(\tX_j)\,dx,
\\
\Pi_{e_j}&:=\int_\R S[\WW] Z(\tX_j)\,dx,
\end{align*}
\medskip
For a fixed $j$ we begin with $\Pi_{f_j}$. We will use
decomposition of $S[\WW]$  similar to the one in section \ref{sec
error}:
\begin{align}
\begin{aligned}
S[\WW]&=\sum_{i=1}^k \big[Tr({\mathbf A}_i\cdot
Hess_{\tX_i,\tZ_i}(w_{\delta_i,\beta_i})+{\mathbf
l}_i\cdot\nabla_{\tX_i,\tZ_i}
w_{\delta_i,\beta_i}\big]\\
&\qquad +\sum_{i=1}^k \Delta_{x,z}
[e_i(z)Z(\tX_i)]\\
&\qquad + [f(\WW)-\sum_{i=1}^k f(w_{\delta_i, \beta_i})]
\\
&:=S_{1j}[\WW]+S_{2j}[{\bf e},Z(\tX_1),\dots, Z(\tX_k)]+
S_{3j}[\WW].
\end{aligned}
\label{expsww}
\end{align}
Using (\ref{wdelta}) we obtain:
\begin{align*}
Hess_{\tX_i, \tZ_i}(w_{\delta_i,\beta_i})&=\left[\begin{array}{cc}
w''(\tX_i) & 0
\\
0& 0\end{array}\right]\\&\quad+\delta_i\left[\begin{array}{cc}
Z''(\tX_i)\cos(\sqrt{\lambda_1}\tZ_i) &
-\sqrt{\lambda_1}Z'(\tX_i)\sin(\sqrt{\lambda_1}\tZ_i)\\
-\sqrt{\lambda_1}Z'(\tX_i)\sin(\sqrt{\lambda_1}\tZ_i) & -\lambda_1
Z(\tX_i)\cos(\sqrt{\lambda_1}\tZ_i)\end{array}\right]\\
&\quad+O(\delta_i^2) e^{\,-|\tX_i|}.
\end{align*}
From the formulas for $\tX_i, \tZ_i$ and their derivatives we also
get
\begin{align*}
{\mathbf A}_i&=\left[\begin{array}{cc} a_{i1}^2\big(
(f'_i)^2-\beta_i^2\eta^2\big) &
\beta_i\eta-f'_i+a_{i1}^2\big((f'_i)^2-\beta_i^2\eta^2\big)
\\
\beta_i\eta-f'_i+a_{i1}^2\big((f'_i)^2-\beta_i^2\eta^2\big) &
2\beta_i\eta(\beta_i\eta-f_i')+\big((f'_i)^2-\beta_i^2\eta^2)+2\alpha\beta_i^2
z \eta'\eta
\end{array}\right]\\
&\qquad\qquad\qquad+\left[\begin{array}{cc} 0 & 0
\\
0 & 2a_0\alpha\beta_i\eta'
\tX_i\end{array}\right]+\left[\begin{array}{cc} O(\alpha^3)&
O(\alpha^3)\\
O(\alpha^3)& O(\alpha^3)\end{array}\right](1+\tX_i^2+\tX_i^4)
e^{\,-\theta_0\alpha|z|}\\
&\equiv {\mathbf A}_{i0}+{\mathbf
A}_{i1}+{\mathbf A}_{i2}.
\end{align*}
Similarly we have
$$
\nabla_{\tX_i,\tZ_i}w_{\delta_i,\beta_i}=(w'(\tX_i),0)+\delta_i\big(Z'(\tX_i)\cos(\sqrt{\lambda_1}\tZ_i),
-\sqrt{\lambda_1}Z(\tX_i)\sin(\sqrt{\lambda_1}\tZ_i)\big)+O(\delta_i^2)
e^{\,-|\tX_i|},
$$
and
$$
{\mathbf
l}_i=-f''_ia_{i1}(1,\beta_i\eta)+O(\alpha^3)(1+|f'_i-\beta_i\eta|+|\tX_i|)e^{\,-\theta_0\alpha|z|}.
$$
We observe that the entries of ${\bf A}_i$ and ${\bf l}_i$ are
polynomials with respect to $\tX_i$ and they depend on $f_i$ only
through $\tX_i$.
In the sequel we will use also expansion
$$
\cos(\sqrt{\lambda_1}\tZ_i)=\cos(\sqrt{\lambda_1}a_{i0}
z)+O(\alpha)|\tX_i|, \quad
\sin(\sqrt{\lambda_1}\tZ_i)=\sin(\sqrt{\lambda_1}a_{i0}
z)+O(\alpha)|\tX_i|.
$$

Since functions $w''$, $Z$ and $Z''$ are even functions of its
arguments and $w'$ is an odd function therefore, changing
variables $x\mapsto \tX_j$, we get
\begin{align*}
\int_\R Tr\big({\mathbf A}_j\cdot Hess(w_{\delta_j, \beta_j})\big)
w'(\tX_j)\,dx&=-2\sqrt{\lambda_1}a_{j0}\delta_j\big[\beta_j\eta-f'_j+a^2_{j1}\big((f'_j)^2-\beta_j^2\eta^2\big)\big]\\
&\qquad\qquad\times
\sin(\sqrt{\lambda_1}a_{i0}z)\int_R w'(s)Z'(s)\,ds
\\
&\qquad
+\delta_jO(\delta_j+\beta_j)(f_j'-\beta_j\eta)\\
&\qquad+\alpha^2O(\delta_j+\alpha)e^{\,-\theta_0\alpha|z|}.
\end{align*}
When $i\neq j$ we have
$$
\int_\R w'(\tX_j) w(\tX_i)\,dx=O(e^{\,-|f_j-f_i|}), \qquad \int_\R
Z(\tX_j) w(\tX_i)\,dx=O(e^{\,-(1+\mu)|f_j-f_i|}),
$$
with similar formulas for the integrals involving $Z'(\tX_j)$ and
$Z''(\tX_j)$. It follows then
\begin{align*}
\int_\R Tr\big({\mathbf A}_i\cdot Hess(w_{\delta_i, \beta_i})\big)
w'(\tX_j)\,dx&=\big((f'_i)^2-\beta_i^2\eta^2\big)O(e^{\,-|f_j-f_i|})\\
&\qquad+(f'_i-\beta_i\eta)(\alpha+\delta_i)\delta_i
O(e^{\,-|f_j-f_i|})\\
&\qquad+O(\alpha^3+\delta_i\alpha^2)
e^{\,-\theta_0\alpha|z|}.
\end{align*}
We have also
\begin{align*}
\int_\R {\mathbf l}_j\cdot\nabla_{\tX_j,\tZ_j} w_{\delta_j,
\beta_j}w'(\tX_j)\,dx&=-f''_j\Big[\int_\R
(w')^2+\delta_j\cos(\sqrt{\lambda_1}a_{i0}z)\int_\R Z'
w'\\
&\qquad\qquad+\delta_jO(\delta_j+\beta_j)\Big]\\
&\qquad+O(\alpha^3)(f_j'+1)e^{\,-\theta_0\alpha|z|},
\end{align*}
and similarly for $i\neq j$ we get
$$ \int_\R
{\mathbf l}_j\cdot\nabla_{\tX_i,\tZ_i} w_{\delta_i,
\beta_i}w'(\tX_j)\,dx=-f''_iO(e^{\,-|f_j-f_i|})+O(\alpha^3)(f_i'+1)e^{\,-\theta_0\alpha|z|}.
$$
To compute the corresponding projections of $S_{2j}$ we notice
that
\begin{align}
\begin{aligned}
\Delta_{x,z} \big(e_i(z)Z(\tX_i)\big)&=e''_i(z)
Z(\tX_i)+2e_i'(z)Z'(\tX_i)\tX_{iz}+e_i(z)\Delta_{x,z}Z(\tX_i)\\
&=e''_iZ(\tX_i)+\big[2e_i'(-f_i'a_{i1}+\tX_i
a_{i2})-e_if''_ia_{i1}\\
&\qquad\qquad+e_i(f_i'+\tX_i)e^{\,-\theta_0\alpha|z|}O(\alpha^3)\big]Z'(\tX_i)\\
&\qquad+e_i\big[1+a_{i1}^2\big((f'_i)^2-\beta_i^2\eta^2\big)+O(\alpha^3)(1+\tX^2_i+\tX_i^4)e^{\,-\theta_0\alpha|z|}\big]Z''(\tX_i).
\end{aligned}\label{lapez1}
\end{align}
Then we get
\begin{align*}
\int_\R S_{2j}
w'(\tX_j)\,dx&=-(2e'_{j}f_j'+e_jf_j'')\int_\R Z'(s)w'(s)\,ds\\
&\qquad +\sum_{i\neq
j}\big\{e_i''-(2e_if_i'a_{i1}+e_if_i''a_{i1})\\
&\qquad\qquad+e_i\big[1+a_{i1}^2\big((f'_i)^2
-\beta_i^2\eta^2\big)\big]\big\}O(e^{\,-|f_i-f_j|})\\
&\qquad+ \sum_{i\neq j}(e_i'+e_if'_i)O(\alpha^3)
e^{\,-\theta_0\alpha|z|}.
\end{align*}
Using the sets $A_l$, $l=1,\dots,k$ introduced above  the
projection of $S_{3j}[\WW]$ on $w'(\tX_j)$ can be written  as
follows
\begin{align*}
\int_\R S_{3j}[\WW] w'(\tX_j)\,dx&=\int_{A_j}
\Big[f'(w_{\delta_j,\beta_j})\sum_{i\neq j}
w_{\delta_i,\beta_i}-\sum_{i\neq j}
f(w_{\delta_i,\beta_i})\Big]w'(\tX_j)\,dx\\
&\qquad +\sum_{l\neq j}\int_{A_l}
\Big[f'(w_{\delta_l,\beta_l})\sum_{i\neq l}
w_{\delta_i,\beta_i}-\sum_{i\neq l}
f(w_{\delta_i,\beta_i})\Big]w'(\tX_j)\,dx\\
 &\qquad +\sum_{l=1}^k\int_{A_l} f'(w_{\delta_l,\beta_l})
\sum_{i=1}^k e_iZ(\tX_i) w'(\tX_j)\,dx\\
&\qquad + \sum_{l=1}^k\int_{A_l} O\big((\WW-w_{\delta_l,\beta_l})^2\big)w'(\tX_j)\,dx \\
 &:= S_{3j1}+S_{3j2}+S_{3j3}+S_{3j4}.
\end{align*}
For $S_{3j1}$ we  get
$$
S_{3j1}=\sum_{i\neq j}\int_{A_j}
[pw_{\delta_j,\beta_j}^{p-1}w_{\delta_i,\beta_i}-w_{\delta_i,\beta_i}^p]w'(\tX_j)\,dx.$$
For a fixed $i \neq j$, using formulas (\ref{wdelta}), we have
\begin{align*}
&\int_{A_j}
[pw_{\delta_j,\beta_j}^{p-1}w_{\delta_i,\beta_i}-w_{\delta_i,\beta_i}^p]w'(\tX_j)\,dx\\
&=\int_{A_j} [pw(\tX_j)^{p-1}w(\tX_i)-w^p(\tX_i)]w'(\tX_j)\,dx\\
&\qquad +\delta_j\int_{A_j}
[O(e^{\,-(1+\mu)|\tX_j|})w(\tX_i)+O(e^{\,-(1+\mu)|\tX_i|})w(\tX_j)]\,dx\\
&=I_{1i}+I_{2i}.
\end{align*}
To compute $I_{1i}$ let us observe that for any $\lambda\in \R$
and $a<b$  we have:
\begin{align*}
&\int_a^b pw^{p-1}(s)w'(s) w(s-\lambda)-\int_\R w^p(s-\lambda) w'(s)=pw^p(s)w(s-\lambda)\left|_{a}^b\right.
\\
&\qquad\qquad-\int_a^b[w^p(s)w'(s-\lambda)+w^p(s-\lambda) w'(s)]\\
&=\big[pw^p(s)w(s-\lambda)+w'(s)w'(s-\lambda)-w(s)w(s-\lambda)\big]\left|_{a}^b\right.,
\end{align*}
hence, changing variables, we get, with
$s_j=\frac{f_{j-1}-f_j}{2}$, $t_j=\frac{f_{j+1}-f_j}{2}$,
$\lambda_{ij}=f_i-f_j$,
\begin{align*}
I_{1i}&=\int_{s_j}^{t_j}
pw^{p-1}(s)w'(s)w\big(sa_{j0}a_{i1}+(f_j-f_i)a_{i1}\big)\,ds\\
&\qquad-\int_\R
w^p\big(sa_{j0}a_{i1}+(f_j-f_i)a_{i1}\big)w'(s)\,ds\\
&=\big[pw^p(s)w(s-\lambda_{ij})+w'(s)w'(s-\lambda_{ij})-w(s)w(s-\lambda_{ij})\big]\left|_{s_j}^{t_j}\right.\\
&\qquad
+O(\beta_i^2+\beta_j^2)e^{\,-|f_j-f_i|}.\\
\end{align*}
 Similarly,
$$
I_{2i}=O(\delta_j)e^{\,-|f_j-f_i|}.
$$
Using this formula with $i=1,\dots, k$ and noting that for a fixed
$j$ the dominating terms in $S_{3j1}$ come from $i=j-1$ and
$i=j+1$ we get
\begin{align*}
S_{3j1}&=w\big(\frac{f_j-f_{j-1}}{2}\big)^2+w'\big(\frac{f_j-f_{j-1}}{2}\big)^2
-w\big(\frac{f_j-f_{j+1}}{2}\big)^2-w'\big(\frac{f_j-f_{j+1}}{2}\big)^2\\
&\qquad+O(\alpha)\sum_{i\neq
j } e^{\,-|f_j-f_i|}\\
&=C_p
\big(e^{\,-|f_j-f_{j-1}|}-e^{\,-|f_{j+1}-f_{j}|}\big)+O(\alpha)\sum_{i\neq
j } e^{\,-|f_j-f_i|},
\end{align*}
where we have denoted
\begin{align}
C_p=\lim_{s\to \infty}\{e^{\,s}[w(s/2)^2+w'(s/2)^2]\}.
\label{defcp}
\end{align}

The remaining term $S_{3j2}$ can be estimated in a  similar way.
Observing that when $l\neq j$ and $x\in A_l$ we have
$w'(\tX_j)=w'(\tX_l) O(e^{\,-|f_j-f_l|})$, we get that
$$
S_{3j2}=O(\alpha)\sum_{i\neq j } e^{\,-|f_j-f_i|}.$$  To compute
$S_{3j3}$ we observe that
$$
\int_\R f'(s) Z(s) w'(s)\,ds =0,
$$
hence
\begin{align*}
S_{3j3}&=\int_{A_j} f'(w(\tX_j))\sum_{i\neq j} e_iZ(\tX_i) w'(\tX_j)\,dx+e_jO(\delta_j)\\
&\qquad+\sum_{l\neq j}^k\int_{A_l} f'(w_{\delta_l,\beta_l})\sum_{i=1}^k e_iZ(\tX_i) w'(\tX_j)\,dx\\
&=O(e_j\delta_j)+\sum_{i\neq j} e_iO(e^{\,-|f_i-f_j|}).
\end{align*}
Finally we have
$$
S_{4j1}=O(e^2_j)+\sum_{i\neq j} e_i^2 O(e^{\,-|f_i-f_j|}).
$$
Summarizing we get the following expression
\begin{align}
\begin{aligned}
\Pi_{f_j}&=-f''_j[c_0+c_1e_j+O(\delta_j)]+2\sqrt{\lambda_1}a_{j0}\delta_j(f_j'-\beta_j\eta)\sin(\sqrt{\lambda_1}a_{i0}z)\\
&\qquad
+C_p(e^{\,-|f_j-f_{j-1}|}-e^{\,-|f_j-f_{j+1}|})\\
&\qquad-2c_1f'_je'_j+\gamma_{j}(f_j'-\beta_j\eta)
+O(e_j^2)+e_jO(\delta_j)
\\
&\qquad +\sum_{i\neq
j}[f_i''p_{1ij}+(f_i'-\beta_i\eta)p_{2ij}]+\sum_{i\neq
j}[f_i'e_iq_{1ij}+(e_i''+e_i'+e_i+e_i^2)q_{2ij}]\\
&\qquad+\sum_{i\neq j} \big(e_i +O(\alpha)\big)e^{\,-|f_j-f_{i}|},
\end{aligned}
\label{proj fj}
\end{align}
where functions $\gamma_j$, $p_{lij}$ and  $q_{lij}$ are smooth
functions of their arguments satisfying estimates
\begin{align*}
\gamma_j&=\gamma_j(f_j',\beta_j,\delta_j,\alpha, z)=O(\delta_j^2)+O(\alpha^2),\\
p_{lij}&=p_{lij}(f_j, f_i, e_i,\alpha, \delta,
z)=O\big(e^{\,-|f_j-f_{i}|}(1+e_i)\big),
\\
q_{lij}&=q_{lij}(f_j, f_i, e_i,\alpha, \delta,
z)=O(\delta_j)+O(e^{\,-|f_j-f_{i}|})+O(\alpha^3
e^{-\theta_0\alpha|z|}),
\end{align*}
and $C_p$, $c_0$, $c_1$ are constants defined by
$$
C_p=\lim_{t\to\infty}e^{\,t}\big(w(t/2)^2+w'(t/2)^2\big), \quad
c_0=\int_\R w'(s)^2\,ds, \quad c_1=\int_\R w'(s)Z'(s)\,ds.
$$

Now we compute projection of the error on $Z(\tX_j)$ denoted by
$\Pi_{e_j}$ above. We will use the expression for $S[\WW]$,
(\ref{expsww}). We will denote the entries of the matrix ${\mathbf
A}_{0i}$ by $A_{0i,lm}$. First we observe that
\begin{align}\begin{aligned}
\int_\R Tr\big({\mathbf A}_i\cdot
Hess_{\tX_i,\tZ_i}(w_{\delta_i,\beta_i})\big)Z(\tX_i)&=A_{0i,11}a_{i0}\int_\R
w''(s) Z(s)\,ds\\
&+\delta_i \Big[
A_{0i,11}\cos(\sqrt{\lambda_1} a_{i0} z)\int_\R Z''(s)Z(s)\,ds\\
 &
-2A_{0i,12}a_{i0}\sqrt{\lambda_1}\sin(\sqrt{\lambda_1} a_{i0}
z)\int_\R\big(Z'(s)\big)^2\,ds\\
&-\lambda_1 A_{0i,22} a_{i0}\cos(\sqrt{\lambda_1}
a_{i0} z)\int_\R Z^2(s)\,ds
\\ &\qquad +\delta_iO\big(\|f_i'-\beta_i\eta\|_{\theta_0\alpha}+\alpha^2\big)e^{\,-\theta_0\alpha|z|}\Big]\\
&\equiv h_{0i}.\end{aligned} \label{def h0j}
\end{align}
It is convenient to write
$$
{\mathbf l}_i=-f''_ia_{i1}(1,\beta_i\eta)+\mathbf{\tilde l}_i.
$$
Notice that
$$
\|\mathbf{\tilde l}_i\|_\alpha+\alpha^{-1}\|\mathbf{\tilde
l}_i'\|_\alpha\leq
C\alpha^3(1+\|f_i'-\beta_i\eta\|_\alpha+|\tX_i|).
$$
Then we have
\begin{align*}
\int_\R{\mathbf l_i}\cdot\nabla_{\tX_i,
\tZ_i}w_{\delta_i,\beta_i}Z(\tX_i)&=\delta_i
f_i''a_{i1}\beta_i\eta\sin(\sqrt{\lambda_1} a_{i0}z)\Big[\int_\R
Z^2+O(\alpha)\Big]\\
&\qquad +
O(\alpha^3)(1+\|f_i'-\beta_i\eta\|_\alpha+\delta_i)e^{\,-\theta_0\alpha|z|}.
\end{align*}

When $i\neq j$ we get
\begin{align*}& \int_\R Tr\big({\mathbf A}_i\cdot
Hess_{\tX_i,\tZ_i}(w_{\delta_i,\beta_i})\big)Z(\tX_j)\\&
=O\big(\|f_i'-\beta_i\eta\|_\alpha+\alpha^2\big)e^{\,-\theta_0\alpha|z|}
(1+\delta_i)e^{\,-|f_i-f_j|},
\end{align*}
and
$$
\int_\R{\mathbf l_i}\cdot\nabla_{\tX_i,
\tZ_i}w_{\delta_i,\beta_i}Z(\tX_j)=O\big(\|f''_i\|_\alpha+\|f_i'-\beta_i\eta\|_\alpha+\alpha^2\big)e^{\,-\theta_0\alpha|z|}(1+\delta_i)e^{\,-|f_i-f_j|}.
$$

Now we compute the projection of $S_{2j}[\WW]$. We notice first
that
\begin{align*}
\int_\R \Delta_{x,z}\big(e_jZ(\tX_j)\big)
Z(\tX_j)\,dx&=e''_jd_0+\big[2e'_ja_{j2}+e_je^{\,-\theta_0\alpha|z|}O(\alpha^3)\big]d_1
\\
&\qquad+e_j\big[1+a_{j1}\big((f_j')^2-\beta_j^2\eta^2\big)\big]d_2,
\end{align*}
where
$$
d_0=\int_\R Z^2(s)\,ds, \quad d_1=\int_\R sZ'(s)Z(s)\,ds, \quad
d_2=\int_\R Z''(s)Z(s)\,ds.
$$
On the other hand
\begin{align*}
&\sum_{i\neq j}\int_\R \Delta_{x,z}\big(e_iZ(\tX_i)\big)
Z(\tX_j)\,dx\\&=\sum_{i\neq j}
e_i''O(e^{\,(1+\mu)|f_i-f_j|})+\sum_{i\neq j}
e_i'\big[f_i'+\alpha^3O(|f_i-f_j|)\big]O(e^{\,(1+\mu)|f_i-f_j|})\\
&\qquad + \sum_{i\neq j} e_i\big[1+
f_i''a_{i1}+\big((f'_i)^2-\beta_i^2\eta^2\big)\big]O(e^{\,(1+\mu)|f_i-f_j|})\\
&\qquad +\sum_{i\neq j }e_i\alpha^3\big(O(|f_i-f_j|)+f_i'\big)
O(e^{\,-\theta_0\alpha|z|})O(e^{\,(1+\mu)|f_i-f_j|}).
\end{align*}
To calculate the projection of $S_{3j}[\WW]$ we observe that
$$
\int_\R S_{3j}[\WW] Z(\tX_j)\,dx=\sum_{l=1}^k\int_{A_l}
f'(w_{\delta_l,\beta_l}) \sum_{i=1}^k e_iZ(\tX_i) Z(\tX_j)\,dx +
\sum_{i\neq j} O(e^{\,-|f_i-f_j|}).
$$
We have
$$
\int_{A_j} f'(w_{\delta_j,\beta_j})
e_jZ^2(\tX_j)\,dx=e_j\big(d_3+O(\delta_j+\alpha^3)\big), \quad
d_3=\int_\R f'\big(w(s)\big) Z^2(s)\,ds.
$$
On the other hand when $l\neq j$ then
$$
\int_{A_l}e_if'(w_{\delta_l,\beta_l})Z(\tX_i)Z(\tX_j)\,dx=\left\{\begin{array}{ll}
e_je^{\,-(1+\mu)|f_j-f_l|}, &\quad i=j,
\medskip\\
e_ie^{\,-(1+\mu)|f_j-f_i|}, &\quad i\neq j.
\end{array}
\right.
$$
Summarizing we get
\begin{align}
\Pi_{e_j}=(e_j''+\lambda e_j)d_0+e_j g_{0j}+e_j'
g_{1j}+\sum_{i\neq j} (e_i r_{0ij}+e_i' r_{1ij}+e''_i r_{2ij})+h_j,
\label{proj ej} \end{align}
where functions $g_{lj}$, $r_{lij}$
are smooth functions of their arguments that satisfy
\begin{align*}
g_{0j}&=g_{0j}(f_j',\beta_j,
\delta_j,z)=O(\delta_j+\alpha^3)+(f_j'-\beta_j\eta)O(|f_j'|+\beta_j),\\
g_{1j}&=g_{1j}(z)=O(\alpha^3) e^{\,-\theta_0\alpha|z|},\\
r_{lij}&=r_{lij}(f_j, f_i,
f_i',f_i'',\beta_i,z)=[\alpha^3O(|f_i-f_j|)+O(|f'_i|+|f_i''|+1)]e^{\,-(1+\mu)|f_i-f_j|},\\
\end{align*}
Function $h_j$ satisfies
\begin{align*}
h_j&=h_{0j} + \delta_j f_j''a_{j1}\beta_i\eta\sin(\sqrt{\lambda_1}
a_{j0}z)\Big[\int_\R Z^2
+O(\alpha)\Big]\\
&\qquad+O\big((\alpha^3(1+\|f_i'-\beta_i\eta\|_\alpha)\big)(1+\delta_j) \\
&\qquad+\sum_{i\neq
j}O\big(\|f''_i\|_\alpha+\|f_i'-\beta_i\eta\|_\alpha+\alpha^2\big)e^{\,-\theta_0|z|}(1+\delta_i)O(e^{\,-|f_i-f_j|}),
\end{align*}
with $h_{0j}$ defined in (\ref{def h0j}).

\setcounter{equation}{0}
\section{Derivation of the reduced problem}

In this section we will derive conditions for ${\mathbf f}$,
${\mathbf e}$ and $\delta$ which imply that
\begin{equation}
c_j(z)\equiv 0, \qquad d_j(z)\equiv 0, \qquad j=1,\dots, k,
\label{red 1} \end{equation} and consequently lead to a solution
of our original problem. We recall that with
$$
\rho_j(x,z)=\eta_a^b\Big( \frac{|\tX_j|}{d_*}\Big), \quad
a=\frac{2^5-1}{2^5}, b = \frac{2^6-1}{2^6},
$$
we have
$$
\int_\R \phi(x,z) w_{j,x}(x,z)\rho_j(x,z)\,dx = 0= \int_\R
\phi(x,z) Z_j(x,z)\rho_j(x,z)\,dx, \quad\hbox{for all }z\in \R,
$$ where $\phi$ is  the solution of
\begin{equation}
\LL(\phi) = E - N(\phi) +\sum_{j=1}^k  c_j(z) \eta_jw_{j,x} +
d_j(z) \eta_j Z_j  \quad\hbox{in }\R^2, \label{red1}
\end{equation}
(see (\ref{nl1}).

To derive the reduced problem we will multiply (\ref{red1}) by
$\rho_j w_{j,x}$ and $\rho_j Z_j$, $j=1, \dots, k$ and integrate
over $\R$ with respect to $x$. To begin with let us observe that
$$
\int_\R w^2_{j,x}\eta_j\rho_j\,dx=\delta_i^j\int_\R
(w')^2+O(\alpha), \qquad \int_\R
Z^2_{j}\eta_j\rho_j\,dx=\delta_i^j\int_\R Z^2+O(\alpha),
$$
hence for (\ref{red1}) to be satisfied we need for $j=1, \dots,
k$:
\begin{align}
\int_\R [E-N(\phi)-{\mathcal L}(\phi)]w_{j,x} \rho_j\,dx&\equiv
0,\label{red2}\\
\int_\R [E-N(\phi)-{\mathcal L}(\phi)] Z_j\rho_j\,dx&\equiv
0,\label{red3}
\end{align}
where $\phi$ is the solution found in Proposition \ref{proposition
3}. These equations can be written also in the form
\begin{align}
\Pi_{f_j}&=\int_\R [N(\phi)+{\mathcal L}(\phi)] w_{j,x}\rho_j\,dx +\int_\R E[w'(\tX_j)-w_{j,x}\rho_j]\,dx\equiv Q_j({\bf f}, {\bf e}),\label{red3a}\\
\Pi_{e_j}&=\int_\R [N(\phi)+{\mathcal L}(\phi)]
Z_j\rho_j\,dx+\int_\R E[Z(\tX_j)-Z_{j}\rho_j]\,dx\equiv P_j({\bf
f}, {\bf e}), \label{red3b}
\end{align}
hence, by the results of the previous section, it only remains to
compute the right hand sides of the above expressions.

We observe first that $|N(\phi)|\leq C|\phi|^2$ hence we have
\begin{align}
\left|\int_\R N(\phi) w_{j,x}\rho_j\,dx\right|+ \left|\int_\R
N(\phi) Z_j\rho_j\,dx\right|\leq C\|\phi\|^2_{\sigma, \alpha, *}
e^{\,-\alpha|z|}\leq
C\alpha^{4-4\sigma'}e^{\,-2\theta_0\alpha|z|}, \quad
 \label{red4}
\end{align}
$\sigma'\in(\sigma, 3/2\sigma)$.
To estimate the term involving ${\mathcal L}(\phi)$ we use
integration by parts and the orthogonality condition
(\ref{phiort2}) to get
\begin{align*}
\int_\R {\mathcal L}(\phi) w_{j,x}\rho_j\,dx&=\int_\R \phi[2
w_{j,xx}\rho_{j,x}+w_{j,x}\rho_{j,xx}]\,dx\\
&\qquad
-\int_\R[\phi_z(w_{j,x}\rho_j)_z+\phi(w_{j,x}\rho_j)_{zz}]\,dx\\
&\qquad +\int_\R [f'(\WW)-f'(w_j)]\phi w_{j,x}\rho_j\,dx
\\
&\qquad = I + II + III.
\end{align*}
To estimate the first term above we notice that in the support of
$\rho_{j,x}$ and $\rho_{j,xx}$ we have
\begin{align}
\frac{(2^6-1)d^*}{2^6}>|x-f_j|>\frac{(2^5-1)d^*}{2^5},
\label{supprho1}
\end{align}
hence
$$
|I|\leq C\|\phi\|_{\sigma, \alpha, *} e^{\,-\theta_0{\alpha}|z|}
\alpha^{15/16}\leq
C\alpha^{2-2\sigma'+15/16}e^{\,-\theta_0{\alpha}|z|}.
$$
The second term is estimated similarly, using additionally
$|f'_j|=O(\alpha)$, so that
$$
|II|\leq C\alpha^{2-2\sigma'+15/16}e^{\,-\theta_0{\alpha}|z|}.
$$
To estimate the last term we notice that
$$
|[f'(\WW)-f'(w_j)]\rho_j|\leq
C(\alpha+\|\vec{\delta}\|+\|e\|_{\theta_0\alpha}),
$$
hence
$$
|III|\leq
C\alpha^{2-2\sigma'}(\alpha+\|\vec{\delta}\|+\|e\|_{\theta_0\alpha}).
$$

To estimate the last term on the right hand side of the equation
(\ref{red3a}) we observe that
\begin{align*}
|w'(\tX_j)-w_{j,x}\rho_j|&\leq
|w'(\tX_j)-w'(x-f_j)|\rho_j+(1-\rho_j)w'(\tX_j)\\
&\leq C\alpha^2 |w'(x-f_j)|\rho_j+C(1-\rho_j)\alpha^{31/32},
\end{align*}
hence
\begin{align*}
\left|\int_\R E[w'(\tX_j)-w_{j,x}\rho_j]\,dx\right|\leq
C\|E\|_{\sigma, \alpha, *}\alpha^{31/32}
e^{\,-\theta_0\alpha|z|}\leq
C\alpha^{2-2\sigma'+15/16}e^{\,-\theta_0\alpha|z|}.
\end{align*}
Analogous argument applies for the components of the right hand
side of (\ref{red3b}). Taking now $\sigma<2^{-10}$ we get that
there exists a $\mu>7/8$ such that
\begin{align}
|Q_j({\bf f}, {\bf e}, \vec{\delta})|\leq C\alpha^{2+\mu}, \qquad
|P_j({\bf f}, {\bf e},\vec{\delta})|\leq C\alpha^{2+\mu}.
\label{red5}
\end{align}
Now we will consider Lipschitz dependence of $Q_j$ and $P_j$ on
their parameters. To do this we can argue in a way similar to the
way we proved Lemma \ref{lem lipsch} and Proposition
\ref{proposition 3a}. This means we need to introduce parameters
${\bf g}_i$, ${\bf h}_i$, ${\bf \tilde f}_i$, ${\bf \tilde e}_i$,
${\bf e}_i$, ${\bf f}_i$ and $\vec{\delta}_i$ and consider
Lipschitz dependence with respect to them making use of Lemma
\ref{lem lipsch} and Proposition \ref{proposition 3a}. In fact,
when ${\bf f}$ is fixed then by a similar argument as above we get
(using the familiar notation)
\begin{align}
\begin{aligned}
\|Q_j^{(1)}-Q_j^{(2)}\|_{\theta_0\alpha}+\|P_j^{(1)}-P_j^{(2)}\|_{\theta_0\alpha}&\leq
C\alpha^\mu(\|{\bf g}_1''-{\bf g}_2''\|_{\theta_0\alpha}+\|{\bf
h}_1''-{\bf
h}_2''\|_{\theta_0\alpha})\\
&\qquad +C\alpha^\mu(\alpha\|{\bf \tilde f}_1'-{\bf \tilde
f}_2'\|_{\theta_0\alpha}+\|{\bf \tilde e}_1'-{\bf \tilde e}_2'\|_{\theta_0\alpha})\\
&\qquad +C\alpha^\mu(\|{\bf e}_1-{\bf
e}_2\|_{\theta_0\alpha}+\alpha^2\|\vec{\delta}_1-\vec{\delta}_2\|).
\end{aligned}
\label{red6}
\end{align}
When all parameters except ${\bf f}$ are fixed then the argument
is similar. Let us for example consider one typical term in the
projection $Q_j$:
\begin{align*}
\left|\int_R
[\phi^{(1)}w^{(1)}_{j,xx}\rho^{(1)}_{j,x}-\phi^{(2)}w^{(2)}_{j,xx}\rho^{(2)}_{j,x}]\,dx\right|&\leq
\int_\R|\phi^{(1)}-\phi^{(2)}||w^{(1)}_{j,xx}\rho^{(1)}_{j,x}|\,dx\\
&\qquad +
\int_\R|\phi^{(2)}||\phi^{(1)}w^{(1)}_{j,xx}\rho^{(2)}_{j,x}-w^{(2)}_{j,xx}\rho^{(2)}_{j,x}|\,dx\\
& \leq C\alpha^{31/32}
e^{\,-\theta_0\alpha|z|}(\|\phi^{(1)}-\phi^{(2)}\|_{\sigma,
\theta_0\alpha, {\bf f}_1}\\
&\qquad+\|{\bf f}_1-{\bf f}_2\|_{\infty}
\|\phi^{(2)}\|_{\sigma, \theta_0\alpha, {\bf f}_2})\\
&\leq C\alpha^{2+\mu}e^{\,-\theta_0\alpha|z|}\|{\bf f}_1-{\bf
f}_2\|_{\infty}.
\end{align*}
Other terms can be handled in a similar way. In summary we will
get then, taking into account the Lipschitz dependence on ${\bf
f}$ only:
\begin{align}
\begin{aligned}
\|Q_j^{(1)}-Q_j^{(2)}\|_{\theta_0\alpha}+\|P_j^{(1)}-P_j^{(2)}\|_{\theta_0\alpha}&\leq
C\alpha^{2+\mu}e^{\,-\theta_0\alpha|z|}\|{\bf f}_1-{\bf
f}_2\|_{\infty}. \label{red7}
\end{aligned}
\end{align}

We will now consider the right hand sides of (\ref{red3a}) and
(\ref{red3b}). Let us go back to formulas (\ref{proj fj}) and
(\ref{proj ej}). Using the same notation we will write now
(\ref{red3a}) and (\ref{red3b}) in the form:
\begin{align}
\begin{aligned}
&-c_0f''_j+2\sqrt{\lambda_1}a_{j0}\delta_j(f_j'-\beta_j\eta)\sin(\sqrt{\lambda_1}a_{i0}z)
+C_p(e^{\,-|f_j-f_{j-1}|}-e^{\,-|f_j-f_{j+1}|})\\
&\qquad=\widetilde\Pi_{f_j}+Q_j,\end{aligned}
\label{red8}
\\
(e_j''+\lambda_1 e_j)d_0+h_{0j}=\widetilde\Pi_{e_j}+P_j,
\label{red9}
\end{align}
where we have
$$
\widetilde\Pi_{f_j}=\widetilde\Pi_{f_j}({\bf f}'', {\bf e}'', {\bf
f}', {\bf e}', {\bf f}, \vec{\delta}), \qquad
\widetilde\Pi_{e_j}=\widetilde\Pi_{e_j}({\bf f}'', {\bf e}'', {\bf
f}', {\bf e}', {\bf f}, \vec{\delta}).
$$
Functions $\widetilde\Pi_{f_j}$ and $\widetilde\Pi_{e_j}$ are
smooth functions of their parameters and under the assumptions
(\ref{co1})--(\ref{cobeta}) and (\ref{co2a})--(\ref{co5}) we have:
\begin{align}
\|\widetilde\Pi_{f_j}\|_{\theta_0\alpha}+\|\widetilde\Pi_{e_j}\|_{\theta_0\alpha}\leq
C\alpha^3.\label{red10}
\end{align}
The rest of this paper will be devoted to solving system
(\ref{red8})--(\ref{red9}).

\setcounter{equation}{0}

\section{The  Toda system}

In this section we will determine the main order behavior of
function ${\bf f}$. Intuitively it should be given by solving the
following system
\begin{align}
-c_0f''_j +C_p(e^{\,-|f_j-f_{j-1}|}-e^{\,-|f_j-f_{j+1}|})=0,\quad
j=1,\dots,k. \label{tod1}
\end{align}

Let us consider first the case $k=2$. Remembering that in this
case $f_3=\infty$ and $f_0=-\infty$ we can reduce (\ref{tod1}) to
a single equation for $u=f_1-f_2$:
\begin{align}
u'' + 2c_p e^{-u} = 0, \quad u(0)=a_1-a_2\quad c_p = \frac {C_p}{\int_\R w_x^2},
\quad\mbox{where $C_p$ is defined in (\ref{defcp})}, \label{tod2}
\end{align}
and $f_j(0)=a_j$, see (\ref{k2}).
For  $\alpha>0$  a family (parametrized by $\alpha$) of explicit   even solutions of (\ref{tod2}) is  given by:
\begin{align}
u_\alpha (z):= 2\log \frac 1\alpha +  \log {2\lambda^{-2}c_p} -\log
\Big(\frac{-1+\tanh^2(\lambda\alpha z/2)}{2}\Big) , \quad \alpha
>0,
\label{tod3}
\end{align}
where $\lambda=\sqrt{\frac{2c_p}{e^{\,a_2-a_1}}}$.
Starting with this family it is natural to define
\begin{align}
f_1(z)=-\frac{1}{2}u_\alpha(z),\quad f_2(z)=\frac{1}{2}
u_\alpha(z). \label{tod4}
\end{align}
Let us observe that in particular  $f_1(z)+f_2(z)\equiv 0$ and
\begin{align}
\| f''_i\|_\alpha=O(\alpha^2), \quad \beta_1=-\frac{\lambda}{2}\alpha,
\quad\beta_2=\frac{\lambda}{2}\alpha, \quad f_2-f_1\geq
2\log\frac{1}{\alpha}+ \log {2\lambda^2c_p}, \label{tod5}
\end{align}
 as needed. In this case we take $\theta_0=\frac{\lambda}{4}$ and $\theta_1=\frac{\lambda}{2}$.
 In the sequel we
 will denote ${\bf f}_0(z)=(f_1(z), f_2(z))$.

We will assume now $k>2$ since the case $k=2$ has just been
treated above. It is convenient to  consider our problem in a
slightly more general framework then that of the system
(\ref{tod1}).
Thus for given functions $q_j(t), p_j(t)$, $j=1,\dots, k$ such
that
$$
\sum_{j=1}^{k} q_j=\sum_{j=1}^k p_j=0,
$$
 we
define the Hamiltonian
$$
H=\sum_{j=1}^k \frac{p_j^2}{2}+V, \quad V=\sum_{j=1}^{k-1}
e^{\,(q_j-q_{j+1})}.
$$
We consider the following Toda system
\begin{align}
\begin{aligned}
\frac{dq_j}{dt}&=p_j,\\
\frac{dp_j}{dt}&=-\pd{H}{q_j}\\
q_j(0)&=q_{0j},\quad p_j(0)=0\quad  \,\quad j=1, \dots, k .
\end{aligned}\lbl{ts3}
\end{align}
Solutions to (\ref{ts3}) are of course even. Observe that also
that the location of their center of mass remains fixed. Thus to
mode out translations we will assume that
\begin{align}
\sum_{j=1}^k q_{0j}=0. \label{tsmass}
\end{align}
We will now give a more precise of these solutions and in
particular their asymptotic behavior as $z\to\pm\infty$. To this
end   we will often make use of classical results of Konstant
\cite{konst} and in particular we will use the explicit formula
for the solutions of (\ref{ts3}) (see formula (7.7.10) in
\cite{konst}).

We will first introduce some notation. Given numbers $w_1, \dots,
w_{k}\in \R$ such that
\begin{align}
\sum_{j=1}^{k} w_j=0, \quad\mbox{and}\ w_j>w_{j+1},\quad j=1,
\dots, k \label{tsw0} \end{align}
we define the matrix
$$
{\bf w}_0=\mbox{diag}\,(w_1, \dots, w_k).
$$
Furthermore, given numbers $g_1, \dots, g_k\in \R$ such that
\begin{align}
\prod_{j=1}^k g_j=1, \quad \mbox{and}\ g_j>0, \quad j=1, \dots, k,
\label{tsg0}
\end{align}
we define the matrix
$$
{\bf g}_0=\mbox{diag}\,(g_1, \dots, g_N).
$$
The matrices $ {\bf w}_0$ and ${\bf g}_0$ can be parameterized by
introducing the following two sets of parameters
\begin{equation}
c_j= w_j -w_{j+1}, \quad   d_j= \log g_{j+1} - \log g_{j},\quad
j=1,\ldots, k.
\end{equation}
Furthermore, we  define functions $\Phi_j({\bf g}_0, {\bf
w}_0;t)$, $t\in \R$, $j=0, \dots, k$, by
\begin{eqnarray}
& &\Phi_0=\Phi_k\equiv 1\nonumber\\
& & \Phi_j ({\bf g}_0, {\bf w}_0;t)=\ \lbl{ts33}\\
& & (-1)^{j(k-j)}\sum_{1\leq i_i<\dots<i_j\leq k} r_{i_1\dots
i_j}({\bf w}_0)g_{i_1}\dots g_{i_j}\exp[-t(w_{i_1}+\dots+w_{i_j})],
\nonumber
\end{eqnarray}
where $r_{i_1\dots i_j}({\bf w}_0)$ are rational functions of the
entries of the matrix ${\bf w}_0$. It is proven in \cite{konst}
that all solutions of (\ref{ts3}) are of the form
\begin{equation}
q_j(t)=\log\Phi_{j-1}({\bf g}_0, {\bf w}_0;t)-\log\Phi_j({\bf
g}_0, {\bf
 w}_0;t),\quad  j=1,
\dots, k, \lbl{ts33a} \end{equation} Namely, given initial
conditions in (\ref{ts3}) there exist matrices ${\bf w}_0$ and
${\bf g}_0$ satisfying (\ref{tsw0})--(\ref{tsg0}).
According to
Theorem 7.7.2 of \cite{konst}, it holds
\begin{equation}
\label{ts55} q_j^{'} (+\infty)= w_{k+1-j}, \quad q_j^{'}(-\infty)=
w_{j},\quad j=1, ..., k.
\end{equation}

We introduce variables
\begin{equation}
u_j=q_{j}-q_{j+1}.
\end{equation}
 In terms of ${\bf u}=(u_1, \dots, u_{k-1})$
system (\ref{ts3}) becomes
\begin{align}
\begin{aligned}
{\bf u}''-Me^{\,{\bf u}}&=0,\\
u_j(0)=q_{0j}-q_{0 j+1}&,\quad u'_j(0)=0, \quad j=1, \dots, k-1,
\end{aligned}\lbl{ts4}
\end{align}
where
$$
M=\left(\begin{array}{cccc} 2 &
-1& 0\cdots&0\\
-1 & 2& -1 \cdots &0\\
& & \ddots& \\
0& \cdots& 2& -1\\
0& \cdots& -1 & 2
\end{array}
\right),\quad e^{\,-{\bf u}}=\left(\begin{array}{c} e^{\,u_1}\\
\vdots\\
e^{\,u_{k-1}}
\end{array}
\right).
$$
As a consequence of (\ref{ts33}) all solutions to (\ref{ts4}) are
given by
\begin{eqnarray}
u_j(t)=q_{j}(t)-q_{j+1}(t)&=&-2\log\Phi_j({\bf g}_0,
{\bf w}_0;t)+\log\Phi_{j-1}({\bf g}_0, {\bf w}_0;t)\nonumber\\
&& \quad+\log\Phi_{j+1}({\bf g}_0, {\bf w}_0;t). \lbl{ts44a}
\end{eqnarray}
Our first goal is to prove the following:
\begin{lemma}\lbl{lem toda 1}
Let $ {\bf w}_0$ be such that
\begin{equation}
\label{wchoice} \min_{j=1,\dots,k-1} (w_j -
w_{j+1})=\vartheta>0.
\end{equation}
Then there holds
\begin{equation}
u_j(t)=\left\{\begin{array}{l} -  c_{k-j} t - d_{k-j} + \tau_j^{+}
( {\bf c}) +  O(e^{\,- \vartheta{ |t|}
} ), \mbox{as} \ t \to +\infty, \quad j=1, \dots, k-1, \medskip\\
  c_{j} t + d_{j} + \tau_j^{-} ( {\bf c}) +  O(e^{\,- \vartheta{ |t|}}),
\mbox{as} \ t \to -\infty, \quad j=1, \dots, k-1,
  \end{array}\right.
\lbl{ts5}
\end{equation}
where $ \tau_j^\pm ({\bf c})$ are smooth functions of the vector
${\bf c}=(c_1, \ldots , c_{k-1})$.
\end{lemma}
\noindent{\bf Proof.} This Lemma has been proved in
\cite{dkw}. We include a proof here for completeness.  Let
$q_j$, $j=1, \dots, k$ be a solution of the system (\ref{ts3})
depending on the (matrix valued) parameters ${\bf w}_0$, ${\bf
g}_0$ and defined in (\ref{ts33a}). We  need to study the
asymptotic behavior of $\Phi_j({\bf w}_0, {\bf g}_0;t)$ as $t\to
\pm \infty$ with the entries of ${\bf w}_0$ satisfying
(\ref{wchoice}) and still undetermined ${\bf g}_0$.

\medskip
By (\ref{ts33}) and (\ref{tsw0}), we get that as $t \to -\infty$
$$
\Phi_j=(-1)^{j(k-j)} r_{1\dots j}({\bf w}_0)g_1\dots g_j e^{ -(w_1
+ \cdots + w_j) t}(1+O(e^{\,-(w_j-w_{j-1})t})\, ,
$$
hence
\begin{equation}
\frac{\Phi_{j+1}\Phi_{j-1}}{\Phi^2_j}=\frac{g_{j+1}r_{1\dots
(j-1)}({\bf w}_0)r_{1\dots (j+1)}({\bf w}_0)e^{-c_jt}}{g_j
r_{1\dots j}^2({\bf w}_0)}(1+O(e^{-\vartheta{|t|} }))\, .
\lbl{ts10}
\end{equation}
It follows that as $t\to-\infty$
\begin{eqnarray}
\lefteqn{u_j(t)=\log\left(\frac{\Phi_{j+1}\Phi_{j-1}}{\Phi^2_j}\right)}
\nonumber\\
&=&- c_j t+\log\left(\frac{g_{j+1}r_{1\dots (j-1)}({\bf
w}_0)r_{1\dots
(j+1)}({\bf w}_0)}{g_j r_{1\dots j}^2({\bf w}_0)}\right)+O(e^{-\vartheta{|t|}}) \\
&=& c_j t + d_j  + \tau_k^{-} ({\bf c}) +O(e^{-\vartheta{|t|}}),
\lbl{ts11}
\end{eqnarray}
where
\[ \tau_j^{-} ({\bf c}) =
\log\left(\frac{r_{1\dots (j-1)}({\bf w}_0)r_{1\dots (j+1)}({\bf
w}_0)}{ r_{1\dots j}^2({\bf w}_0)}\right).\] Similarly,  as $t\to
+\infty$ we get
\begin{eqnarray}
\lefteqn{u_j
(t)=\log\left(\frac{\Phi_{j+1}\Phi_{j-1}}{\Phi^2_j}\right)}
\nonumber\\
&=& -c_{k-j}t- d_{k-j} +\tau_{j}^{+} + O(
e^{-\vartheta{|t|}}),\lbl{ts111}
\end{eqnarray}
where
\[ \tau_j^{+} ({\bf c})= \log\left(\frac{ r_{k+2-j\dots k}({\bf w}_0)r_{k-j\dots
k}({\bf w}_0)}{ r_{k+1-j\dots
 k}^2({\bf w}_0)}\right).\]
This ends the proof. \qed

\medskip





To find a family of solutions parameterized by $\alpha$ starting
from a solution of (\ref{ts3}) we use functions $u_j$ and set
\begin{align}
u_{j\alpha}(z)=u_j(\alpha z)-2\log\frac{1}{\alpha}-\log c_p.
\label{tsualpha} \end{align} Then functions $f_{j\alpha}(z)$ are obtained
from the relations
\begin{align}
\begin{aligned}
u_{j\alpha}(z)=f_{j\alpha}(z)-f_{j+1\alpha}(z),\\
\sum_{j=1}^k f_{j\alpha}(z)=0.\end{aligned} \label{tsfalpha}
\end{align}
Observe that as a consequence of Lemma \ref{lem toda 1} we get
that there exist $w_j, g_j$, $j=1,\dots, k$ such that
(\ref{tsw0})--(\ref{tsg0}) holds, that
$$
\min_{j=1,\dots, k} (w_j-w_{j+1})=\vartheta>0,$$ and functions
$f_j$ satisfy
\begin{align}
\begin{aligned}
&\|f_j\|_{\vartheta\alpha}=\|f_j''
e^{\,\vartheta\alpha|z|}\|_\infty\leq C\alpha^2,
\\
&f_j'(\infty)=\beta_j=f'_j(-\infty), \quad\mbox{where}\
\beta_{j+1}-\beta_{j}=(w_{j}-w_{j+1})\alpha>\vartheta \alpha,\\
&f_j(z)-f_{j+1}(z)\geq 2\log\frac{1}{\alpha}+\log c_p.
\end{aligned}\label{tsf0}
\end{align}
In this case we take $\theta_0=\frac{1}{4}\vartheta$ and $\theta_1=\frac{1}{2}\vartheta$.

\setcounter{equation}{0}
\section{Bounded  solvability of some equations in the
line}\label{sec solvability}
In this section will continue with
preliminaries needed to solve (\ref{red8})--(\ref{red9}). First we
will study the linearization of the system (\ref{red8}) around the
solution ${\bf f}_0$ of the Toda system  defined in (\ref{tod4})
and (\ref{tsf0}). We will always assume that $\alpha>0$ is small
and $\theta_0>0$ has value defined in the previous section.

%

Again, we will consider the case $k=2$ first. Let $p_j(z)$,
$j=1,2$ be given even and continuous functions. The linearized
Toda system takes form
\begin{align}
\begin{aligned}
\varphi_1''+c_pe^{\,-u_\alpha}(\varphi_1-\varphi_2)&=p_1(z),\\
\varphi_2''-c_pe^{\,-u_\alpha}(\varphi_1-\varphi_2)&=p_2(z),\\
\varphi_j(0)=x_j, \quad \varphi'(0)&=0, \quad j=1,2, \end{aligned}
\label{lleqn0}
\end{align}
which can be reduced to a single ODE for $h=\varphi_1-\varphi_2$
with $p(z)=p_1(z)-p_2(z)$:
\begin{align}
\begin{aligned}
\LL_\alpha (h) := h'' + 2c_pe^{-u_\alpha }h  &=  p(z), \quad z\in \R, \\
h(0)=h_0, \quad h'(0)&=0. \end{aligned} \label{lleqn1}\end{align}
Notice that we do not impose the initial condition in
(\ref{lleqn0}) and (\ref{lleqn1}). In fact we will require that
the solution to (\ref{lleqn0}) and (\ref{lleqn1}) is {\it bounded
and even} and this determines uniquely the initial data.

We will assume that function $p(z)$ is even and satisfies
\begin{align}
\|p\|_{\theta_0\alpha}\leq C\alpha^{2+\mu}, \quad\mbox{some}\
\mu>0, \label{ll1aaa}
\end{align}
and look for a solution to (\ref{lleqn1}) in the space of $C^2$
even functions such that
\begin{align}
\|h''\|_{\theta_0\alpha}+\alpha\|h'\|_{\theta_0\alpha}+\alpha^2\|h\|_{\infty}<\infty.
\label{ll2}
\end{align}
We will denote the space of such functions by  ${\mathcal X}$.
Observe that evenness of solutions to (\ref{lleqn1}) is guaranteed
by the initial conditions and the evenness of the right hand side.

To solve  problem (\ref{lleqn1}) we will construct a suitably
bounded left inverse for the linearized operator $\LL_\alpha$. By
$\psi_{j\alpha}$ we will denote  the elements of the fundamental
set of $\LL_\alpha$. They are known explicitly:
$$
\psi_{1\alpha}(z)=u_\alpha'(z), \quad
\psi_{2\alpha}(z)=zu'_\alpha(z)+2,
$$
and their Wronskian is $W(\psi_{1\alpha},
\psi_{2\alpha})=\alpha^{2}$. Then  there exists a unique even and
bounded solution  solution to (\ref{lleqn1}). It is given by:
\begin{align}
\begin{aligned}
h(z)&=-\frac{1}{\alpha^2}\psi_{1\alpha}(z)\int_0^z
\psi_{2\alpha}(s)
p(s)\,ds+\frac{1}{\alpha^2}\psi_{2\alpha}(z)\int_0^z
\psi_{1\alpha}(s)
p(s)\,ds\\
&\qquad-\frac{1}{\alpha^2}\psi_{2\alpha}(z)\int_0^\infty u'_\alpha(s) p(s)\,ds\\
&=\frac{1}{\alpha^2}\int_0^z u'_\alpha(z)
u'_\alpha(s)(z-s)p(s)\,ds +\frac{1}{\alpha^2}\int_0^z \big(
u'_\alpha(s)-u'_\alpha(z)\big)p(s)\,ds\\
&\qquad-\frac{1}{\alpha^2}\psi_{2\alpha}(z)\int_0^\infty
u'_\alpha(s) p(s)\,ds.\end{aligned} \label{ll3}
\end{align}
Directly examining this formula we get that
$$
\|h''\|_{\theta_0\alpha}\leq C\|p\|_{\theta_0\alpha}, \qquad
|h(\infty)|\leq C\alpha^{-2}\|p\|_{\theta_0\alpha}\leq
C\alpha^{\mu},
$$
and using $h'(0)=0$ we obtain the following estimate for the
inverse of $\LL_\alpha$:
\begin{align}
\|h''\|_{\theta_0\alpha}+\alpha\|h'\|_{\theta_0\alpha}+\alpha^2\|h\|_{\infty}\leq
C\|p\|_{\theta_0\alpha}. \label{ll4}
\end{align}

These considerations provide us with a framework needed to solve
(\ref{lleqn1}). We have the following result.

\begin{lemma}\label{ll lema}
There exists a constants $\mu>0, C>0$ independent of $\alpha$ such
that if $p$ is an even function with
\begin{equation}
\| p\, \|_{ \theta_0\A}  \le C\alpha^{2+\mu},
\label{hh}\end{equation} then Problem \ref{lleqn1} possesses an
even solution $h$ of the form (\ref{ll3} such that the estimate
(\ref{ll4}) holds. Denoting this solution for a given $p$ by
$\RR_\alpha [p]$ we have additionally the following estimate
\begin{align*}
&\|\RR_\alpha [p_1]'' - \RR_\alpha [p_2]'' \|_{\theta_0 \A}+\alpha
\|\RR_\alpha [p_1]' - \RR_\alpha [p_2]' \|_{\theta_0
\A}+\alpha^2\|\RR_\alpha [p_1] - \RR_\alpha [p_2] \|_{\infty}\\
&\qquad\le\, C \|p_1- p_2 \, \|_{\theta_0\A},
\end{align*}
for all $p_1$, $p_2$ satisfying $\equ{hh}$.
\end{lemma}
The proof of this result is left to the reader.

%

Now we will consider the general case $k>2$. We are  lead to  the
following linear system
\begin{align}
\begin{aligned}
&\vec{\phi}''-\left(\begin{array}{cccc} 2e^{\,u_{1\alpha}} &
-e^{\,u_{2\alpha}}& 0\cdots&0\\
-e^{\,-u_{1\alpha}} & 2e^{\,u_{2\alpha}}& -e^{\,u_{3\alpha}} \cdots &0\\
& & \ddots& \\
0& \cdots& 2e^{\,u_{k-2\alpha}}& -e^{\,u_{k-1\alpha}}\\
0& \cdots& -e^{\,u_{k-2\alpha}} & 2e^{\,u_{k-1\alpha}}
\end{array}\right)\vec{\phi}^T=\vec{p}, \\
&
\quad \vec{\phi}=(\phi_1, \dots,
\phi_{k-1}), \quad \vec{p}=(p_1, \dots, p_{k-1}),\end{aligned} \lbl{ts14}
\end{align}
where $\vec{p}$ is an even function such that
\begin{align}
\|\vec{p}\|_{\theta_0\alpha}\leq C\alpha^{2+\mu}. \label{ll5}
\end{align}
We will analyze the solvability of this problem in the space of
even $C^2$ functions $\vec{\phi}$ such that
\begin{align}
\|\vec{\phi}''\|_{\theta_0\alpha}+\alpha\|\vec{\phi}'\|_{\theta_0\alpha}+\alpha^2\|\vec{\phi}\|_{\infty}<\infty.
\label{ll6}
\end{align}
Thus in addition to (\ref{ts14}) we will require that
\begin{align}
\vec{\phi}(0)=\vec{x}, \quad\vec{\phi}'(0)=0. \label{ts14a}
\end{align}

We first observe that
$$
g_j\pd{u_{m\alpha}}{g_j}=\alpha\left\{\begin{array}{rl}  1, &\quad
j=m+1,\quad
z\to\infty,\\
-1, &\quad j=m,\quad
z\to\infty\\
1, &\quad j=k+2-m,\quad
z\to -\infty,\\
-1, &\quad j=k+1-m,\quad
z\to -\infty,\\
0, &\qquad \mbox{otherwise}.
\end{array}
\right.
$$
Hence by a transformation we can find a set of linearly
independent solutions to the homogenous version of (\ref{ts14})
$$
\psi_{1j\alpha}(z)=\left\{\begin{array}{rl} \alpha\vec{e}_j, &\quad z\to\infty,\\
-\alpha\vec{e}_j, &\quad z\to-\infty.
\end{array}
\right.
$$
Similarly, considering derivatives of $u_k$ with respect to $w_j$
we can find solutions of (\ref{ts14}),  $\psi_{2j}(z)$,
$j=1,\dots, k-1$ such that
$$
\psi_{2j\alpha}(t)=\alpha\vec{e}_j|z|+O(1).
$$
The functions $\{(\psi_{1j}(z),\psi_{1j}'(z)) $, $(\psi_{2j}(z),
\psi_{2j}'(z)) \}$ form a fundamental set for the system
(\ref{ts14}), whose behavior as $z\to\pm\infty$ is analogous to
that of the functions $(\psi_1(z)\psi_1'(z))$, $(\psi_2(z),
\psi_2'(z))$, respectively. Let us denote the fundamental matrix
of the system (\ref{ts14}) described above by $\Psi_{\alpha}(z)$
and the right hand side of the transformed system by $\vec{q}$. We
observe that as $z\to \infty$, matrices $\Psi_\alpha,
\Psi^{-1}_\alpha$ are block matrices of the form
\begin{align}
\Psi_\alpha(z)=\left(\begin{array}{ll}\alpha I+o(1) & \alpha zI+O(1)\\
{\bf 0} & \alpha I+o(1)\end{array}\right), \quad \Psi_\alpha^{-1}=\left(\begin{array}{ll}\alpha I+o(1) & -\alpha zI+O(1)\\
{\bf 0} & \alpha I +o(1)\end{array}\right). \label{ll8}
\end{align}
Let us denote these blocks by $\Psi_{mn\alpha}$,
$\Psi^{-1}_{mn\alpha}$, $m, n=1,2$, respectively.  Then, from
variation of parameters formula we get that the solution of our
problem has form
\begin{align}
(\vec{\phi}(z), \vec{\phi}'(z))
=\Psi_\alpha(z)\cdot\int_0^z\Psi^{-1}_\alpha(s)\cdot({\bf 0},
\vec{q}(s))^T\,ds-\Psi_{12\alpha}(z)\int_0^\infty
\Psi_{22\alpha}^{-1}(s)\vec{q}(s)\,ds. \label{ll7}
\end{align}
Using this we can directly estimate
$$
\|\vec{\phi}''\|_{\theta_0\alpha}\leq C\|{\vec
q}\|_{\theta_0\alpha},
$$
form which it follows that (using the notation of Lemma \ref{ll
lema} and going back to the original variables)
\begin{align}
\|\RR[\vec{p}]''\|_{\theta_0\alpha}+\alpha\|\RR[\vec{p}]'\|_{\theta_0\alpha}+\alpha^2\|\RR[\vec{p}]\|_{\infty}&\leq
C\|\vec{p}\|_{\theta_0\alpha}, \label{ll9}\\
\|\RR[\vec{p}_1]''-\RR[\vec{p}_2]''\|_{\theta_0\alpha}+\alpha\|\RR[\vec{p}_1]'-\RR[\vec{p}_2]'\|_{\theta_0\alpha}
+\alpha^2\|\RR[\vec{p}_1]-\RR[\vec{p}_2]\|_{\infty}&\leq
C\|\vec{p}_1-\vec{p}_2\|_{\theta_0\alpha}. \label{ll10}
\end{align}
\medskip



A second problem, important for our purposes is that of finding an
even
solution of the equation
\begin{equation}
e''  + \la_1 e  = q(z),\quad z\in \R, \label{lleqn2}\end{equation}
where, again,  $q$ is an even function with $ \|q\|_{ \A} <+\infty
$. This time we want $e$  to satisfy $$
\|e''\|_{\theta_0\alpha}+\|e'\|_{\theta_0\alpha}+
\|e\|_{\theta_0\A} <+ \infty.$$ We need to assume the following
solvability condition on $q$ for this to be the case:
\begin{equation}
\int_0^\infty q(z) \cos (\sqrt{\la_1} z ) \, dz\, =\, 0.
\label{ortp}\end{equation} Under this assumption the solution
turns out to be unique. Denoting the solution of (\ref{lleqn2}) by
$\SS[q]$ we get  explicitly
\begin{align}
\SS[q ] = \frac 1{\sqrt{\la_1}} \sin(\sqrt{\la_1} z )
\int_z^\infty
 q(t) \cos (\sqrt{\la_1} t ) \, dt
-
\frac 1{\sqrt{\la_1}} \cos(\sqrt{\la_1} z ) \int_z^\infty
 q(t) \sin (\sqrt{\la_1} t )\, dt .
\label{ll11} \end{align}
Clearly this operator is bounded in the
sense that $e= S[q]$ satisfies
\begin{equation}
 \|  e \|_{\theta_0 \A}  +  \|  e' \|_{\theta_0 \A}  + \|  e'' \|_{\theta_0 \A}
\le C\alpha^{-1} \|
 q\, \|_{\theta_0 \A}.
\label{e1}\end{equation} On the other hand, a better estimate is
available in case that we know, in addition that $ \|
 q' \|_{\theta_0 \A}  < +\infty$.
Indeed, in this case integrating by parts in the formula
(\ref{ll11}) we get
 \begin{equation}
 \|  e \|_{\theta_0 \A}  +  \|  e' \|_{ \theta_0\A}  + \|  e'' \|_{\theta_0 \A}
\le C(\, \| q \|_{\theta_0 \A}  + \alpha^{-1} \| q' \|_{\theta_0
\A} \, ). \label{e2}\end{equation} Finally it is clear that the
operator $\SS[q]$ is Lipschitz in the norms used above.

\setcounter{equation}{0}

\section{Solving the reduced system for $(\ff, \ee, \vec{\delta})$}

We will now go back to  (\ref{red8})--(\ref{red9}). We will solve
this system using a fixed point argument around and approximate
solution ${\bf f}={\bf f}_0$, ${\bf e}=0$ and $\vec{\delta}=0$,
where ${\bf f}_0$ is a solution to Toda system  will be determined
shortly. We let ${\bf f}_0$ to be an even solution to (\ref{tsf0})
satisfying the initial data
$$
f_{0j}(0)=x_j, \quad \sum_{j=1}^k x_j=0, \qquad f_{0j}'(0)=0,
\quad j=1, \dots, k.
$$
We will set
\begin{align}
\beta_j=f_{0j}(\pm\infty) \label{srs-1}.
\end{align}
In the sequel we will often use the fact that ${\bf f}_0$ is a
smooth function which in addition satisfies
\begin{align}
\alpha^{-1}\|{\bf f}_0'''\|_{\theta_0\alpha}+\|{\bf
f}_0''\|_{\theta_0\alpha}+\alpha\|{\bf f}_0'\|_{\theta_0\alpha}
\leq C_0\alpha^2.\label{srs-05}
\end{align}

 With ${\bf f}_0$ fixed now we will will
assume that the solution of our problem is of the form
$$
{\bf f}={\bf f}_0+{\bf f}_1,
$$
where
\begin{align}
\|{\bf f}_1''\|_{\theta_0\alpha}+\alpha\|{\bf
f}_1'\|_{\theta_0\alpha}+\alpha^2\|{\bf f}_1\|_{\infty}\leq
C\alpha^{2+\nu},  \label{srs0a}
\end{align}
and for ${\bf e}$ and $\vec{\delta}$ we will assume that
\begin{align}\begin{aligned}
\|{\bf e}''\|_{\theta_0\alpha}+\|{\bf
e}'\|_{\theta_0\alpha}+\|{\bf e}\|_{\theta_0\alpha}&\leq
C\alpha^{2},
\\
|\vec{\delta}|&\leq C\alpha^{1+\nu}, \qquad \mbox{some}\ \nu>0.
\end{aligned}
\label{srs0b}
\end{align}
Let $\nu$ be a fixed small number to be specified later and let
\begin{align*} {\mathcal Y} =\left\{({\bf f}_1, {\bf e},
\vec{\delta})\left|\begin{aligned}\|({\bf f}_1, {\bf e},
\vec{\delta})\|_{\mathcal Y}&=\alpha^{-2}(\|{\bf
f}_1''\|_{\theta_0\alpha}+\alpha\|{\bf
f}_1'\|_{\theta_0\alpha}+\alpha^2\|{\bf f}_1\|_{\infty})\\
&\quad+\alpha^{-2+\nu}(\|{\bf e}''\|_{\theta_0\alpha}+\|{\bf
e}'\|_{\theta_0\alpha}+\|{\bf
e}\|_{\theta_0\alpha})+\alpha^{-1}|\vec{\delta}|<\infty
\end{aligned}\right.\right\}.
\end{align*}
Let also
\begin{align*}
{B}_R=\{({\bf f}_1, {\bf e}, \vec{\delta})\mid \|({\bf f}_1, {\bf
e}, \vec{\delta})\|_{\mathcal Y}<R\}.
\end{align*}
We want to use Banach fixed point theorem to solve
(\ref{red8})--(\ref{red9}) in $B_{R\alpha^\nu}$, some $\nu>0$. It
is convenient to linearize this system around the approximate
solution first. We will denote $f_j=f_{0j}+f_{1j}$ and:
\begin{align*}
\LLL_j({\bf f}_1)&=-c_1
f_{j1}''+C_p(e^{\,-(f_{0j-1}-f_{0j})}(f_{1j}-f_{1j-1})-e^{\,-(f_{0j}-f_{0j+1})}(f_{1j+1}-f_{1j}),
\\
\MM_j({\bf f}, \vec{d}
)&=-C_p\big[(e^{\,-|f_j-f_{j-1}|}-e^{\,-|f_j-f_{j+1}|})
-(e^{\,-(f_{0j-1}-f_{0j})}(f_{1j}-f_{1j-1})\\
&\qquad\qquad+e^{\,-(f_{0j}-f_{0j+1})}(f_{1j+1}-f_{1j})\big]
\\
&\qquad-2\sqrt{\lambda_1}a_{j0}d_j(f_{j}'-\beta_j\eta)\sin(\sqrt{\lambda_1}a_{j0}z),\\
\KKK_j({\bf e})&=(e_j''+\lambda_1 e_j)d_0,\\
B_{0j}&= [A_{0j,11}({\bf f}_0)\int_\R Z''(s)Z(s)\,ds -\lambda_1
 A_{0j,22}({\bf
f}_0)\int_\R Z^2(s)\,ds]\cos(\sqrt{\lambda_1} a_{i0} z),\\
{\tilde h}_{0j}({\bf f}_0+{\tilde f}, \vec{d})&=-h_{0j}({\bf
f}_0+{\bf\tilde f}, \vec{d})+d_j B_{0j}.
\end{align*}
Let ${\bf \tilde f}$, ${\bf \tilde e}$ and $\vec{d}\in \R^k$
satisfying (\ref{srs0a})--(\ref{srs0b}) be fixed. We will consider
the linearized system (\ref{red8})--(\ref{red9}) written in  the
form:
\begin{align}
\LLL_j({\bf f}_1)&=\MM_j({\bf\tilde f},
\vec{d})+\widetilde\Pi_{f_{j}}({\bf f}_0+{\bf\tilde f}, {\bf\tilde
e}, \vec{d})+Q_j({\bf \tilde f}, {\bf\tilde  e}, \vec{d}),
\label{srs0c}
\\
\KKK_j({\bf e})+\delta_jB_{0j}&={\tilde h}_{0j}({\bf
f}_0+{\bf\tilde f}, \vec{d})+\widetilde\Pi_{e_j}({\bf
f}_0+{\bf\tilde f}, {\bf\tilde e}, \vec{d})+P_j({\bf \tilde f},
{\bf\tilde  e}, \vec{d}). \label{srs0d}
\end{align}

Our first goal is to show that given ${\bf \tilde f}$, ${\bf\tilde
e}$ and $\vec{d}$ one can chose $\vec{\delta}$  so that the
solvability condition for the equation (\ref{srs0d}) is satisfied.
Thus we need to compute projections of  (\ref{srs0d}) onto
$\cos(\sqrt{\lambda_1}z)$ and show that adjusting $\vec{\delta}$
we can accomplish that   for $j=1,\dots,k$:
\begin{align}
\int_0^\infty B_{0j}\cos(\sqrt{\lambda_1}z)\,dz\neq 0. \label{srs0}
\end{align}

We first notice that given function $q(z)$ such that
$$
\alpha^{-1}\|q'\|_{\theta_0\alpha} +\|q\|_{\theta_0\alpha}<\infty,
$$
we have
\begin{align}
\begin{aligned}
\int_0^\infty\cos(\sqrt{\lambda_1}z)\cos(\sqrt{\lambda_1} a_{j0}
z) q(z)\,dz&=\frac{1}{2}\int_0^\infty q(z)\,dz+\frac{1}{2}\int_0^\infty \cos(2\sqrt{\lambda_1}z)q(z)\,dz\\
&+ \int_0^\infty
z\cos(\sqrt{\lambda_1}z)\sin(\sqrt{\lambda_1}z)O(\beta_j^2\eta^2)
q(z)\,dz\\
&=\frac{1}{2}\int_0^\infty
q(z)\,dz+O(\alpha^{-1}\|q'\|_{\theta_0\alpha})+O(\|q\|_{\theta_0\alpha}).
\end{aligned}\label{srs2}
\end{align}
The first term above which is of order
$O(\alpha^{-1}\|q\|_{\theta_0\alpha})$ is expected to be the
leading order term. Thus to find asymptotic values of the
integrals involved in  (\ref{srs0}) we need to consider
\begin{align*}
I_0&=\int_0^\infty A_{0j,11}({\bf f}_0)\,dz\int_\R
Z''(s)Z(s)\,ds,\\
II_0&=-\lambda_1 \int_0^\infty A_{0j,22}({\bf f}_0)\,dz\int_\R
Z^2(s)\,ds.
\end{align*}
Let us recall that
\begin{align*}
A_{0j,11}({\bf f}_0)&=a_{j1}^2\big(
(f'_{0j})^2-\beta_j^2\eta^2\big),\\
A_{0j,22}({\bf
f}_0)&=2\beta_j\eta(\beta_j\eta-f_{0j}')+\big((f'_{0j})^2-\beta_j^2\eta^2)+2\alpha\beta_j^2
z \eta'\eta.
\end{align*}
We will first prove the following
\begin{lemma}\label{lema delta}
Let  ${\bf f}_0$ be the solution of the Toda system described above. For any $B>0$ there exists a cut-off function
$\eta(t)=0$, $t<T_1$, $\eta(t)=0$, $t>T_2$ such that
\begin{align}
|I_0+II_0|\geq B\alpha.
\label{srs3}
\end{align}
\end{lemma}
\proof{}
Let us denote $b_1=-\int_0^\infty (Z')^2$, $b_2=-\lambda_1\int_0^\infty Z^2$.
Since $a_{j1}^2=1+O(\beta_j^2\eta^2)$, therefore
\begin{align*}
I_0&=b_1\int_0^\infty\big(
(f'_{0j})^2-\beta_j^2\eta^2\big)\,dz+O(\alpha^2)\\
&=b_1\int_0^\infty\big(
(f'_{0j})^2-\beta_j^2\big)\,dz+b_1\beta^2_j\alpha^{-1}\int_0^\infty\big(1-\eta^2(t)\big)\,dt+O(\alpha^2)\\
&=I_{01}+I_{02}+O(\alpha^2).
\end{align*}
Likewise we get
\begin{align*}
II_0&=b_2\int_0^\infty \big[\big(
(f'_{0j})^2-\beta_j^2\big)+2\beta_j\eta(\beta_j-f'_j)\big]\,dz\\
&\qquad+b_2\beta_j\int_0^\infty\big[(1-\eta^2(t))(\beta_j\alpha^{-1}-1)\big]\,dt+2\beta_j\alpha^{-1}\int_0^\infty\eta(t)
(\eta(t)-1)\,dt\\
&=II_{01}+II_{02}+II_{03}.
\end{align*}
We notice that there exists a constant $C_0>0$, independent of  the cut off function $\eta$ such that
\begin{align}
|I_{01}|+|II_{01}|\leq C_0\alpha.
\label{srs4}
\end{align}
To estimate the rest of the integrals we will consider
$$
b_3=b_1\beta_j^2\alpha^2+\beta_2\beta_j\alpha^{-1}.
$$
Notice that $b_3=O(1)$ as $\alpha\to 0$. Let us assume first that
$b_3\neq 0$. Then we can take $T_2=T_1+R$ (see the definition of
$\eta$) with $R$ fixed. Increasing $T_1$ if necessary we can make
then $|I_{02}+II_{02}|$ as large as we wish while $|II_{03}|$ will
remain bounded so that
$$
|I_0+II_0|\geq
|I_{02}+II_{02}|-|II_{03}|-|I_{01}|-|II_{01}|+O(\alpha^2)\geq
B\alpha.
$$
If on the other hand $b_3=0$ then by taking $T_2-T_1=R$  large we
can make $II_{03}$ as large as we want and then
$$
|I_0+II_0|\geq |II_{03}|-|I_{01}|-|II_{01}|+O(\alpha^2)\geq
B\alpha.
$$
The proof of the Lemma is complete. \qed

We are ready now to set up the fixed point argument. We  will
begin by  outlining  the steps that are needed to complete the
proof.
\begin{enumerate}
\item
Given $({\bf \tilde f},{\bf \tilde e}, \vec{d})\in {\mathcal Y}$
we will chose $\delta$ so that the solvability condition for
(\ref{srs0d}) is satisfied.
\item
Using Step 1 we define  an operator ${\mathcal F}:
{B}_{R\alpha^{\nu}}\to {\mathcal Y}$ by
\begin{align*}
{\mathcal F}({\bf \tilde f},{\bf \tilde e}, \vec{d})=({\bf f}_1,
{\bf e}, \vec{\delta}), \quad\mbox{where}\ ({\bf f}_1, {\bf e},
\vec{\delta})\ \mbox{is a solution of
(\ref{srs0c})--(\ref{srs0d})}.
\end{align*}
\item
We will show that there exist $\nu>0$ and $R>0$ such that
\begin{align}
\|{\mathcal F}({\bf \tilde f},{\bf \tilde e}, \vec{d})\|_{\mathcal
Y}\leq R\alpha^{\nu},\label{srs5}
\end{align}
given that
\begin{align}
\|({\bf \tilde f},{\bf \tilde e}, \vec{d})\|_{\mathcal Y}\leq
R\alpha^{\nu}. \label{srs6}
\end{align}
\item
We will show that the operator  ${\mathcal F}$ restricted to
$B_{R\alpha^{\nu}}$ is a contraction.
\end{enumerate}
Once steps (1)--(4) are executed will conclude by applying Banach
fixed point theorem.

Step 1 is an easy consequence of Lemma \ref{lema delta} and give
the following formula
\begin{align}
\begin{aligned}
\delta_j\int_0^\infty
B_{0j}\cos(\sqrt{\lambda_1}z)\,dz&=\int_0^\infty\big[{\tilde
h}_{0j}({\bf f}_0+{\bf\tilde f}, \vec{d})+\widetilde\Pi_{e_j}({\bf
f}_0+{\bf\tilde f}, {\bf\tilde e}, \vec{d})+P_j({\bf \tilde f},
{\bf\tilde  e}, \vec{d})\big]\\
&\qquad\qquad\times\cos(\sqrt{\lambda_1}z)\,dz.\end{aligned}
\label{srs7}
\end{align}
We will need a precise estimate for $\vec{\delta}$ and this is the
computation we will carry on now. We will make use of the
following:
\begin{lemma}\label{lem rho}
Let $q:\R_+\to \R$ be a function such that
$\|q'\|_{\theta_0\alpha}<\infty$. Then for $s\in\N$ we have
\begin{align}
\begin{aligned}
\left|\int_0^\infty q (z)z^s \cos (\sqrt{\la_1} z)\,dz\right|&\leq
\frac{C}{\alpha^{s+1}}\|q'\|_{\theta_0\alpha},\qquad s\geq 0,\\
\left|\int_0^\infty q (z) z^s\sin (\sqrt{\la_1} z)\,dz\right|&\leq
\frac{C}{\alpha^{s+1}}\|q'\|_{\theta_0\alpha},\qquad s\geq 1.
\end{aligned}
\label{est rho}
\end{align}
Similar estimates  hold when $\|q^{''}\|_{\theta_0\alpha}<\infty$,
and and $q'(0)=0$, with $\|q'\|_{\theta_0\alpha}$ above replaced
by $\|q''\|_{\theta_0\alpha}$.
\end{lemma}
\proof{} We have
\begin{align*}
\int_0^\infty q (z) \cos (\sqrt{\la_1}
z)\,dz&=-\frac{1}{\sqrt{\lambda_1}}\int_0^\infty q
'(z)\sin(\sqrt{\lambda_1}
z)\,dz\\
&=-\frac{1}{\lambda_1}\int_0^\infty q'' (z) \cos (\sqrt{\la_1}
z)\,dz, \qquad\mbox{using}\ q'(0)=0,
\end{align*}
hence we get
$$
\left|\int_0^\infty q (z) \cos (\sqrt{\la_1}
z)\,dz\right|\leq\frac{C}{\alpha}\|q'\|_{\theta_0\alpha},
\qquad\left|\int_0^\infty q (z) \cos (\sqrt{\la_1}
z)\,dz\right|\leq\frac{C}{\alpha}\|q''\|_{\theta_0\alpha}.
$$
if $\|q''\|_{\theta_0\alpha}<\infty$.
The first inequality in (\ref{est rho}) follows now by induction.
Similarly we get
\begin{align*}
\int_0^\infty z q (z) \sin (\sqrt{\la_1}
z)\,dz&=\frac{1}{\sqrt{\lambda_1}}\int_0^\infty\big(z q
(z)\big)'\cos(\sqrt{\lambda_1}
z)\,dz-\frac{1}{\lambda_1}\int_0^\infty q' (z) \sin
(\sqrt{\la_1} z)\,dz\\
&=\Big(\frac{2}{\lambda_1^{3/2}}-\frac{1}{\lambda_1}\Big)\int_0^\infty
z q''(z)\cos(\sqrt{\lambda_1} z)\,dz
\\ &\qquad-\frac{2}{\lambda_1^2}\int_0^\infty q
''(z)\cos(\sqrt{\lambda_1}z)\,dz,\\
&
 \qquad\mbox{using}\ q'(0)=0.
\end{align*}
From these identities, using the above and induction, one can show
the second estimate in (\ref{est rho}). \qed

With some abuse of notation we will write
$$
\|({\bf\tilde f}, {\bf 0}, {\bf 0})\|_{\mathcal Y}=\|{\bf \tilde
f}\|_{\mathcal Y}, \qquad \|({\bf 0}, {\bf \tilde e}, {\bf
0})\|_{\mathcal Y}=\|{\bf \tilde e}\|_{\mathcal Y},\qquad
\mbox{etc.}
$$

Will will now compute the integrals on the right hand side of
(\ref{srs7}). We begin with
\begin{align*}
&\int_0^\infty\big[{\tilde h}_{0j}({\bf f}_0+{\bf\tilde
f})\big]\cos (\sqrt{\la_1} z)\,dz\\
&=\int_\R w''(s)
Z(s)\,ds\int_0^\infty
A_{0j,11}a_{j0}\cos (\sqrt{\la_1} z)\,dz\\
&\qquad -2d_j\int_\R\big(Z'(s)\big)^2\,ds\int_0^\infty
A_{0j,12}a_{j0}\sqrt{\lambda_1}\sin(\sqrt{\lambda_1} a_{j0} z)\cos
(\sqrt{\la_1} z)\,dz\\
&\qquad + d_j^2\int_0^\infty O\big(\|f_{0j}'+\tilde f_{j}'
-\beta_j\eta\|_{\theta_0\alpha}+\alpha^2\big)e^{\,-\theta_0\alpha|z|}\cos
(\sqrt{\la_1} z)\,dz
\\&= \Big(\int_\R w''(s) Z(s)\,ds\Big)I_j-2d_j\Big(\int_\R\big(Z'(s)\big)^2\,ds\Big) II_j+d^2_jIII_j.
\end{align*}
Since ${\bf f}_0$ is a smooth function, therefore from
(\ref{srs-05}) we get  using $f_{j0}'(0)=0$ and Lemma \ref{lem
rho}:
\begin{align*}
I_j&=\int_0^\infty a_{j0} \big((f_{0j}'+\tilde
f_{j}')^2-\beta^2_j\eta^2\big)\cos
(\sqrt{\la_1} z)\,dz\\
&=\int_0^\infty
 a_{j0}\big((f'_{0j})^2-\beta^2_j\eta^2)\big)\cos(\sqrt{\la_1}
z)\,dz\\
&\qquad +\int_0^\infty a_{j0}
\big(2f_{j0}'\tilde f_{j}'+(\tilde f_{j}')^2\big)\cos(\sqrt{\la_1})\,dz,\\
\end{align*}
hence
\begin{align*}
|I_j|\leq C(\alpha^3+\alpha^2\|{\bf \tilde f}\|_{\mathcal Y}).
\end{align*}
Similarly we get
\begin{align*}
|II_j|&\leq C\alpha\|d\|_{\mathcal Y}(\alpha^2+\alpha\|{\bf \tilde
f}\|_{\mathcal Y}),\\
|III_j|&\leq C\alpha^2\|d\|^2_{\mathcal Y}(1+\alpha\|{\bf \tilde
f}\|_{\mathcal Y}).
\end{align*}
Some tedious but standard calculations yield:
\begin{align*}
\left|\int_0^\infty\widetilde\Pi_{e_j}({\bf f}_0+{\bf\tilde f},
{\bf\tilde e}, \vec{d})\cos(\sqrt{\lambda_1}z)\,dz\right|& \leq
C\big[\alpha^3+\alpha\|d\|_{\mathcal
Y}(\alpha^2+\alpha\|{\bf\tilde f}\|_{\mathcal
Y}+\alpha^{1-\nu}\|{\bf \tilde e}\|_{\mathcal
Y})\\
&\qquad+\alpha^{2+\frac{1}{2}-\nu}\|{\bf \tilde e}\|_{\mathcal Y}\big].
\end{align*}
Using (\ref{supprho1}) and  $Z(s)\approx e^{\,-\frac{p+1}{2}|s|}$,
$p\geq 2$ we can estimate
\begin{align*}
\left|\int_0^\infty P_j({\bf \tilde f}, {\bf\tilde  e},
\vec{d})\cos(\sqrt{\lambda_1}z)\,dz \right|\leq
C\alpha^{2+\frac{1}{4}-2\sigma'}.
\end{align*}
Since $\sigma'<2^{-10}$ the above estimates and Lemma \ref{lema
delta} yield whenever $\nu<\frac{1}{4}$:
\begin{align}
\|\vec{\delta}\|_{\mathcal Y}\leq
\frac{C}{B}\big[\alpha^{\frac{1}{8}}+R^2\alpha^{\nu}\big],
\label{srs7a}
\end{align}
as long as (\ref{srs5}) is satisfied.

Step 2 follows directly from the results of Section \ref{sec
solvability}. We will now show that the claim of Step 3. From
(\ref{ll11}) we have
\begin{align}
e_j=-\delta_jS[B_{0j}]+S[{\tilde h}_{0j}({\bf f}_0+{\bf\tilde f},
\vec{d})]+S[\widetilde\Pi_{e_j}({\bf f}_0+{\bf\tilde f},
{\bf\tilde e}, \vec{d})]+S[P_j({\bf \tilde f}, {\bf\tilde  e},
\vec{d})], \label{srs8}
\end{align}
where $S$ is the operator defined in (\ref{ll11}). We have, using
estimates (\ref{e1})--(\ref{e2}):
\begin{align*}
\|\delta_jS[B_{0j}]\|_{\theta_0\alpha}\leq
C\alpha^2\|\vec{\delta}\|_{\mathcal Y}.
\end{align*}
Likewise
\begin{align*}
\|S[{\tilde h}_{0j}({\bf f}_0+{\bf\tilde f},
\vec{d})]\|_{\theta_0\alpha}\leq R_0\alpha^2+C\alpha^2\|{\bf
\tilde f}\|_{\mathcal Y}+C\alpha^2\|\vec{d}\|_{\mathcal
Y}\|{\bf\tilde f}\|_{\mathcal Y}.
\end{align*}
where $R_0>0$ is a constant depending on ${\bf f}_0$ only. Using
the formulas for $\tilde\Pi_{e_j}$ and $P_j$ we get
\begin{align*}
\|S[\widetilde\Pi_{e_j}({\bf f}_0+{\bf\tilde f}, {\bf\tilde e},
\vec{d})]\|_{\theta_0\alpha}&\leq
C\alpha^{-1}\big[\alpha^{2-\nu}\|{\bf\tilde e}\|_{\mathcal Y
}(\alpha\|\vec{d}\|_{\mathcal Y}+\alpha^2)\big],
\\
\|S[P_j({\bf \tilde f}, {\bf\tilde  e},
\vec{d})]\|_{\theta_0\alpha} &\leq C\alpha^{2+\frac{1}{8}},
\end{align*}
These estimates can be summarized as follows:
\begin{align}
\|{\bf e}\|_{\mathcal Y}\leq
R_0\alpha^\nu+C\alpha^\nu(\alpha^{\frac{1}{8}}+\alpha^{\nu}).\label{srs7b}
\end{align}
We will now estimate the right hand side of (\ref{srs0c}).
\begin{align*}
\|\MM_j({\bf\tilde f}, \vec{d})\|_{\theta_0\alpha}\leq
R_0\alpha^2\|\vec{d}\|_{\mathcal Y}+C\alpha^2(\|{\bf\tilde
f}\|^2_{\mathcal Y}+\|\vec{d}\|_{\mathcal Y}\|{\bf\tilde
f}\|_{\mathcal Y}).\end{align*} We also get in $B_{R\alpha^\nu}$:
\begin{align*}
\|\widetilde\Pi_{f_{j}}({\bf f}_0+{\bf\tilde f}, {\bf\tilde e},
\vec{d})\|_{\theta_0\alpha}+\|Q_j({\bf \tilde f}, {\bf\tilde  e},
\vec{d})\|_{\theta_0\alpha}\leq
C\alpha^2(\alpha^{\frac{1}{8}}+\alpha^{2\nu}).
\end{align*}
The last two estimates yield
\begin{align}
\|{\bf f}_1\|_{\mathcal Y}\leq R_0\|\vec{d}\|_{\mathcal
Y}+C(\alpha^{\frac{1}{8}}+\alpha^{2\nu}). \label{srs7c}
\end{align}
Combining now (\ref{srs7a})--(\ref{srs7c}) we get
\begin{align}
\|{\mathcal F}({\bf \tilde f},{\bf \tilde e}, \vec{d})\|_{\mathcal
Y}&\leq
R_0\alpha^\nu+\frac{C(1+R_0)}{B}(\alpha^{\frac{1}{8}}+\alpha^{\nu})+C(\alpha^{\frac{1}{8}}+\alpha^{2\nu})
\\
&\leq \frac{R}{2}\alpha^{\nu}. \label{srs7d}
\end{align}
provided that $\frac{1}{8}<\nu<\frac{1}{4}$, $R$ in the definition
of $B_{R\alpha^\nu}$ is taken sufficiently large, $\alpha$ is
taken small and $B$ in the Lemma \ref{lema delta} is taken large.
This shows the claim of Step 3.

It now remain to show the claim of Step 4. This follows from a
rather straightforward application of the theory of Lipschitz
dependence of various terms involved in
(\ref{srs0c})--(\ref{srs0d}). We have for instance
\begin{align}
\|\vec{\delta}^{(1)}-\vec{\delta}^{(2)}\|_{\mathcal Y}\leq
\frac{C}{B}\|({\bf\tilde f}^{(1)}, {\bf\tilde e}^{(1)},
\vec{d}^{(1)})-({\bf\tilde f}^{(2)}, {\bf\tilde e}^{(2)},
\vec{d}^{(2)})\|_{\mathcal Y}, \label{srs8a}
\end{align}
and
\begin{align}
\|{\bf f}_1^{(1)}-{\bf f}_1^{(2)}\|_{\mathcal Y}+ \|{\bf
e}^{(1)}-{\bf e}^{(2)}\|_{\mathcal Y}\leq
C\alpha^{\nu}\|({\bf\tilde f}^{(1)}, {\bf\tilde e}^{(1)},
\vec{d}^{(1)})-({\bf\tilde f}^{(2)}, {\bf\tilde e}^{(2)},
\vec{d}^{(2)})\|_{\mathcal Y}, \label{srs9}
\end{align}
We leave the details to the reader. From these estimates, taking
$B$ larger if necessary we get that ${\mathcal F}$ is a Lipschitz
contraction as claimed. This shows that ${\mathcal F}$ has a fixed
point in $B_{R\alpha^\nu}$. The proof of the theorem is complete.
\qed

\bigskip

\bigskip
{\bf Acknowledgments:}  This work has been partly supported  by
chilean research grants Fondecyt 1070389, 1050311, FONDAP, Nucleus
Millennium grant P04-069-F, an Ecos-Conicyt contract and an
Earmarked Grant from RGC of Hong Kong. The fourth author thanks Professor Fanhua Lin for  a nice conversation on the classification of solutions to (\ref{sch}), which motivated this research.

\end{document}